\documentclass[english]{amsart}

\usepackage[margin=1in]{geometry}
\usepackage{amsmath,amsthm,amssymb}
\usepackage{wasysym}

\usepackage[usenames,dvipsnames,svgnames,table]{xcolor}

\usepackage{hyperref}
\usepackage{graphicx}

\def\eps{\varepsilon}
\def\R{\mathbb R}

\def\N{\mathbb N}

\def\d{\mathrm{d}}
\def\dv{\mathrm{d} v}
\def\dx{\mathrm{d} x}

\def\dz{\mathrm{d} z}
\def\fin{f_0}
\def\T{\mathbb{T}}
\def\Rone{R_1}
\def\rzero{r_1}

\def\cCl{\mathcal{C}_{\ell}}
\def\cClloc{\mathcal{C}_{\ell,\text{loc}}}
\def\cK{\mathcal{K}}
\def\cL{\mathcal{L}}

\def\cLdiv{\mathcal{L}}

\def\Cwhi{C_{\text{whi}}}

\def\Qext{Q_{\mathrm{ext}}}

\newcommand{\classK}{\eqref{e:ellipticity_upper}, \eqref{e:ellipticity_cone} and \eqref{e:symmetry_non-divergence}}

\DeclareMathOperator{\pv}{P.V.}
\DeclareMathOperator{\degk}{deg_{\, kin}}

\newtheorem{thm}{Theorem}[section]
\newtheorem{lemma}[thm]{Lemma}
\newtheorem{cor}[thm]{Corollary}
\newtheorem{proposition}[thm]{Proposition}

\theoremstyle{definition}
\newtheorem{defi}{Definition}[section]

\theoremstyle{remark}
\newtheorem{remark}{Remark}[section]
\newtheorem{assumption}{Assumption}[section]

\numberwithin{equation}{section}

\begin{document}

\title[Regularity for the Boltzmann equation]{Regularity for the Boltzmann equation conditional to macroscopic bounds}
\author{Cyril Imbert}
\address[C.~Imbert]{CNRS \& Department of Mathematics and Applications, \'Ecole Normale Sup\'erieure (Paris) \\
45 rue d'Ulm, 75005 Paris, France}
%\address[C.~Imbert]{CNRS \& Instituto de   matem\'atica pura e aplicada, Estrada Dona Castorina, 110, Jardim   Bot\^anico, CEP 22460-320, Rio de Janeiro, RJ, Brasil}
\email{Cyril.Imbert@ens.fr}
\author{Luis Silvestre}
\thanks{LS is supported by NSF grant DMS-1764285.}
\address[L.~Silvestre]{Mathematics Department, University of Chicago,
  Chicago, Illinois 60637, USA} \email{luis@math.uchicago.edu}

\date{\today}

\begin{abstract}
  The Boltzmann equation is a nonlinear partial differential equation
  that plays a central role in statistical mechanics. From the
  mathematical point of view, the existence of global smooth solutions
  for arbitrary initial data is an outstanding open problem. In the
  present article, we review a program focused on the type of particle interactions known as \emph{non-cutoff}. It is
  dedicated to the derivation of a priori estimates in $C^\infty$,
  depending only on physically meaningful conditions. We prove that
  the solution will stay uniformly smooth provided that its mass,
  energy and entropy densities remain bounded, and away from vacuum.
\end{abstract}

\maketitle

\setcounter{tocdepth}{1}
\tableofcontents

\section{Introduction}

% \subsection{The Boltzmann equation}

In this work, we study a priori estimates on the regularity of
solutions to the Boltzmann equation,
\begin{equation}
  \label{e:boltzmann}
  \partial_t f + v \cdot \nabla_x f = Q(f,f).
\end{equation}

The Boltzmann equation describes the dynamics of a dilute gas. The nonnegative function
$f : [0,\infty) \times \R^d \times \R^d \to [0,\infty)$ encodes the density of particles which at time $t$ have position $x$ and velocity $v$. The left hand side of the equation, pure transport, denotes the fact that particles travel in straight lines in the absence of external forces. The right-hand side, the Boltzmann collision operator $Q(f,f)$, corresponds to the fluctuations in velocity that result from particle interactions. Our regularity estimates apply to a type of
collision operators known as \emph{non-cutoff}, \textit{i.e.} where long range grazing collisions play a crucial role.

The kinetic description of gas dynamics is intermediate between the microscopic scale, associated to the trajectory of each individual particle, and the macroscopic scale associated to fluid mechanics, such as the Euler or Navier-Stokes equations. In kinetic equations, the density function $f$ encodes the state of a gas, the position and velocity of its particles for each time, but only in a statistical way.

The equation was initially studied by Maxwell and Boltzmann by the end of the 19\textsuperscript{th} century. As one of the fundamental equations from statistical mechanics, it has received substantial interest from the mathematical community through the years. There are varied mathematical problems that result from the study of the Boltzmann equation. An important line of research is on the rigorous derivation of the equations of fluid mechanics from it. It is also natural to study the well posedness of the initial value problem, and it involves different problematics. In this program, we are concerned with regularity estimates for the solutions, which are intimately related to the classical solvability of the equation. As we discuss below, finding smooth solutions unconditionally, for any initial data, appears to be completely out of reach with current techniques. Our results are on a priori regularity estimates conditional to certain macroscopic bounds on physically meaningful quantities.

\subsection{The collision operator}

The Boltzmann collision operator has two terms: a gain term and a loss
term. The loss term represents the fact that a particle with velocity
$v$ may collide with another one and change its velocity to something
else. Conversely, the gain term accounts for collisions of other
particles, with different velocities, that result in new particles
with velocity $v$.

The Boltzmann collision operator is an integro-differential operator
acting on the function $f(t,x,\cdot)$ for fixed values of $t$ and
$x$. Its most common form is
\begin{equation} \label{e:collision_operator}
 Q(f,f) = \iint_{\R^d \times \partial B_1} (f(v'_*)f(v')
  -f(v_\ast)f(v)) B(|v_*-v|,\cos \theta) \d \sigma \d v_\ast,
\end{equation}
where $\cos \theta = \frac{v-v_*}{|v-v_*|} \cdot \sigma$. The
nonnegative function $B$ is referred to as the \emph{collision
  kernel}.

The velocities $v$, $v_\ast$, $v'$ and $v_\ast'$ are pre- and
post-collisional velocities. Particles with velocities $v$ and
$v_\ast$ collide and immediately switch their velocities to $v'$ and
$v_\ast'$ with a rate given by $B$. Conversely, particles with
velocities $v'$ and $v_\ast'$ may collide and turn their velocities to
$v$ and $v_\ast$ at the same rate. The kernel $B$ denotes this rate by
which velocities $v$ and $v_\ast$ deviate to $v'$ and $v'_\ast$ after
a collision.

We consider elastic collisions that conserve momentum and energy. It
is represented by the following two relations between pre- and
post-collisional velocities.
\begin{equation}
  \label{e:ellastic_collisions}
  \begin{aligned}
    v+v_\ast &= v' + v_\ast', \\
    |v|^2 + |v_\ast|^2 &= |v'|^2 + |v_\ast'|^2.
  \end{aligned}
\end{equation}

Consequently, the segments $v \, v_\ast$ and $v' \, v_\ast'$ are two
diameters of the same sphere in $\R^d$. For any given $v$, we can
parametrize all possible values of $v_\ast$, $v'$ and $v_\ast'$
through the value of $v_\ast \in \R^d$ and a unit vector $\sigma$
denoting the direction of $v'$ from the center of the sphere,
\begin{align*}
v' := \frac{v+v_\ast}2 + \frac{|v-v_\ast|}2 \sigma, \\
v_\ast' := \frac{v+v_\ast}2 - \frac{|v-v_\ast|}2 \sigma.
\end{align*}
This explains the parametrization of the integral in
\eqref{e:collision_operator}.

Several different collision kernels $B$ can be
considered. A model where particles bounce each other like billiard balls leads to $B = c|v-v_\ast|$ for some constant $c>0$. A
model where particles repel each other by a power law potential
when they are sufficiently close leads to a collision kernel $B$ that is singular
around $\theta=0$. Moreover, in that case, $B$ is not integrable with
respect to $\sigma$. The integral in \eqref{e:collision_operator}
still makes sense thanks to the cancellation in the factor
$(f(v_\ast') f(v') - f(v_\ast) f(v))$ as $v' \to v$. Note that the
gain and loss terms, both integrate to $+\infty$. The non-integrability
of $B$ might be considered a \emph{difficulty} in the analysis of the
equation. It is common to tame this singularity by considering only
collision kernels $B$ that are integrable. This
integrability condition for $B$ takes the name \emph{Grad's cutoff
  assumption}.

For the purpose of this work, it is essential that we \textbf{do not}
make Grad's cutoff assumption. The singularity in $B$ around
$\theta=0$ is what drives the regularization effects in the equation
that we exploit to obtain our a priori estimates. We focus on the
standard family of non-cutoff collision kernels, of the form
\begin{equation}\label{e:B}
B(r,\cos \theta) = r^\gamma b(\cos \theta) \quad \text{ with} \quad b (\cos \theta)
\approx|\sin (\theta/2)|^{-(d-1)-2s}
\end{equation}
with $\gamma > -d$ and $s \in (0,1)$.

Our regularity results are restricted to the range of parameters
$\gamma +2s \in [0,2]$. This includes cases usually referred to as
hard potentials ($\gamma \ge 0$) and moderately soft potentials
($\gamma \in [-2s,0)$). See Section \ref{s:very_soft} for a discussion
of open problems outside this range.

\subsection{Conserved quantities and entropy}

The dynamics of the Boltzmann equation conserve the total mass, energy
and momentum of the solution. That is
\begin{align*}
\langle \text{mass} \rangle &= \iint_{\R^d \times \R^d} f(t,x,v) \, \d v \d  x = \iint_{\R^d \times \R^d} f_0(x,v) \, \d v \d x, \\
\langle \text{momentum} \rangle &= \iint_{\R^d \times \R^d} f(t,x,v) v \, \d v \d x = \iint_{\R^d \times \R^d} f_0(x,v) v \, \d v \d x, \\
\langle \text{energy} \rangle &= \iint_{\R^d \times \R^d} f(t,x,v) |v|^2 \, \d v \d x = \iint_{\R^d \times \R^d} f_0(x,v) |v|^2 \, \d v \d x
\end{align*}
where $\fin (x,v) = f(0,x,v)$.

In addition to these three conserved quantities, a remarkable property
of Boltzmann's collision operator is the fact that the total \emph{entropy} of
solutions decreases with time,
\begin{equation}
  \label{e:entropy}
\langle \text{entropy} \rangle = \iint_{\R^d \times \R^d} f \log f (t,x,v)  \d v \d x \le \iint_{\R^d \times \R^d} f_0 \log f_0 (x,v)  \d v \d x.
\end{equation}

\subsection{The hydrodynamic limit}

We define the mass, momentum, energy and entropy densities as the functions corresponding to the conserved quantities described above but without integrating in the space variable $x$. That is
\begin{align*}
\rho(t,x) &:= \int_{\R^d} f(t,x,v) \d v, \\
\rho(t,x) u(t,x) &:= \int_{\R^d} f(t,x,v) v \d v, \\
e(t,x) &:= \int_{\R^d} f(t,x,v) |v|^2 \d v, \\
h(t,x) &:= \int_{\R^d} f \log f(t,x,v) \d v.
\end{align*}
The temperature density $\theta(t,x)$ is usually defined as
\[ \theta(t,x) := \frac 13 \left( \frac e \rho - |u|^2 \right) = \frac 1 {3 \rho} \int_{\R^d} f(t,x,v) |v - u|^2 \d v.\]

The values of $\rho$, $u$, $e$, $h$ and $\theta$ are the
macroscopic quantities corresponding to a solution $f$ to the
Boltzmann equation. In certain asymptotic regimes, they formally
converge to solutions to classical hydrodynamic equations like the Euler
and Navier-Stokes (see \cite{bardos1991}). More precisely, for a small
parameter $\eps>0$, we consider the Boltzmann equation with enhanced
collisions,
\begin{equation} \label{e:eps-Boltzmann}
 \partial_t f + v \cdot \nabla_x f = \frac 1 \eps Q(f,f)
\end{equation}
and a family of solutions $f^\eps$ to the equation
\eqref{e:eps-Boltzmann} with the same initial data. As $\eps \to 0$,
the solutions are expected to converge to Maxwellian functions in $v$
whose corresponding hydrodynamic quantities satisfy the compressible
Euler equations
\begin{equation} \label{e:compressible_Euler}
\begin{aligned}
\partial_t \rho + \nabla \cdot ( \rho u ) &= 0, \\
\partial_t(\rho u) + \nabla \cdot (\rho u \otimes u) + \nabla (\rho \theta) &= 0, \\
\partial_t(\rho (|u|^2 + 3\theta)) + \nabla \cdot \left( \rho u \left( |u|^2 + 5 \theta \right) \right) &= 0.
\end{aligned}
\end{equation}

Moreover, the limit of the entropy density satisfies $h = \rho \log(\rho^{2/3} / \theta)$ and
\[ \partial_t h + \nabla \cdot \left( uh \right) \leq 0.\]

For $\eps>0$ small (but not quite zero), the hydrodynamic quantities
solve certain local conservation laws. Up to a first order
approximation, they solve a compressible Navier-Stokes system with a
viscosity term depending on $\eps$.

The kinetic function $f$ provides a more detailed description of the state of the fluid than the hydrodynamic quantities alone. Yet, the hydrodynamic quantities are all that might be observable from a macroscopic point of view. Because the compressible Euler and Navier-Stokes equations arise as asymptotic regimes, one may wonder if any mechanism that produces a singularity in the flow of the hydrodynamic equations would be reproduced in the kinetic framework of the Boltzmann equation as well.

It is well known that, in the compressible Euler equations, shock
singularities can emerge from the flow even if the initial data is
smooth and away from vacuum. A \emph{shock} singularity is a
discontinuity in the solution of the equations, in a way that the
values of all functions involved in the system stay bounded uniformly
up to the singularity. The existence of shock singularities is
classical and easy to verify in one dimensional profiles. The density
also stays away from vacuum. In the three dimensional setting, a well
known result by Sideris \cite{sideris1985} shows the development of
singularities for the compressible Euler equations. Even though his
result does not describe the precise type of singularity that emerges,
the proof suggests that it may be a shock. In recent years, the
emergence of shock singularities for the compressible Euler equations
and their stability has been studied in different scenarios (see
\cite{christodoulou2007,Luk2018,buckmaster2019formation,buckmaster2019formation2}). We will not see singularities in kinetic equations as a result of shocks, because the function $f$ allows for different velocities to coexist at the same point in space. Likewise, shock singularities will not be observed in the Navier-Stokes equations since
the diffusion would smooth out any bounded solution that stays away from vacuum.

There are other types of singularities for hydrodynamic equations where at least one of the values $\rho(t,x)$, $u(t,x)$ or $\theta(t,x)$ becomes unbounded in some finite time. Another possible singular behavior would be when $\rho$ flows to zero (creating vacuum) or the temperature $\theta$ becomes zero somewhere (corresponding to unbounded entropy). There is a couple of very recent papers, \cite{merle2019smooth} and \cite{merle2019implosion}, where the authors construct a stable \emph{implosion} singularity for the compressible Euler and Navier-Stokes equations. Recall that the compressible Euler and Navier-Stokes equations correspond to the asymptotic behavior of \eqref{e:eps-Boltzmann} for small $\eps$. Unlike the case of shock singularities, there is no apparent incompatibility between implosion singularities and the flow of the Boltzmann equation. It is conceivable that there may exist some solution to the Boltzmann equation that flows into a singularity, whose associated hydrodynamic quantities resemble the singular solution in \cite{merle2019smooth} near the implosion point. It is currently difficult to assert whether this type of singular solutions exist. Either a rigorous construction or an impossibility proof seems complicated, requires further work, and probably new ideas.

The study of hydrodynamic equations presents great mathematical difficulties. The Boltzmann equation keeps track of more information and is more complex than the compressible Euler or Navier-Stokes equations. It is to be expected that the study of the Boltzmann equation, its well posedness, regularity of solutions, and possible singularities, will include, at least, similar difficulties. Given a smooth initial data $f_0$, away from vacuum, and with appropriate decay as $|v| \to \infty$, the question of existence of a smooth solution of the Boltzmann equation
\eqref{e:boltzmann} is an outstanding open problem for any choice of a physical collision operator $Q$. Given our current understanding of the related hydrodynamic models, and the extra complexity of the Boltzmann equation, the question of well posedness \footnote{Here, we mean unconditional well posedness far from equilibrium.} for \eqref{e:boltzmann} appears to be completely out of reach.

\subsection{The conditional regularity program}

Given the insurmountable difficulties of studying the existence of
smooth solutions to \eqref{e:boltzmann} in general, it makes sense to
study their \emph{conditional regularity}. In the program reviewed in
this note, we assume point-wise bounds on the hydrodynamic quantities:
mass density, energy density and entropy density, described in the
previous subsection. Conditional to these bounds, we prove $C^\infty$
estimates for the solution $f$ to the Boltzmann equation
\eqref{e:boltzmann}. By assuming that the hydrodynamic quantities stay
bounded, we remove some of the difficulties of hydrodynamic equations,
and we focus on the new regularization effects that are characteristic
of the kinetic scale.

More precisely, we make the following assumption: there exist positive
constants $m_0,M_0,E_0,H_0$ such that for all $(t,x)$,
\begin{gather} \tag{H}
  \label{e:hydro}
  \left\{
    \begin{aligned}
      m_0 \le \int_{\R^d} f(t,x,v) \dv &\le M_0, \\
      \int_{\R^d} f(t,x,v) |v|^2 \dv &\le E_0, \\
      \int_{\R^d} f \log f (t,x,v) \dv & \le H_0.
    \end{aligned}
\right.
\end{gather}

The lower bound $m_0$ on the mass density prevents the formation of
vacuum regions. The upper bound $H_0$ on the entropy density prevents
the concentration of the function $f$ as a singular measure. In
particular, the upper bound on the entropy density, together with the
other upper bounds, prevents the temperature to reach absolute zero.

We stress that a justification that the assumption \eqref{e:hydro}
holds for general solutions of \eqref{e:boltzmann} seems to be far out
of reach at the moment. If an implosion singularity with a similar
structure as in \cite{merle2019implosion} was possible for the
Boltzmann equation, then the assumption \eqref{e:hydro} may fail for
such solutions.

We will see that as long as \eqref{e:hydro} is true, solutions of the
Boltzmann equation \emph{remain smooth}. It is the reason why we claim
that a singularity is always macroscopically observable. We now state
our main theorem, which we proved in \cite{imbert2019regularity} and
whose proof relies on the series of articles
\cite{imbert2020weak,imbert2020decay,imbert2020schauder}.
%-------------------------------------------
\begin{thm}[Global regularity estimates]
  \label{t:main}
  Let $f$ be a solution to the Boltzmann equation
  \eqref{e:boltzmann}. Assume that $f$ is periodic in space, the
  collision kernel is of the form \eqref{e:B} and
  $\gamma+2s \in [0,2]$. If \eqref{e:hydro} holds, then for any
  multi-index $k \in \mathbb N^{1+2d}$, any time $\tau>0$ and any decay rate $q>0$,
\begin{equation} \label{e:main_estimate}
 \| (1+|v|)^q D^k f \|_{L^\infty([\tau,T) \times \R^d \times \R^d)} \leq C_{k,q,\tau}.
\end{equation}
Here $D^k$ is an arbitrary derivative of $f$ of any order, in $t$, $x$ and/or $v$.

When $\gamma>0$, the constants $C_{k,q,\tau}$ depend only on $k$, $q$
and $\tau$, and the constants $m_0$, $M_0$, $E_0$ and $H_0$ from
\eqref{e:hydro}, and the parameters $s$, $\gamma$ from the assumption
\eqref{e:B} on the collision kernel, and dimension $d$.

When $\gamma\leq 0$, the constants $C_{k,q,\tau}$ depend in addition
on the pointwise decay of the initial data $\fin$. That is, on the constants
$N_r$ with $r \ge 0$, given by
\begin{equation} \label{e:initial_decay_assumption}
 N_r := \sup_{x,v} (1+|v|)^r f_0(x,v) \qquad \text{for each } r\geq 0.
\end{equation}
\end{thm}
% ---------------------------------------------

Theorem \ref{t:main} provides an a priori estimate on the smoothness
and decay of the solutions of the Boltzmann equation without cutoff,
provided that \eqref{e:hydro} holds.

Note the difference between the cases $\gamma>0$ (hard potentials) and
$\gamma\leq 0$ (soft potentials). In the case of hard potentials, all
the bounds depend on the quantities in \eqref{e:hydro} and the
parameters of the equation only. Both the smoothness and the decay
estimates are self-generated for positive time. In the case of soft
potentials, the equation does not force a fast decay at infinity, but
it only propagates it. Our estimates depend on the (pointwise) decay
of the initial data. Note also that in both cases our estimates remain
uniform as $t \to \infty$.

The assumption that $f$ is periodic in space is a convenient way for
us to avoid the extra difficulties of analyzing the boundary behavior
of kinetic equations in bounded domains. Our estimates do not depend
on the size of the period. It would be straight forward to reproduce
the estimates in Theorem \ref{t:main} if instead of periodicity in
$x$, we assume that $f(t,x,v)$ converges to a fixed Maxweillian as
$|x| \to \infty$.

Theorem \ref{t:main} provides a priori estimates for classical
(smooth) solutions. In order to avoid gratuitous technical
difficulties, we work with a very strong notion of solution. We start
with a function $f$ that is $C^\infty$ with respect to all variables
and, moreover, for every $q>0$, $(1+|v|)^q f(t,x,v)$ converges to zero
as $|v| \to \infty$ uniformly in $x$ and $t$. This is merely a
qualitative assumption. The estimates in Theorem~\ref{t:main} are
independent of any norm quantifying smoothness or decay for $f$. The
question of whether the estimates of Theorem~\ref{t:main} would hold
for any weaker notion of solution is discussed in Section
\ref{s:extensions}.

\subsection{Previously known regularity estimates}
\label{s:previous_results}

The first progress toward understanding the regularization effect of the Boltzmann equation without cutoff appeared in the form of entropy dissipation estimates. As we mentioned before, the total entropy is nonincreasing in time. The entropy dissipation equals minus its derivative
\[ D(t) := -\frac{\d}{\d t} \iint_{\R^d\times \R^d} f \log f \, \d v \d x.\]
In \cite{advw}, they showed that the entropy dissipation is bounded below by a weighted fractional Sobolev norm minus lower order terms.
\[ D(t) \gtrsim \|\sqrt f\|_{L^2_x H^s_v (|v| < R)}^{2} - C \|f\|_{L^1_2}^2\]
where $L^1_2$ denotes the weighted norm $L^1( \R^d \times \R^d, (1+|v|^2) \dx \dv)$.
This entropy dissipation estimate effectively provides the regularity estimate $\sqrt f \in L^2_{t,x} H^s_v$, at least for small velocities. It is the first indication of a regularization effect of the Boltzmann equation in the non-cutoff case.

In \cite{advw}, they also obtain some form of a coercivity estimate for the Boltzmann collision operator $Q(f,f)$ (under rather restrictive assumptions on $B$), and the cancellation lemma which will be used later in this paper to compute the lower order term in \eqref{e:boltzmann_rewritten}.

The coercivity estimates for the Boltzmann collision operator were subsequently improved in several papers including
\cite{AlexandreLittlewood1,AlexandreLittlewood2,MouhotExpl,AlexandreI,AlexandreSome,ChenHe2011,he2016sharp,AlexandreReview,gressman2011}.

The coercivity estimate is the main tool to obtain global smooth solutions for the non-cutoff Boltzmann equation in the \textbf{space homogeneous} regime. The first result of this kind was by Desvillettes and Wennberg \cite{desvillettes2005smoothness} for a rather restrictive family of collision kernels $B$. Later on,  He improves the result in \cite{he2012wellposedness} to include the natural non-cutoff collision kernels provided only $\gamma+2s > 0$. See also \cite{AlexandreLittlewood1, AlexandreLittlewood2, huo2008regularity, morimoto2009regularity,desvillettes2009stability,morimoto2010gevrey,alexandre2012smoothing,ChenHe2011,zhang2012gevrey} for other results in the space homogenous regime. Note that due to the conservation of mass and energy, and the monotonicity of entropy, our condition \eqref{e:hydro} always holds true in the space homogeneous setting, for any initial data with finite mass, energy and entropy.

For the space in-homogeneous Boltzmann equation without cutoff, we only know the existence of smooth solutions if the initial data is sufficiently close to a Maxwellian with respect to a suitable norm. See \cite{gs2011,AlexandreI,AlexandreII,amuxy2011ultimate}.
A sharp analysis on the asymptotic behavior as $|v| \to \infty$ for the coercivity estimate plays a key role in the proof of \cite{gs2011}.

The well posedness of the equation for short time is obtained in \cite{morimoto2009regularity} for initial data $f_0$ with Gaussian decay as $|v| \to \infty$. In \cite{morimoto2015local} and \cite{henderson2019local}, they develop the local well posedness theory for initial data $f_0$ that decays only algebraically as $|v| \to \infty$.

Regularity results for the inhomogeneous non-cutoff Boltzmann equation far from
equilibrium and beyond the coercivity estimates are very scarce. The
most relevant results in the literature are given in \cite{amuxy2010}
and \cite{ChenHe2012}. They prove $C^\infty$ regularity estimates for
any solution $f$ to the Boltzmann equation \eqref{e:boltzmann} whose
mass density is bounded below and with five derivatives (in all
directions with respect to $x$ and $v$) in a weighted $L^2$ space with
infinite moments. Our condition \eqref{e:hydro} is naturally much less
restrictive and arguably more physically meaningful.

Concerning \emph{unconditional} regularity estimates for solutions to the non-cutoff Boltzmann equation \eqref{e:boltzmann}, far from equilibrium, the only result we are aware of is by Arsenio and Masmoudi \cite{arsenio2012}, in a fractional Sobolev space with a low order of differentiability. They show that $f/(1+f) \in W_{t,x,v,loc}^{s,p}$, for every $p \in [1,d/(d-1))$ and $s>0$ depending on $p$.

\subsection{Other related results}

Here, we discuss briefly a few other results for the Boltzmann equations where the conditional regularity of solutions is applied.

There is a well known result by Desvillettes and Villani \cite{dv2005} about the long term behavior of solutions to the Boltzmann equation, in the non-cutoff case. It says that solutions converge to equilibrium as $t \to \infty$ provided that they stay uniformly smooth and bounded below by a fixed Maxwellian function. Any conditional regularity result, like originally the one in \cite{amuxy2010}, improves the result in \cite{dv2005} by reducing its regularity assumption. After our Theorem \ref{t:main}, the condition to apply the theorem of Desvillettes and Villani in \cite{dv2005} becomes simply that \eqref{e:hydro} holds. We discuss it further in Section \ref{s:extensions}.

There is a recent result by Duan, Liu, Sakamoto and Strain \cite{duan2019global} where they construct global mild solutions to the Boltzmann equation without cutoff, whose initial data is close to equilibrium in a suitable norm. Interestingly, these solutions are not a priori very smooth, but one can verify that our assumption \eqref{e:hydro} holds. Thus, we can deduce the $C^\infty$ smoothness and decay of the solutions in \cite{duan2019global} as a consequence of Theorem~\ref{t:main}. We also discuss this application in Section \ref{s:extensions}.

\subsection{Notation}
We use this section both to clarify our choices of notation, and also as a form of glossary of symbols. We list the places in which most symbols used throughout this paper are defined.

We typically use the letters $t \in \R$ for time, $x \in \R^d$ for
space and $v \in \R^d$ for velocity. We also use $z \in \R^{1+2d}$ for
$z=(t,x,v)$.

We write $a \lesssim b$ to denote that there exists a constant $C$
(depending on the parameters that are appropriate for each scope) so
that $a \leq C b$. We write $a \approx b$ to express that
$a \lesssim b$ and $a \gtrsim b$.

In the rest of this section we refer to objects that will be defined later on in this text. It makes little sense to read it linearly. It can be useful to come back here whenever the reader wants to remember the meaning of some symbol.

The integro-differential operator $\cL_K$ is defined in
\eqref{e:integro-differential_operator} with respect to some kernel
$K(t,x,v,v')$. It is a fractional order diffusion operator in the
velocity variable. It is equal to some integral involving the values
of $f(t,x,v')$ and $K(t,x,v,v')$ for fixed $t$ and $x$. The same
formula makes sense for $f=f(v)$ and $K=K(v,v')$ independent of $t$
and $x$.

The class of kinetic intego-differential equations of order $2s$ is
invariant by a special scaling $S_R$ and the Galilean Lie group
structure $(\R^{1+2d}, \circ)$ defined in Section~\ref{s:invariance}.

The H\"older norms $\cCl^\alpha$ are invariant by the scaling $S_R$
and the Galilean group structure. They are defined in Section
\ref{s:holder}.

The kinetic cylinders $Q_1 = (-1,0] \times B_1 \times B_1$ and
$Q_R(z_0) = \{ z_0 \circ S_R (z): z \in Q_1\}$ are defined in Section~\ref{s:cylinders}.

\subsection{Organization of the article}

After this introduction, in Section \ref{s:diffusion}, we analyze the
structure of the Boltzmann equation. We describe the collision
operator $Q(f,f)$ as a nonlinear integro-differential operator, and we
relate it to the study of nonlocal elliptic operators. We also compare
the Boltzmann equation with some basic hypoelliptic equations.

In Section \ref{s:warm_up}, we outline all the steps involved in the proof of Theorem \ref{t:main}. Section \ref{s:generic} sets up the foundation for working with kinetic integro-differential equations that will be needed later. Sections \ref{s:DG} to \ref{s:bootstrap} provide some details for each of the ingredients listed before in Section \ref{s:warm_up}. We provide references to the original papers in each case.

In Section \ref{s:extensions}, we describe two implications of Theorem \ref{t:main}: a continuation criteria and an improvement of the conditions for convergence to equilibrium. We also discuss its applicability to weak solutions. We finish the paper with a collection of related open problems in Section \ref{s:open_problems}.

\section{Structure of the diffusion}
\label{s:diffusion}

\subsection{A  \emph{hypoelliptic}  structure}
\label{s:hypoelliptic}

The study of regularity properties of solutions of the Boltzmann equation necessarily starts with the study of the structure of the diffusion.  The regularization effect of the equation is driven by the fact that the collision operator is \emph{diffusive} with respect to the velocity variable. For non-cutoff collision kernels as in
\eqref{e:B}, the collision operator turns out to be a nonlinear
integro-differential operator of order $2s$. The collision operator
produces some regularization with respect to the velocity
variable. This diffusion in ``$v$'' combined with the transport in
``$x$'' gives rise to a \emph{hypoelliptic} structure that regularizes
the solution in all variables.

In order to understand how the diffusion in velocity combined with the
transport in space produces a regularization effect in all directions,
it is better to start with a simpler toy model. If we replace the
collision operator $Q(f,f)$ with the Laplacian of $f$ with respect to
$v$, we arrive at the following equation
\begin{equation} \label{e:kolmogorov}
  \partial_t f + v \cdot \nabla_x f = \Delta_v f.
\end{equation}
In 1934, Kolmogorov computed explicitly in \cite{kolm} the fundamental solution
for this equation. This kernel is a $C^\infty$ function with rapid
decay at infinity both in $v$ and $x$. It turns out that the
\emph{Kolmogorov equation} \eqref{e:kolmogorov} enjoys similar
regularization properties as the usual heat equation.  Such an effect
is a consequence of the combination of two mechanisms. The diffusion
in the velocity variable $v$ regularizes the solution in this variable
and the \emph{free streaming} operator
$\partial_t + v \cdot \nabla_x$ transfers regularity
in $v$ into regularity in the $(t,x)$ variables. This remarkable
observation is the starting point of the hypoellipticity theory
developed by H\"ormander from 1967.

In the Boltzmann equation, the collision operator is an
integro-differential nonlinear diffusion. Its closest linear analog
would be the \emph{fractional Kolmogorov equation},
\begin{equation}
  \label{e:fractional_kolmogorov}
  \partial_t f + v \cdot \nabla_x f + (- \Delta)^s_v f = 0.
\end{equation}
The fractional Kolmogorov equation also enjoys similar regularizing
properties. Its corresponding fundamental solution is $C^\infty$, but it decays
polynomially at infinity, much like the fundamental solution corresponding to
the fractional heat equation. The Boltzmann equation is a nonlinear
variant of \eqref{e:fractional_kolmogorov}. From this point of view,
it is natural to expect regularity estimates by using tools from
hypoelliptic and integro-differential equations.

\subsection{Non-local diffusions}
\label{s:nonlocal_diffusion}

Motivated by probabilistic models involving discontinuous stochastic
processes, since the beginning of the 21\textsuperscript{st}
century there was an explosion of results in the area of nonlocal
diffusions. Basically, a linear parabolic integro-differential
equation is an equation of the form
\begin{equation} \label{e:pide}
 \partial_t f (t,v) = \int_{\R^d} [ f(t,w) - f(t,v) ] K(t,v,w) \d w.
\end{equation}
Here $K$ is a nonnegative kernel function satisfying some
nondegeneracy and symmetry assumptions. This type of equations can be
studied in the context of parabolic equations. They satisfy similar
characteristic properties: maximum principles, energy dissipation
inequalities, regularization effects, \textit{etc}. Some of the landmark results that make it possible to study the
regularity of nonlinear parabolic equations were reproduced in the
nonlocal setting. They include
\begin{itemize}
\item the Harnack inequality of De Giorgi, Nash and Moser
  \cite{komatsu1995, barlow2009non,basslevin2002,caffarelli2010drift,kassmann2009priori,Felsinger2013,chan2011,kassmann2013regularity};
\item the Krylov-Safonov Harnack inequality  \cite{basslevin2002,bass2002harnack,song2004,bass2005holder,bass2005harnack,%
silvestre2006holder,caffarelli2009regularity,silvestre2011differentiability,davila2012nonsymmetric,davila2014parabolic,kassmann2013intrinsic,bjorland2012,rang2013h,schwab2016};
\item the Schauder estimates \cite{MP,tj2015,serra2015,imbert2016schauder,tj2018}.
\end{itemize}

The diversity of results in this area is explained
in part by the richness of the family of equations. While a classical
second order diffusion is characterized simply by a positive definite
matrix of coefficients at each point, the integro-differential
diffusion is defined in terms of a whole kernel function
$K(t,v,\cdot)$, giving much more flexibility in terms of possible
structural assumptions.

The primal example of an integro-differential operator is the
fractional Laplacian $(-\Delta)^{s}$, which corresponds to the kernel
$K(t,v,w) = c_{d,s}|v-w|^{-d-2s}$ for some positive constant $c_{d,s}$
only depending on dimension and $s$. In the context of general
integro-differential equations, one must start by making sense of the
notions of uniform ellipticity, smoothness of coefficients, divergence form,
non-divergence form, weak solutions, viscosity solutions,
\textit{etc}. The following dictionary provides a basic understanding
of the different assumptions for integro-differential diffusions
that correspond to common structural conditions for second order
elliptic operators.
\begin{itemize}
\item \textbf{Uniform ellipticity} of order $2s$ corresponds to the
  bounds
  \[
    \lambda |v-w|^{-d-2s} \le K(t,v,w) \le \Lambda |v-w|^{-d-2s}
  \]
  for two positive constants $\lambda, \Lambda$.
  That is, the kernel should be comparable with that of the fractional
  Laplacian.
\item Equations in \textbf{divergence} form correspond to the symmetry
  condition $K(t,v,w) = K(t,w,v)$. In this case, the diffusion
  operator is self-adjoint in $L^2$.
\item Equations in \textbf{non-divergence} form correspond to the
  different symmetry condition $K(t,v,v+h) = K(t,v,v-h)$. In this
  case, the diffusion operator evaluates to zero when applied to an
  affine function, and evaluates to a locally bounded value when applied to a smooth function.
\end{itemize}
These definitions are the starting point for understanding
parabolic integro-differential equations. However, as we will see in
the rest of this article, the conditions on the kernel $K$ can be
significantly relaxed for each of the three notions above.

\subsection{The non-local diffusion of the Boltzmann equation}

The collision operator $Q(f,f)$ in the Boltzmann equation is a
nonlinear integro-differential diffusion in the velocity
variable. However, looking at the expression
\eqref{e:collision_operator}, there is no apparent similarity with the
equation \eqref{e:pide}. We resort to \emph{Carleman coordinates}
\cite{MR0098477,MR1555365} in order to rewrite $Q(f,f)$ as an
integro-differential diffusion with an $f$-dependent kernel, plus a
lower order term.

To avoid clutter, we omit writing explicitly the time $t$ and space
$x$ variables in every formula. Every expression we write is evaluated
for each fixed value of $t$ and $x$.

We reparametrize the integral in \eqref{e:collision_operator} in terms
of $w := v_\ast'-v$ and $v'$. The identities
\eqref{e:ellastic_collisions} are equivalent to
\begin{align*}
w &\perp (v' - v), \\
v_\ast &= v' + w.
\end{align*}
\begin{figure}[ht]
  \begin{center}
    \setlength{\unitlength}{1cm}
    \begin{picture}(4,5)
      \put(-3,0){\includegraphics[height=4cm]{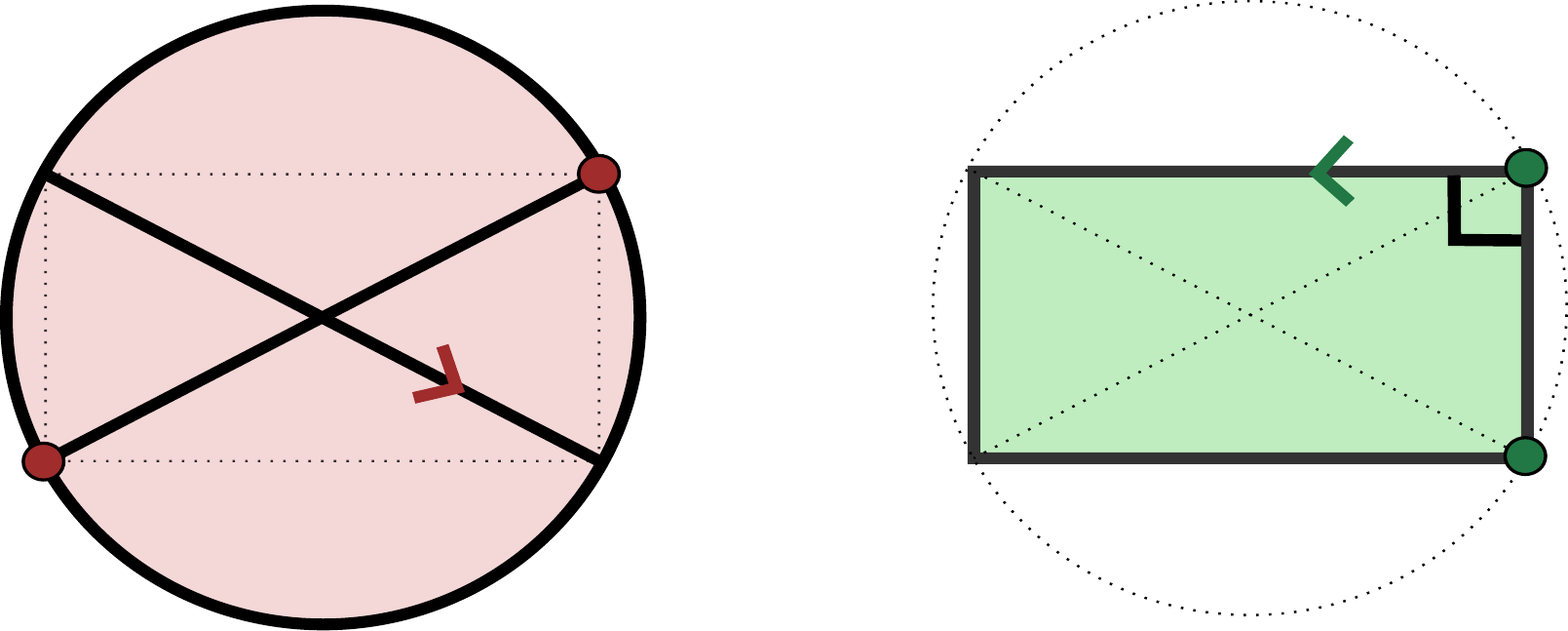}}
      \put(1.1,3.0){$v$}
      \put(7.0,3.0){$v$}
      \put(5.4,3.3){$w$}
      \put(7.0,1.0){$v'$}
      \put(-3.3,1.0){$v_\ast$}
      \put(-.6,1.2){$\sigma$}
    \end{picture}
  \end{center}
  \caption{On the left hand side, pre- and post-collisional velocities
    are parametrized by $v, v_\ast \in \R^d$ and
    $\sigma \in \partial B_1$. On the right hand side, Carleman
    coordinates are shown: velocities are parametrized by
    $v,v' \in \R^d$ and $w \perp (v'-v)$.}
  \label{fig:coord}
\end{figure}
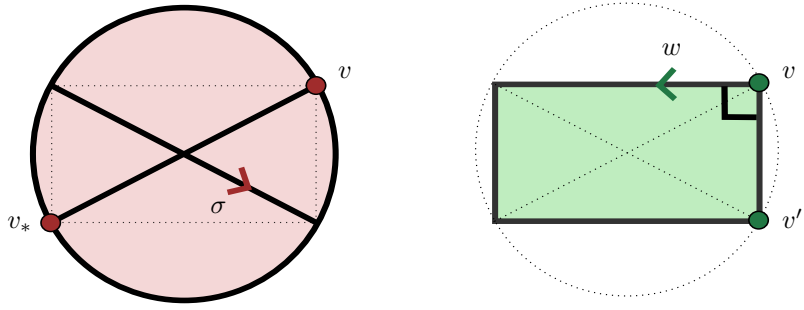
In terms of these variables, and taking into account the Jacobian of this change of variables, the operator $Q(f,f)$ from
\eqref{e:collision_operator} becomes
\begin{align*}
 Q(f,f) &= \int_{\R^d} \left( \int_{w \perp \{v'-v\}} (f(v') f(v'_\ast) - f(v) f(v_\ast) ) B(|v-v_\ast|,\cos \theta)  \frac{2^{d-1}}{|v'-v|}|v-v_\ast|^{-d+2} \d w \right) \d v'_\ast ,\\
&= \int_{\R^d} \bigg(f(v') K_f (v,v') - f(v) K_f(v',v) \bigg)  \dv'
\end{align*}
where the kernel $K_f$ depends on the function $f$ through the following
explicit formula (see \cite{silvestre2016new}),

\begin{equation}\label{e:Kf_exact}
 K_f (v,v') = \frac{2^{d-1}}{|v'-v|} \int_{w \perp v'-v} f(v+w) B(r,\cos \theta) r^{-d+2} \d w \quad \text{
with }
\begin{cases}
r^2 =|v'-v|^2 + |w|^2, \\
\cos \theta = \frac{w-(v-v')}{|w-(v-v')|}\cdot \frac{w+(v'-v)}{|w+(v'-v)|} .
\end{cases}
\end{equation}

Under the non-cutoff
assumption \eqref{e:B}, the kernel $K_f$ satisfies,
\begin{equation}
  \label{e:kf}
  \begin{aligned}
    K_f (v,v') &\approx |v-v'|^{-d-2s} \left\{ \int_{w \perp v'-v}  f(v+w) |w|^{\gamma +2s +1} \d w \right\}.
  \end{aligned}
\end{equation}
The sign $\approx$ means that the kernel $K_f$ is bounded from
  below and from above by the right hand side, up to a constant only
  depending on $B$.

  We observe in the formula \eqref{e:kf} that if the factor inside the
  brackets on the right hand side is $\simeq 1$, that is to say is
  bounded from above and below by positive constants, then the kernel
  $K_f$ is uniformly elliptic. However, we cannot ensure such a
  property based only on the hydrodynamic conditions in
  \eqref{e:hydro}. Our kernel $K_f$, a priori, may be a lot more
  degenerate than those considered in the earlier literature on
  integro-differential equations. There is no pointwise lower or upper
  bound for $K_f$ that can be deduced from \eqref{e:hydro}.

The well known \emph{cancellation lemma}, which appeared for the first time in \cite{advw}, tells us that
\[
  \int_{\R^d} (K_f(v',v) -  K_f (v,v')) \dv' = c_b \int_{\R^d}  f(v+w)|w|^\gamma \d w
\]
for some constant $c_b>0$ only depending on the collision kernel $B$. The operator $Q(f,f)$ can thus be rewritten under the following form
\[ Q(f,f) = \mathcal{L}_{K_f} f + f (f \ast c_b |\cdot |^\gamma) \]
with $K_f$ given by \eqref{e:kf} above. We write $\mathcal{L}_K$ to denote an integro-differential operator in $v$ associated to a kernel $K$,
\[ \mathcal{L}_K f (v) = \int_{\R^d} (f(v')-f(v)) K (v,v') \dv'. \]

The term $\mathcal L_{K_f} f$ is a nonlinear (since $K_f$ depends on
$f$) integro-differential diffusion. The term
$f (f \ast c_b |\cdot |^\gamma)$ is of lower order. The integro-differential diffusion
leads the smoothing effect of the equation.

In order to apply ideas from the area of integro-differential
equations in the context of the Boltzmann equation, there are several
difficulties that we must overcome. In particular, we must answer the
following questions.
\begin{enumerate}
\item In what way is the kernel $K_f$  elliptic? Can we generalize
  regularity results for integro-differential equations to possibly
  degenerate kernels like the ones for the Boltzmann equation?
\item Is the Boltzmann kernel $K_f$ in divergence or non-divergence form?
\item Is it possible to generalize the regularity results for
  integro-differential parabolic equations to the hypoelliptic setting
  of kinetic equations?
\end{enumerate}

In the next subsection, we will discuss the precise way in which the
kernel $K_f$ is elliptic, only in terms of the parameters of
\eqref{e:hydro}. A first regularity result for parabolic
integro-differential equations with such irregular kernels appeared in
\cite{schwab2016} in the form of a Krylov-Safonov type theorem. A
version of De Giorgi-Nash-Moser theorem for kinetic
integro-differential equations with possibly degenerate kernels
appeared in \cite{imbert2020weak}.

The Boltzmann kernel $K_f$ is naturally in non-divergence form, since
the identity $K(v,v+h) = K(v,v-h)$ is evident from its explicit
formula. The difference $K(v,v') - K(v',v)$ satisfies
cancellation conditions that allow us to also work with $\mathcal L_K$ as an
integro-differential operator in divergence form plus a lower order
correction.

In \cite{imbert2020weak} and \cite{imbert2020schauder}, we develop
integro-differential \textbf{kinetic} versions of the De Giorgi - Nash
- Moser theorem and of the Schauder estimates, making use of the
hypoelliptic relationship between the integral diffusion and the
transport terms. These are results for general kinetic
integro-differential equations that were developed with the explicit
purpose of applying them to our program of conditional regularity for
the Boltzmann equation. We explain them in Sections \ref{s:DG} and
\ref{s:schauder}

 \subsection{Non-degeneracy cones}
\label{s:cone}

The lower bound for the mass density of $f$ in \eqref{e:hydro}, combined with
the upper bound for the energy density, tells us that there is certain amount
of mass inside a ball $\{|v|<R\}$ (for $R$ depending on $m_0$ and
$E_0$). Moreover, the upper bound $H_0$ on entropy tells that
this mass cannot concentrate in a set of measure zero. A
quantification of this reasoning gives us that, for every value of $t$ and $x$,
\[
  \exists \ell , \mu , R >0 \quad / \quad |\{ f \ge \ell\}\cap \{|v|<R\}| \ge \mu .
\]

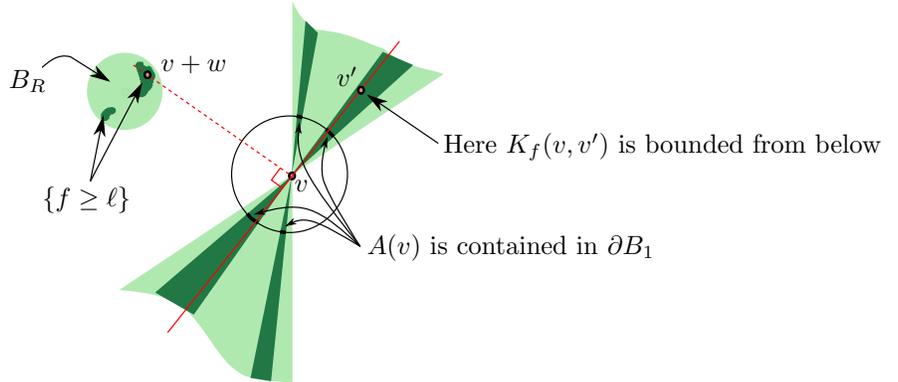
\begin{figure}[ht]
  \begin{center}
    \setlength{\unitlength}{1in}
    \begin{picture}(2.11111,2.2)
      \put(0,0){\includegraphics[height=2in]{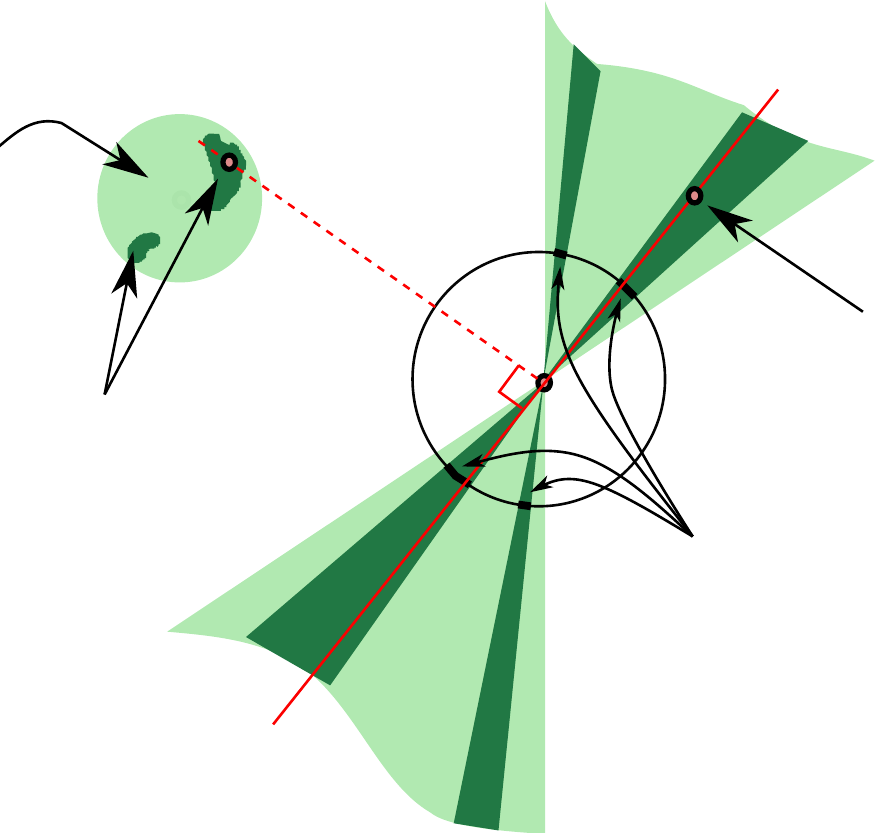}}
      \put(0.0,0.92){$\{f \geq \ell\}$}
      \put(-0.174,1.54){$B_R$}
      \put(1.32,1.){$v$}
      \put(1.54,1.55){$v'$}
      \put(0.62,1.63){$v+w$}
      \put(1.7,0.67){$A(v)$ is contained in $\partial B_1$}
%      \put(1.46,0.41){$A(v)$ lies in the strip}
%      \put(1.46,0.216){of width at most $C/|v|$}
      \put(2.1,1.20){Here $K_f(v,v')$ is bounded from below}
    \end{picture}
  \end{center}
  \caption{Non-degeneracy cone. It is made of lines passing through $v$ of velocities $v'$ such that $K_f (v,v')$ is bounded from
    below. In particular, it is centered at $v$.  It is generated in terms of the dark green set $\{f \ge \ell \} \cap B_R$. A generic line of the cone (shown in red on the figure) is made of velocities $v'$ such that the perpendicular plane $\{v+w : w \perp v'-v\}$ (shown in dashed red on the figure) significantly intersects the dark green set $\{ f \ge \ell \} \cap B_R$. See \eqref{e:kf}.  This figure is adapted from \cite{silvestre2016new}, see also \cite{imbert2019regularity}.}
  \label{fig:cone}
\end{figure}

In words, for every value of $t$ and $x$, there exists a set of positive measure and localized around
the origin where the function $f$ is bounded from below. This set
allows us to obtain a lower bound for the diffusion kernel $K_f$ in
some directions. Indeed, from any arbitrary point $v$, there is a symmetric cone
of directions whose perpendicular planes will intersect the set
$\{f \ge \ell\}$ on a set with $\mathcal{H}^{d-1}$ positive Hausdorff
measure. This is what we call the \emph{nondegeneracy cone} of $K_f$ at $v$
(See Figure \ref{fig:cone}). This cone is characterized by the fact that
$(v'-v)/|v'-v| \in A(v)$ for a certain subset $A(v)$ of the unit sphere
$S^{d-1}$. We have
\begin{equation} \label{e:cone_nondegeneracy}
 K_f(v,v') \geq \lambda(v) \,  |v-v'|^{-d-2s} \text{ whenever } v' \in v + \R A(v),
\end{equation}
for $\lambda(v) \approx (1+|v|)^{\gamma+2s+1} >0$ depending on $v$ and
the constants from \eqref{e:hydro}. Note that the nondegeneracy cone
depends on the function $f$ and on the point $v \in \R^d$. The
nondegeneracy cones rotate and stretch in some directions when we move
the point $v$.

The lower bound $K(v,v') \geq \lambda |v-v'|^{-d-2s}$, in the usual
ellipticity condition for integro-differential equations, holds only
when $v'$ belongs to the nondegeneracy cone emanating from $v$. This
is the only non-degeneracy condition that we can deduce from
\eqref{e:hydro}.

The usual uniform ellipticity condition for integro-differential
equations consists also of an upper bound for the form
$K(v,v') \leq \Lambda |v-v'|^{-d-2s}$. A pointwise upper bound like
that cannot be deduced from \eqref{e:hydro}. Instead, we can deduce an
upper bound in the following averaged sense. For every $r>0$, we have
\begin{equation} \label{e:K-upper}
 \int_{B_r} K(v,v+w) |w|^2 \d w \leq \Lambda(v) r^{2-2s}.
\end{equation}
Here, the value of $\Lambda(v) \approx (1+|v|)^{\gamma+2s}$ depends
only $M_0$ and $E_0$ in \eqref{e:hydro}, provided that
$\gamma+2s \in [0,2]$. Note that both parameters $\lambda(v)$ and
$\Lambda(v)$ in \eqref{e:cone_nondegeneracy} and \eqref{e:K-upper}
depend on $v$. They are locally uniformly bounded, but they do not stay bounded
as $|v| \to \infty$.

The conditions \eqref{e:cone_nondegeneracy} and \eqref{e:K-upper} are
weaker than the usual pointwise bounds
$\lambda |v-v'|^{-d-2s} \leq K(v,v') \leq \Lambda |v-v'|^{-d-2s}$ that
appear in the earlier literature as a notion of uniform ellipticity
appropriate for integro-differential equations.

\subsection{Structure of the Boltzmann equation}

The nonlinear integro-differential structure of the Boltzmann equation that we discussed so far can
be summarized in the following formula,
\begin{equation} \label{e:boltzmann_rewritten}
  \underset{\text{\normalsize free transport}}{\underbrace{\partial_t f + v\cdot \nabla_x f}} =
  \underset{\begin{minipage}{37mm} non-local diffusion, non-degenerate in several directions
    \end{minipage}}{\underbrace{ \mathcal{L}_{K_f} f }}
  + \underset{\text{\normalsize lower order term}}{\underbrace{f (f \ast c_b |\cdot |^\gamma)}}.
\end{equation}

We know that the diffusion is non-degenerate only in the cone of
directions discussed in Section~\ref{s:cone}. As we will see, this
cone of non-degeneracy is large enough to produce a regularizing
effect in the velocity variable $v$. We stress the fact that it is
Condition \eqref{e:hydro} that ensures the existence of such cones of
non-degeneracy.

\section{Global regularity estimates: a warm-up tour}
\label{s:warm_up}

The present section is dedicated to the presentation, in a few
paragraphs, of the steps involved in the proof of Theorem
\ref{t:main}. We describe the path from the assumption \eqref{e:hydro}
all the way to the a priori estimates for all derivatives of $f$.

Our program combines several ingredients. We summarize them in the following list.

\begin{itemize}
\item Pointwise upper bounds that show an arbitrarily fast polynomial
  decay as $|v| \to \infty$.
\item A weak Harnack inequality for kinetic integro-differential
  equations: it gives us local $C^\alpha$ estimates for some $\alpha>0$ (possibly small).
\item The Schauder estimates for kinetic integro-differential equations: they give us local $C^{2s+\alpha}$ estimates for some (small) $\alpha>0$.
\item A change of variables to adjust the ellipticity of the operators
  as $|v| \to \infty$: it turns our local H\"older and Schauder
  estimates into global ones.
\item An iterative gain in regularity, by successively applying the
  Schauder estimates to derivatives of solutions,
  to obtain $C^\infty$ estimates.
\end{itemize}

In this section, we briefly outline what each step does, and how they
combine to produce a proof of Theorem \ref{t:main}. We give a more
detailed explanation of each of these ingredients in later sections.

\subsection{Pointwise decay estimates in the velocity variable}

\label{s:warm-up_bounds}

Before studying the derivatives of a solution, we are interested in
its fast decay with respect to the velocity variable.

In \cite{silvestre2016new}, the second author of this paper showed that the
condition \eqref{e:hydro} suffices to obtain an a priori estimate in
$L^\infty$ for the solution $f$, in terms of the constants from
\eqref{e:hydro} only, provided that $\gamma+2s \geq 0$ and $\gamma \leq 2$.
For the range of parameters $\gamma+2s \in [0,2]$, in our joint work with Cl\'ement Mouhot \cite{imbert2020decay}, we improve the upper bound by
proving that solutions decay at an arbitrary algebraic rate as
$|v| \to \infty$. More precisely, given any solution $f$ of
\eqref{e:boltzmann}, with a non-cutoff collision kernel as in
\eqref{e:B}, if \eqref{e:hydro} holds, then for any value of
$q \geq 0$ and $\tau>0$, there is a constant $C_0$ such that
\begin{equation}
  \label{e:upper_bound}
  \sup_{(t,x,v) \in [\tau,T] \times \R^d \times \R^d} (1+|v|)^q | f
  (t,x,v)| \le C_0.
\end{equation}
The constant $C_0$ depends only on the parameters of \eqref{e:hydro}
when $\gamma>0$ (the so-called ``hard potentials'' case). In the case of
$\gamma \leq 0$, it is necessary to impose that the initial datum
$\fin (x,v) = f(0,x,v)$ enjoys a fast decay in $v$. The corresponding
constant $C_0$ then also depends on the constants measuring this
decay. In that sense, we say that the bounds \eqref{e:upper_bound} are
self generated for hard potentials ($\gamma>0$), and they are
propagated from the initial data for (moderately) soft potentials
($\gamma \leq 0$).

The structure of the proof is inspired by the classical idea of
barrier functions in elliptic PDEs. We set up a contradiction by
evaluating the equation at the first point of contact between the
solution $f$ and an upper barrier. However, the upper barrier is not a
super-solution of any particular equation. The contradiction is based
on purely nonlocal considerations in the spirit of
\cite{silvestre2006holder} or the nonlinear maximum principles in
\cite{constantin2012}.

We give a more detailed description of this result in Section \ref{s:upper_bounds}.

\subsection{A H\"older modulus of continuity}

Once the bound \eqref{e:upper_bound} is established, we are interested in estimating some H\"older modulus of continuity for the solution, with respect to all its variables. This will be done through a weak Harnack inequality, in the style of De Giorgi, for
general kinetic integro-differential equations. In \cite{imbert2020weak}, we study equations of the form,
\begin{equation}
  \label{e:lin-kin}
  \partial_t f + v \cdot \nabla_x f = \cLdiv_K f + h
\end{equation}
where
\[
  \cLdiv_K f(v) = \pv \int_{\R^d} [f(t,x,v')-f(t,x,v)] K(t,x,v,v') \d  v'.
\]
The kernel $K$ needs to satisfy only some mild notion of ellipticity
(implied by \eqref{e:cone_nondegeneracy} and \eqref{e:K-upper} in the
case of Boltzmann equation) and cancellation conditions. The source
function $h$ needs to be bounded. The H\"older regularity of the
solution $u$ is obtained by following a variant of De Giorgi's
method. Neither the kernel $K$, nor the source function $h$, need to
be linked to the solution $f$ by any additional formula. Because of
that, a local H\"older regularity estimate for an equation of the
form \eqref{e:lin-kin} is more general, and implies in particular a
local H\"older regularity estimate for the Boltzmann equation
\eqref{e:boltzmann_rewritten}.

Not long before our work in \cite{imbert2020weak}, the Harnack
inequality, and consequential H\"older estimates, as in the theorem of
De Giorgi, Nash and Moser, were obtained in the context of kinetic
equations with second order diffusion with rough coefficients (see
\cite{pascucci2004,wang2009,wang2011,gimv}). They apply to equations
of the form
\begin{equation} \label{e:kinetic_2nd_order}
 \partial_t f + v \cdot \nabla_x f = \partial_{v_i} \left[ a_{ij}(t,x,v) \partial_{v_j} f \right] +h,
\end{equation}
where the coefficients $a_{ij}$ are assumed to be uniformly elliptic
(i.e. $\lambda \mathrm{I} \leq \{a_{ij}\} \leq \Lambda \mathrm{I}$ for
some constants $\Lambda \geq \lambda > 0$), but they are not required
to satisfy any smoothness condition.

The equation \eqref{e:kinetic_2nd_order} would be hypoelliptic in the
sense of H\"ormander if its coefficients were smooth. In the present
case, the coefficients $a_{ij}$ are merely bounded and
measurable with respect to $t$, $x$ and $v$. Our result in
\cite{imbert2020weak} about an equation of the form \eqref{e:lin-kin}
was motivated by these results for \eqref{e:kinetic_2nd_order}, and in
particular by \cite{gimv}.

As we mentioned it in Section \ref{s:nonlocal_diffusion}, there are
several results available for integro-differential versions of the De
Giorgi-Nash-Moser theorem. Our result in \cite{imbert2020weak} can be
described as an integro-differential version of the main result in
\cite{gimv}. However, we face some very specific difficulties.

\begin{enumerate}
\item We deal with kernels that are significantly more singular than
  in the earlier literature. The conditions
  \eqref{e:cone_nondegeneracy} and \eqref{e:K-upper} only provide a
  very mild form of ellipticity.
\item The integral operator in the Boltzmann equation has the symmetry
  structure characteristic of equations in non-divergence form. In
  order to adapt to that, we substitute the natural symmetry
  condition, $K(t,x,v,v') = K(t,x,v',v)$, by a the following weaker cancellation
  condition.
\begin{align}
 \text{For all } r \in (0,1), \qquad \left \vert \int_{B_r(v)} \left( K(t,x,v,v') - K(t,x,v',v) \right) \d v' \right\vert &\leq C r^{-2s} \ \text{ and }  \label{e:cancellation1} \\
 \left \vert \int_{B_r(v)} (K(t,x,v,v') - K(t,x,v',v)) (v'-v)  \d v' \right\vert &\leq C r^{1-2s}.  \label{e:cancellation2}
\end{align}
In the context of the Boltzmann equation, the fact that the first
inequality holds is a reformulation of the classical cancellation
lemma from \cite{advw}. The second cancellation property seems to be
new.
\item The compactness argument which is at the core of the proof in
  \cite{gimv} cannot possibly be
  adapted to integro-differential equations of order less than one. In that case, we construct special barrier functions and use some ideas that originated in the work of Krylov and Safonov for parabolic equations in non-divergence form. More precisely, we adapt their well known covering technique that they call \emph{the growing ink spots lemma}.
\end{enumerate}

In a later paper \cite{stokols2019}, Logan Stokols shows that a
considerably simpler proof can be given when the kernels $K$ are
assumed to be symmetric (i.e. $K(t,x,v,v') = K(t,x,v',v)$) and
uniformly elliptic in the most classical sense:
$\lambda |v-v'|^{-d-2s} \leq K(t,x,v,v') \leq \Lambda
|v-v'|^{-d-2s}$. Such a result, however, cannot be used to derive
H\"older estimates for the Boltzmann equation.

The precise statement of our result, and more details on its proof,
will be given in Section \ref{s:DG}.

\subsection{Gain of $2s$ derivatives through Schauder's approach}

The classical Schauder theory gives us an estimate of the
$C^{2+\alpha}$ norm of solutions to linear uniformly elliptic
equations with $C^\alpha$ coefficients and $C^\alpha$ source terms.

Linear kinetic equations with second order diffusion are a particular
case of the more general theory of ultraparabolic equations of
Kolmogorov type. The Schauder theory for this type of equations was
developed mostly in the late 1990's and early 2000's. See
\cite{sayro1971,manfredini, lunardi1997, eidelman1998,
  dfp,radkevich2008}, and the survey article
\cite{sergio2004recent}. In particular, it applies to equations of the
form
\begin{equation} \label{e:kinetic_non_divergence}
 \partial_t f + v \cdot \nabla_x f = a_{ij}(t,x,v) \partial_{v_i v_j} f  +h,
\end{equation}
when the coefficients $a_{ij}$ are uniformly elliptic and these
  functions together with the source term $h$ are H\"older
  continuous.

The Schauder estimates are a powerful tool to bootstrap higher
regularity estimates once an initial H\"older estimate is established
for a nonlinear parabolic equation. A nonlinear equation can often be
written as a parabolic or elliptic equation whose coefficients depend
on the solution itself. If we know a H\"older modulus of continuity
for the solution, it may imply a H\"older bound for the coefficients. The Schauder estimates then give us a $C^{2+\alpha}$
estimate for the solution. We thus  know that the coefficients are even
more regular, and we iterate. In other words, the Schauder estimate allows
one to \emph{gain two derivatives}: starting from the control of the
modulus of continuity of the solution, we reach a control of the
modulus of continuity of second order derivatives.

In \cite{imbert2020schauder}, we obtain Schauder estimates for kinetic
integro-differential equations of the form \eqref{e:lin-kin}. The
method of our proof is very different from the earlier work on
ultraparabolic equations. Instead, we borrow ideas from the
\emph{blow-up} technique developed by Joaquim Serra \cite{serra2015} for integro-differential equations.

Classical Schauder theory for parabolic equations involves H\"older
spaces whose definition encodes the parabolic scaling, through the
introduction of a parabolic distance. Turning to kinetic equations
with nonlocal diffusion in the velocity variable, a new scaling
appears naturally that is different for time, space and velocity. The
class of equations is not translation invariant anymore, but rather
Galilean invariant. The H\"older spaces must take these invariances
into consideration. The success of an appropriate Schauder theory
depends a great deal on finding the right definition for kinetic
H\"older spaces and an appropriate kinetic distance. These concepts,
together with the precise formulation of the kinetic
integro-differential Schauder estimates, will be given in Sections
\ref{s:cylinders_and_holder} and \ref{s:schauder}.

Applying this Schauder theory to the Boltzmann equation is then
possible \cite{imbert2020schauder}. Like we described in the general
framework above, the Schauder estimates are applied iteratively to
gain higher regularity estimates starting from the initial H\"older
estimate. In each application of the kinetic integro-differential
Schauder estimates, we gain $2s$ derivatives in velocity, $1$
derivative in time, and $2s / (1+2s)$ derivatives in space.

\subsection{Bootstrap}

The application of the kinetic integro-differential version of the De
Giorgi-Nash-Moser theory developed in \cite{imbert2020weak} gives us
localized H\"older estimates for the solution of the Boltzmann
equation \eqref{e:boltzmann}. They apply provided that $v$ is
contained in some bounded ball $B_R$. This is because the ellipticity
properties of the kernel $K$ given in \eqref{e:cone_nondegeneracy} and
\eqref{e:K-upper} degenerate as $|v| \to \infty$.

In order to improve this initial regularity estimate by successive
application of the Schauder estimates, it is necessary to turn our
local H\"older estimates into global ones. Inspired by an idea in
\cite{cameron2018}, we devised in \cite{imbert2019regularity} a change
of variables that transforms the kinetic integro-differential equation
into one whose ellipticity parameters are uniform in $v$. This change
of variables allows us to turn the local regularity estimates, from De
Giorgi and Schauder theories, into global ones that hold for all
velocities. It provides a key ingredient for the proof of Theorem
\ref{t:main} in \cite{imbert2019regularity}.

The application of the Schauder estimates, together with the change of
variables, gives us a global gain of regularity. We gain $2s$
derivatives in velocity, $1$ derivative in time, and $2s/(1+2s)$
derivatives in space. Then we compute an equation for discrete
incremental quotients of the solution, and apply the Schauder
estimates again. An iteration of this procedure leads to the
$C^\infty$ estimates in Theorem \ref{t:main}. In each iteration, we
gain a certain fractional number of derivatives in $v$, $x$ and
$t$. However, we also loose a decay power. More precisely, an upper
bound for $(\partial_t)^m f$ that decays like $\lesssim (1+|v|)^{-q}$
will depend on an earlier upper bound on $(\partial_t)^{m-1} f$ that
decays like $\lesssim (1+|v|)^{-\tilde q}$, for some $\tilde q >
q$. Since we start with upper bounds that decay arbitrarily fast in
\eqref{e:upper_bound}, the iteration continues forever.

We explain the change of variables in Section \ref{s:change_of_vars} and the iteration procedure finalizing the proof of Theorem \ref{t:main} in Section \ref{s:bootstrap}.

\section{Generic kinetic equations with integral diffusion}
\label{s:generic}

The study of general kinetic equations with integral diffusion plays a
important role in the derivation of the global regularity estimates
stated in Theorem~\ref{t:main}.  A kinetic equation with integral
diffusion takes the following form,
\begin{equation}
  \label{e:lin-kin-2}
  \partial_t f + v \cdot \nabla_x f = \cL_K f + h,
\end{equation}
where the integral diffusion $\cL_K$ depends on a kernel $K$ and is
given by the formula,
\begin{equation}
  \label{e:integro-differential_operator}
  \cLdiv_K f(v) = \pv \int_{\R^d} [f(t,x,v')-f(t,x,v)] K(t,x,v,v') \d v'.
\end{equation}

The integral diffusion operator $\cL_K$ acts similarly as a classical
diffusion operator in the $v$ variable but it is nonlocal, of
fractional order. It can be compared with a second order operator in
divergence form $\partial_{v_i} (a_{ij}(t,x,v) \partial_{v_j} f)$ or
one in non-divergence form $a_{ij}(t,x,v) \partial_{v_j v_j}
f$. Classical second order elliptic operators are studied using
different tools depending on these two structures. Test functions and
estimates in Sobolev spaces are typical of operators in divergence
form, whereas barrier functions and comparison principles are typical
of equations in non-divergence form. For integro-differential
equations, the operators are always written with the same formula
above. The divergence \textit{vs} non-divergence structures are
determined by two alternative symmetry assumptions on the kernel $K$:
\begin{itemize}
\item Divergence: $K(t,x,v,v') = K(t,x,v',v)$.
\item Non-divergence: $K(t,x,v,v+w) = K(t,x,v,v-w)$.
\end{itemize}
The kernel $K_f$ corresponding to the Boltzmann equation (as in
\eqref{e:Kf_exact}) is naturally in non-divergence form. In general,
we have $K_f(t,x,v,v') \neq K_f(t,x,v',v)$. This is an obstruction in
order to apply \emph{divergence} techniques. However, the Boltzmann
kernel satisfies the cancellation conditions
\eqref{e:cancellation1} and \eqref{e:cancellation2}, that are a weaker form of the \emph{divergence} symmetry assumption. It turns out that
these cancellation conditions suffice in order to derive the most
crucial estimate for equations in divergence form: the De
Giorgi-Nash-Moser (weak) Harnack inequality.

The reader should not be misled into thinking that any result for an
equation of the form \eqref{e:lin-kin-2} would apply to linear
equations only. A nonlinear kinetic equation, like the Boltzmann
equation, also has the form \eqref{e:lin-kin-2} with the bonus piece
of information that the kernel $K$ and the source term $h$ are related
to the solution $f$ by certain formulas. Any a priori estimate for
solutions of \eqref{e:lin-kin-2}, that depends on minor assumptions on
$K$ and $h$, will provide us with estimates in particular for the
Boltzmann equation as part of a more general family of equations. The
key is to analyze the equation \eqref{e:lin-kin-2} with minimalistic
assumptions on the kernel $K$ and the source term $h$, so that we can
verify those assumptions in the case of the Boltzmann equation
depending only on the hydrodynamic bounds in \eqref{e:hydro}.

We develop two fundamental regularity techniques in the context of
generic kinetic integral equations like \eqref{e:lin-kin-2}. The weak
Harnack inequality in the spirit of De Giorgi gives us H\"older
continuity estimates for $f$ in terms of quantitative bounds for $K$
and $h$ only. It applies to equations in \emph{divergence form}, that
in our context appears in the form of the cancellation conditions
\eqref{e:cancellation1} and \eqref{e:cancellation2}. The Schauder
estimates give us regularity estimates in higher order H\"older
spaces, once we have an initial H\"older continuity estimate for the
source term and the kernel. Schauder estimates apply naturally to
equations in non-divergence form.

\subsection{Ellipticity conditions}

The notion of ellipticity for integro-differential operators is more
subtle than for second order differential equations. A classical
second order operator involves a matrix of coefficients multiplying
the second derivatives of a function, and could be written in
divergence or non-divergence form. The ellipticity in the classical
case reduces simply to upper and lower bounds on the matrix of
coefficients. An integro-differential operator like $\cL_K$ is defined
in terms of a nonnegative kernel $K$, which is typically singular
around the origin. It should be thought of as a nonlocal diffusion
operator in the variable $v$, which is applied for every fixed value
of $t$ and $x$. With this point of view, we state the ellipticity
assumptions for kernels $K(v,v')$ depending on $v$ and $v'$
only. Ultimately, we will require our full kernel $K(t,x,v,v')$ to
satisfy these conditions uniformly in $t$ and $x$.

The first notion of ellipticity that appeared in the literature of
nonlocal equations consist in a pointwise comparability between
$K(v,v')$ and the kernel of the fractional Laplacian. The majority of
the results on regularity of nonlocal equations that appeared between
the years 2000 and 2015 depend on the assumption:
$\lambda |v-v'|^{-d-2s} \leq K(v,v') \leq \Lambda |v-v'|^{-d-2s}$. In
the case of the Boltzmann equation, there is no way to establish these
bounds for the kernel $K_f$ of \eqref{e:Kf_exact} in terms of the
parameters of \eqref{e:hydro} only. So, we are forced to consider more
general ellipticity conditions, that are harder to work with.

Instead of the pointwise upper bound
$K(v,v') \leq \Lambda |v-v'|^{-d-2s}$, we require that this bound
holds \emph{on average} only. That is, we require that for all $r>0$,
\begin{equation} \label{e:ellipticity_upper}
 \int_{B_r(v)} K(v,v') |v-v'|^2 \d v' \leq \Lambda r^{2-2s},
\end{equation}
for every $v$ in the domain of the equation.

In terms of the lower bound, as we described in Section \ref{s:cone},
the pointwise lower bound for the Boltzmann equation holds on a
symmetric cone of nondegeneracy emanating from each point $v$. That
is, for every $v$ in the domain of the equation, there exists a subset
$A \subset S^{d-1}$ so that
\begin{equation} \label{e:ellipticity_cone}
|A|_{\mathcal H^{d-1}} \geq \mu, A = -A \qquad
K(v,v+w) \geq \lambda |w|^{-d-2s} \ \text{ if } w/|w| \in A.
\end{equation}

One may argue that \eqref{e:ellipticity_cone} does not seem to be such
a weaker replacement of the lower bound
$K(v,v') \geq \lambda |v-v'|^{-d-2s}$ as \eqref{e:ellipticity_upper}
is of the upper bound $K(v,v') \leq \Lambda |v-v'|^{-d-2s}$. Indeed,
it is not clear what the sharpest notion of ellipticity should be. The
important feature of \eqref{e:ellipticity_upper} and
\eqref{e:ellipticity_cone} is that they can be verified to hold for
the Boltzmann kernel $K_f$ of \eqref{e:Kf_exact} with parameters
$\lambda>0$, $\Lambda$ and $\mu>0$ depending on the constants of
\eqref{e:hydro} only.

Requiring that \eqref{e:ellipticity_upper} and
\eqref{e:ellipticity_cone} hold with uniform constants $\mu>0$,
$\lambda>0$ and $\Lambda$, for every $(t,x,v)$ in the domain of the
equation, is a notion of uniform ellipticity that is weaker than the
usual assumptions in the literature. It is the assumption we will work
with, in order to be able to apply our general estimates to the
Boltzmann equation.

The ellipticity conditions \eqref{e:ellipticity_upper} and
\eqref{e:ellipticity_cone} must be accompanied with symmetry
assumptions depending on whether we want to apply methods from
divergence or non-divergence equations. In the case of non-divergence
equations, we would require the following symmetry condition
\begin{equation} \label{e:symmetry_non-divergence}
 \text{non-divergence symmetry condition: } K(v,v+w) = K(v,v-w).
\end{equation}
Naturally, the condition \eqref{e:symmetry_non-divergence} should hold
for every $t$, $x$ and $v$ in the domain of the equation. We recall
that we are omitting writing the $t$ and $x$ dependence on $K$ since
the operator $\cL_K$ is applied for each fixed value of $t$ and $x$.

Luckily, the expression \eqref{e:Kf_exact} for the kernel of the
Boltzmann equation satisfies the non-divergence symmetry condition
\eqref{e:symmetry_non-divergence}. The Schauder estimate (described
below in Section \ref{s:schauder}) requires this symmetry condition,
and will apply to the Boltzmann equation.

The natural symmetry condition to apply methods for equations in
divergence form would be $K(v,v') = K(v',v)$. Indeed, a characteristic
property of second order operators in divergence form (of the form
$\partial_i a_{ij} \partial_j f$) is that they are self-adjoint in
$L^2$. It is easy to see that $\cL_K$ will be self-adjoint if and only
if $K(v,v') = K(v',v)$. Unfortunately, this symmetry condition does
not hold for the Boltzmann kernel $K_f$ of \eqref{e:Kf_exact}. Thus,
we are forced to consider a more general condition, that is
naturally harder to work with. To that aim, we state the following
cancellation conditions: there exist a constant $\Lambda$ so that for
all $v$ in the domain of the equation the following inequalities hold
for all $r \in (0,1)$,
\begin{equation} \label{e:symmetry_cancellation}
\begin{aligned}
\left\vert \int_{B_r(v)} \left( K(v,v') - K(v',v) \right) \, \d v' \right\vert &\leq \Lambda r^{-2s},  \\
\text{if $s \geq 1/2$, } \ \left\vert  \int_{B_r(v)} (K(v,v') - K(v',v))(v-v') \, \d v' \right\vert &\leq \Lambda r^{1-2s}.
\end{aligned}
\end{equation}

Note that for $s \in (0,1/2)$ there is only one cancellation
condition, whereas for $s \in [1/2,1)$ both inequalities are supposed
to hold. In fact, for $s \in (0,1/2)$, the second cancellation
condition follows as a consequence of the first one, combined with
\eqref{e:ellipticity_upper}.

There are a few extra inequalities that are a consequence of the
ellipticity conditions \eqref{e:ellipticity_upper} and
\eqref{e:ellipticity_cone} and the symmetry conditions
\eqref{e:symmetry_cancellation}. First of all, the upper bound
\eqref{e:ellipticity_upper} can be rephrased (by adjusting the
constant $\Lambda$ as necessary) in any of the following equivalent
alternative formulations,
\begin{align*}
\int_{B_{2r}(v) \setminus B_r(v)} K(v,v') \, \d v' \leq \Lambda r^{-2s}, \\
\int_{\R^d \setminus B_r(v)} K(v,v') \, \d v' \leq \Lambda r^{-2s}.
\end{align*}
Combining \eqref{e:ellipticity_upper} with the first inequality in
\eqref{e:symmetry_cancellation}, we get also, for all $r \in (0,1)$,
\[ \int_{B_r(v)} K(v',v) |v-v'|^2 \, \d v' \leq \tilde \Lambda r^{2-2s},\]
for a constant $\tilde \Lambda$ depending on $\Lambda$.

The conditions \eqref{e:ellipticity_upper} and
\eqref{e:symmetry_cancellation} imply that the operator $\cL_K$ maps
$H^s$ into $H^{-s}$. The following result is proved in
\cite{imbert2020weak}.

\begin{proposition} \label{p:Hs2H-s} Assume $K$ is a kernel for which
  \eqref{e:ellipticity_upper} and \eqref{e:symmetry_cancellation} hold
  for every $v \in \R^d$. Then, for any pair of functions
  $f, g \in H^s(\R^d)$,
\[ \int_{\R^d} (\cL_K f) g \, \d v \leq C \|f\|_{H^s(\R^d)} \, \|g\|_{H^s(\R^d)},\]
for a constant $C$ depending on $\Lambda$, $s$ and dimension only.
\end{proposition}

Proposition \ref{p:Hs2H-s} indicates that it is fair to think of the
operator $\cL_K$ as a nonlocal operator of order $2s$. Its ellipticity
is justified by the following proposition.
\begin{proposition} \label{p:coercivity} Assume $K$ is a kernel for
  which \eqref{e:ellipticity_cone} holds for every $v \in \R^d$. Then,
  for every $f \in H^s(\R^d)$, we have
\[ \iint_{\R^d \times \R^d} |f(v') - f(v)|^2 K(v,v') \, \d v' \d v \geq c \|f\|_{\dot H^s(\R^d)}^2,\]
for a constant $c>0$ depending on $\mu$, $\lambda$, $s$ and dimension only.
\end{proposition}

The following identity follows by a straight forward manipulation of the integral expression of \eqref{e:integro-differential_operator} and applying Fubini's theorem.
\[ -\int (\cL_K f) f \, \d v = \frac 12 \iint_{\R^d \times \R^d} |f(v') - f(v)|^2 K(v,v') \, \d v' \d v + \frac 12 \int_{\R^d} \left( \int_{\R^d} K(v',v) - K(v,v') \, \d v' \right) f(v)^2 \d v. \]
Note that the last term on the right hand side vanishes if the kernel satisfies the symmetry condition $K(v,v')=K(v',v)$. Otherwise, the cancellation condition \eqref{e:symmetry_cancellation} together with the upper bound \eqref{e:ellipticity_upper} allow us to bound this term as $\lesssim \|f\|_{L^2}^2$. Thus, Proposition \ref{p:coercivity} provides a coercivity estimate for the operator $\cL_K$ in terms of $\|f\|_{\dot H^s}$ minus a lower order correction.

Proposition \ref{p:coercivity} is proved in \cite{chaker2019coercivity} under more general conditions.

Naturally, when the conditions \eqref{e:ellipticity_upper},
\eqref{e:ellipticity_cone} and \eqref{e:symmetry_cancellation} hold
only on some subdomain of $\R^d$, then  appropriately localized
versions of Propositions \ref{p:Hs2H-s} and \ref{p:coercivity} hold as
well.

Note that Propositions \ref{p:Hs2H-s} and \ref{p:coercivity} apply to generic integro-differntial operators of the form \ref{e:integro-differential_operator}. We can verify that when $K$ is the kernel associated to the Boltzmann equation (as in \eqref{e:Kf_exact}), the assumptions \eqref{e:ellipticity_upper}, \eqref{e:ellipticity_cone} and \eqref{e:symmetry_cancellation} hold locally only in terms of the constants in \eqref{e:hydro}, but they do not hold uniformly for all velocities in $\R^d$.

\section{Cylinders and  H\"older spaces}
\label{s:cylinders_and_holder}

\subsection{Invariant transformations}
\label{s:invariance}

Here, we study the transformations that keep the class of equations of
the form \eqref{e:lin-kin-2} invariant. We describe two types of
transformations: scaling and Galilean translations.

We first describe the scaling of the
equation. Given any $r>0$, let us define
$S_r: \R^{2d+1} \to \R^{2d+1}$ by the following formula
\[ S_r(t,x,v) = (r^{2s} t, r^{1+2s} x, rv).\]

Suppose that $f$ is a solution of the equation \eqref{e:lin-kin-2} in
some domain. Then, we verify by a direct computation that for any
constants $a, r>0$, the function
\[ f_{a,r}(t,x,v) := a f (S_r(t,x,v)) \]
solves an equation of the same form in an appropriately scaled domain
with the modified kernel
\[ K_{r}(t,x,v,v')  := r^{d+2s} K(r^{2s} t, r^{1+2s} x, r v , r v'),\]
and the modified source term
\[ h_{a,r}(t,x,v)  := a h(S_r(t,x,v)).\]

The importance of the choice of exponents in $S_r$ is that if the
kernel $K$ satisfies the ellipticity conditions
\eqref{e:ellipticity_upper} and \eqref{e:ellipticity_cone}, then
$K_{r}$  also satisfies the same conditions with the same constants. In
the space-homogeneous case (that is, when $f$ does not depend on $x$),
if $s=1$, the scaling $S_r$ coincides with the usual parabolic scaling
$(t,v) \mapsto (r^2 t, rv)$. The scaling exponent we write here is
properly adjusted to operators of order $2s$ and kinetic equations.

The cancellation condition \eqref{e:symmetry_cancellation} is not
exactly preserved by the scaling since the restriction $r \in (0,1)$ in \eqref{e:symmetry_cancellation}
would become $r \in (0,1/\tilde r)$ after the transformation $S_{\tilde r}$. In fact, the
condition \eqref{e:symmetry_cancellation} is \emph{subcritical} since
it becomes stronger as we focus on small scales with $\tilde r \ll 1$.

Because of the $v$-dependence in the second term in
\eqref{e:lin-kin-2}, the class of equations is not translation
invariant in the usual way. Instead, it is Galilean invariant. For a
given $z_0=(t_0,x_0,v_0) \in \R^{1+2d}$, let us consider the Lie group
operator $z_0 \circ (t,x,v) = (t_0+t,x_0+x+tv_0,v_0+v)$. The
correction term $tv_0$ in the $x$ variable accounts for the change of
the position coordinate when passing from a motionless frame to
another one moving at a constant speed $v_0$. If $f$ is a solution of
the equation \eqref{e:lin-kin-2} in some domain, we consider its
Galilean translation. The function
\[ f_{z_0}(t,x,v) = f( z_0 \circ (t,x,v) )\]
solves an equation of the same form with the modified kernel
\[ K_{z_0}(t,x,v,v')  := K(z_0 \circ (t,x,v) , v_0+v'),\]
and the modified source term
\[ h_{z_0}(t,x,v) = h( z_0 \circ (t,x,v)).\]
Again, the domain of the equation has to be translated
accordingly. When the kernel $K$ satisfies any one of the conditions
\eqref{e:ellipticity_upper}, \eqref{e:ellipticity_cone},
\eqref{e:symmetry_non-divergence} and/or
\eqref{e:symmetry_cancellation}, the same holds for $K_{z_0}$.

\subsection{Cylinders}
\label{s:cylinders}

Given $r>0$ and $z_0 \in \R^{1+2d}$ and in view of the scaling
$S_r\colon \R^{1+2d} \to \R^{1+2d}$ and the Galilean invariance
$z \mapsto z_0 \circ z$, both defined above, it is natural to define \emph{cylinders}
$Q_r (z_0)$ centered at $z_0$ of radius $r>0$ as follows,
\[
  Q_r (z_0) = \{ z_0 \circ S_r (z) : z \in (-1,0] \times B_1 \times B_1 \}.
\]
If $z_0= (t_0,x_0,v_0)$, it is equivalent to the following definition
\[
  Q_r (z_0) = \{ (t,x,v) \in \R^{1+2d}: -r^{2s} < t-t_0 \le 0, |x-x_0-(t-t_0)v_0| < r^{1+2s}, |v-v_0| < r \} .
\]
Like in parabolic theory, the reference point $z_0$ for the cylinder is at the final time $t_0$. The cylinder $Q_r(z_0)$ includes points at earlier times than $t_0$ but not on its future.

\subsection{Kinetic H\"older spaces}
\label{s:holder}

Here, we describe an appropriate notion of H\"older space that is adapted to the scaling and Galilean invariance of the equations.

In order to provide a proper definition of H\"older spaces with any exponent $\alpha>0$, we must start with a modified notion of
degree for polynomials in $\R [t,x,v]$. Given a nonzero monomial
$m \in \R [t,x,v]$, the kinetic degree $\degk m$ is the number
$\kappa$ so that for all $z \in \R^{1+2d}$, $r>0$, we have
$m (S_r (z))=r^\kappa m(z)$. For a general polynomial $p = \sum m_k$,
the kinetic degree $\degk p$ is defined as the maximal kinetic degree
of the monomials $m_k$.

Roughly speaking, every exponent of $t$ counts as $2s$, every exponent of $x_i$ counts as $1+2s$ and
every exponent of $v_i$ counts as $1$. A monomial $m(t,x,v) =  a t^{k_0} x_1^{k_1} \dots x_d^{k_d} v_1^{k_{d+1}} \dots v_d^{k_d}$ has kinetic degree equal to $2s k_0 + (1+2s) (k_1 + \dots + k_d) + (k_{d+1} + \dots + k_{2d})$.

We notice that the kinetic degree of a nonzero polynomial can be any number of the discrete
set $\N + (2s) \N$. We adopt the convention that the kinetic degree of the zero polynomial equals $-\infty$.

With the definition of kinetic degree at hand, we can now define
kinetic H\"older spaces.
%----------------------------------------------------------------------
\begin{defi}[H\"older spaces]
  Given $\alpha \in (0,+\infty)$ and a open set $D \subset \R^{1+2d}$,
  we say that a function $f \colon D \to \R$ is $\alpha$-H\"older
  continuous in $D$ if there is some constant $C$ so that for any
  kinetic cylinder $Q_r(z_0)$, there exists some polynomial $p$ of
  kinetic degree strictly smaller than $\alpha$ such that
  \[
    |f(z) - p(z)| \le C r^\alpha \qquad \text{for all } z \in Q_r(z_0) \cap D.
  \]

  The set of $\alpha$-H\"older continuous functions $f\colon D \to \R$
  is denoted by $\cCl^\alpha (D)$.

  The least positive constant $C$ so that the inequality above holds
  is denoted by $[f]_{\cCl^\alpha (D)}$.
\end{defi}
% ----------------------------------------------------------------------
Note that with the definitions above, given any continuous function
$f : D \to \R$, the seminorm $[f]_{\cCl^0 (D)}$ is precisely the
supremum norm: $[f]_{\cCl^0 (D)} = \sup_{D} |f|$. We also define the
norms
$\|f\|_{\cCl^\alpha(D)} := [f]_{\cCl^\alpha(D)} + [f]_{\cCl^0(D)}$.

We write the subindex $\ell$ in $\cCl^\alpha$ to stress the fact that these norms are \textbf{l}eft-invariant by the action of the Galilean group, and not right-invariant.

The $\cCl^\alpha$ semi-norms encode a H\"older continuity behavior
that is compatible with the scaling and the Galilean invariances of the
equation described above. For a function $f(v)$ that depends on $v$
only, they would coincide with the usual $C^\alpha$ norms. For a
function $f(t)$ that depends only on $t$, it would rather correspond
to the $C^{\alpha/(2s)}$ norm. And for a function $f(x)$, depending
only on $x$, it would correspond to the $C^{\alpha/(1+2s)}$ norm. The
$\cCl^\alpha$ norm is left-invariant by the action of the Galilean Lie
group. It is \textbf{not} right-invariant.

Kinetic H\"older norms satisfy many of the same formal relationships
as the usual H\"older norms. For example, the following interpolation
inequality (proved in \cite{imbert2020schauder}) looks very much like the classical one.

\begin{proposition}[Interpolation inequalities]
  \label{p:interpol}
  Given $0 \le \alpha_1 < \alpha_2 < \alpha_3$ so that
  $\alpha_2 = \theta \alpha_1 + (1-\theta) \alpha_3$, a cylinder
  $Q_r(z_0)$ and a function $f \in \cCl^{\alpha_3} (Q_r(z_0))$,
  \[
    c[f]_{\cCl^{\alpha_2} (Q_r (z_0))} \leq [f]_{\cCl^{\alpha_1 }(Q_r
        (z_0))}^\theta [f]_{\cCl^{\alpha_3} (Q_r (z_0))}^{1-\theta} +
    r^{-(\alpha_2-\alpha_1)} [f]_{\cCl^{\alpha_1} (Q_r (z_0))}
  \]
  for some constant $c$ only depending on dimension.
\end{proposition}

It is also true that if $f \in \cCl^\alpha$, then derivatives of $f$ will
belong to a H\"older space with a smaller exponent. In this case, we
must account for differential operators that are left invariant by the
Lie group, and their kinetic degree should be properly accounted for.

\begin{proposition} \label{p:holder_derivatives}
Let $f \in C^\alpha(Q)$ for some kinetic cylinder $Q$. Then
\begin{itemize}
\item If $\alpha \geq 2s$, $(\partial_t + v \cdot \nabla_x) f \in C^{\alpha-2s}$ and
\[ [(\partial_t + v \cdot \nabla_x) f]_{C^{\alpha-2s}(Q)} \lesssim [f]_{C^\alpha(Q)}.\]
\item If $\alpha \geq 1+2s$, $\partial_{x_i} f \in C^{\alpha-1-2s}$ and
\[ [\partial_{x_1} f]_{C^{\alpha-1-2s}(Q)} \lesssim [f]_{C^\alpha(Q)}.\]
\item If $\alpha \geq 1$, $\partial_{v_i} f \in C^{\alpha-1}$ and
\[ [\partial_{v_1} f]_{C^{\alpha-1}(Q)} \lesssim [f]_{C^\alpha(Q)}.\]
\end{itemize}
\end{proposition}

Note that the statement of Proposition \ref{p:holder_derivatives}
involves the operator $(\partial_t + v \cdot \nabla_x)$ and not the
plain time derivative $\partial_t$. The time derivative is not
left-invariant by the Lie group structure (it is right-invariant). In
practice, one can still compute an estimate for the H\"older norm of
$\partial_t f$ in a bounded cylinder by combining the first two bullet
points in Proposition \ref{p:holder_derivatives} and the triangle
inequality. However, such an estimate would depend on the size of the
cylinder, and the H\"older exponent will not be better than
$\alpha - 1 -2s$.

These basic properties of kinetic H\"older spaces are stated and
proved in \cite{imbert2020schauder}. Some further analysis of H\"older
norms is continued in \cite{imbert2019regularity}.

\section{H\"older estimates via De Giorgi's method}
\label{s:DG}

In the study of regularity estimates for elliptic and parabolic quasilinear equations (away from equilibrium), the first step is usually to apply the method of De Giorgi, Nash and Moser. The key of this theorem is that it has no regularity requirement for the coefficients of the equation. In a quasilinear equation, the coefficients depend on the solution itself. Thus, at the beginning of our analysis we typically have little or no information on how smooth these coefficients might be. The Boltzmann equation, as described in \eqref{e:boltzmann_rewritten}, can be thought of as a quasilinear evolution equation. It is however not an ordinary parabolic equation, but a kinetic and nonlocal one.

In this section we describe a local H\"older estimate, in the style of
the classical theorems of De Giorgi, Nash, and Moser, but this time
for non-local kinetic equations. The precise statement is the
following.
\begin{thm}[Local H\"older estimate] \label{t:holder} Let $f$ be a
  bounded function that solves \eqref{e:lin-kin-2} in
  $Q_1 = (-1,0] \times B_1 \times B_1$. Assume that the kernel $K$ is
  a nonnegative function defined in
  $(-1,0] \times B_1 \times B_2 \times \R^d$ so that
  \eqref{e:ellipticity_upper}, \eqref{e:ellipticity_cone} and
  \eqref{e:symmetry_cancellation} hold. Then, $f$ is H\"older
  continuous in the  cylinder $Q_{1/2}$ with
  \[ [f]_{\cCl^\alpha(Q_{1/2})} \leq C \left( \|f\|_{L^\infty((-1,0]
        \times B_1 \times \R^d)} + \|h\|_{C^0(Q_1)} \right),\] for
  some constants $C$ and $\alpha>0$ depending on dimension, (a lower
  bound for) $s$, and the ellipticity parameters $\mu$, $\lambda$ and
  $\Lambda$.
\end{thm}

Note that the $L^\infty$ norm of $f$ in the right hand side is
evaluated in the set $(-1,0] \times B_1 \times \R^d$, eventhough the
equation needs to hold in $Q_1 = (-1,0] \times B_1 \times B_1$
only. This is a common theme in nonlocal equations. The operator
$\mathcal L$ takes into account the values of $f$ for every
$v' \in \R^d$. The first term in the right hand side of the inequality
allows us to control the tails of the integral in the expression
\eqref{e:integro-differential_operator} for $\cL_K$.

The symmetry condition \eqref{e:symmetry_cancellation} together with
the upper bound on the kernel \eqref{e:ellipticity_upper} allow us to
properly understand the weak solutions of the equation in the sense of
distributions via Proposition \ref{p:Hs2H-s}. In this case, it would
not be a problem to state Theorem \ref{t:whi} for bounded weak
solutions $f$.

Theorem \ref{t:holder} should be considered as a kinetic,
integro-differential, version of the well known theorem by De Giorgi,
Nash and Moser. We can also say that it is an integro-differential
version of the more recent estimate in \cite{gimv}. The cancellation
condition \eqref{e:symmetry_cancellation} is a weaker form of the
symmetry condition that corresponds to equations in divergence form.

Note that no smoothness assumption is imposed on either the kernel $K$
or the source function $h$. The H\"older estimate for $f$ depends only
on quantitative conditions on $K$ and $h$. As we discussed above (see
Section \ref{s:warm_up}), the kernel $K_f$ of the Boltzmann equation
given in \eqref{e:Kf_exact} automatically satisfies
\eqref{e:ellipticity_upper}, \eqref{e:ellipticity_cone} and
\eqref{e:symmetry_cancellation} for $v \in B_R$, for any bounded value
of $R$, with parameters $\mu$, $\lambda$ and $\Lambda$ depending on
the constants in \eqref{e:hydro}. As a consequence, Theorem
\ref{t:holder} implies that any solution of the Boltzmann equation
that satisfies the hydrodynamic condition \eqref{e:hydro} will be
locally H\"older continuous. We need further work (described in Section
\ref{s:change_of_vars}) in order to obtain a H\"older estimate that
holds uniformly for large velocities.

The rest of Section \ref{s:DG} is devoted to outline the steps of the
proof of Theorem \ref{t:holder}. It is somewhat technical and it may
be difficult to follow for the readers that are unfamiliar with (at
least) De Giorgi's method for parabolic equations
\cite{zbMATH03277871}. A reader that is willing to take Theorem
\ref{t:holder} for granted, can safely skip to
Section~\ref{s:schauder} at this point.

\subsection{The weak Harnack inequality}

Theorem \ref{t:holder} is a consequence of the following weak Harnack
inequality.
% -----------------------------------
\begin{thm}[Weak Harnack inequality]
  \label{t:whi}
  There exist radii $0 < r_0 < 1 < R_0$, only depending on dimension and $s$,
  and two positive constants $\eps$ (small) and $C$ (large), only
  depending on dimension, $s$ and ellipticity constants
  $\mu, \lambda,\Lambda$, such that any non-negative super-solution $f$
  of \eqref{e:lin-kin-2} in
  $\Qext := (-1,0] \times B_{R_0^{1+2s}} \times B_{R_0}$,
  \[
      \partial_t f + v \cdot \nabla_x f \ge \cL_K f + h \qquad \text{ in } \Qext,
  \]
  where $K$
  satisfies \eqref{e:ellipticity_upper}, \eqref{e:ellipticity_cone}
  and \eqref{e:symmetry_cancellation} and $h \in L^\infty (\Qext)$,
  satisfies
  \[
    \left( \int_{Q^-} f^\eps (z) \dz  \right)^{\frac1\eps} \le \Cwhi \left( \inf_{Q^+} f + \| h\|_{L^\infty(\Qext)} \right)
  \]
  where $Q^+= (-r_0^{2s},0] \times B_{r_0^{1+2s}} \times B_{r_0}$ and $Q^- = (-1,-1+r_0^{2s}] \times B_{r_0^{1+2s}} \times B_{r_0}$.
\end{thm}
% -----------------------------------

\begin{figure}
\setlength{\unitlength}{1in}
\begin{center}
\begin{picture}(2.903 ,1.000)
\put(0,0){\includegraphics[height=1.000in]{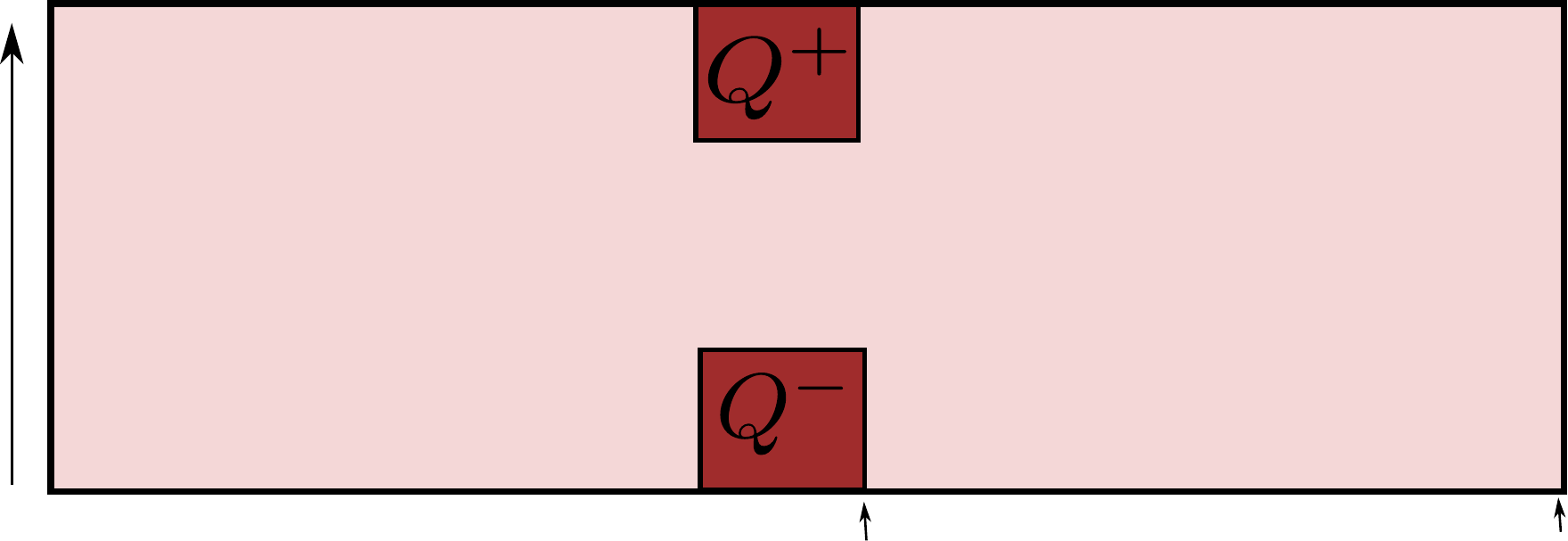}}
\put(1.57,-0.083){$r_0$}
\put(2.8,-0.102){$R_0$}
\put(-0.65,0.479){$t \in (-1,0]$}
\end{picture}
\end{center}
\caption{The geometric setting of the weak Harnack inequality}
\label{fig:whi}
\end{figure}
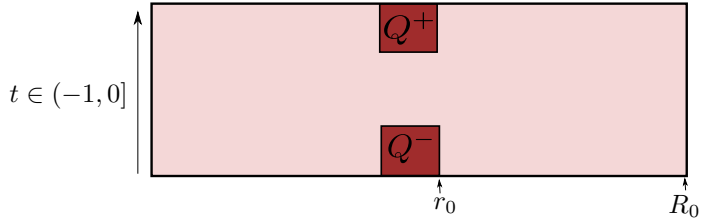

It is well known, in different contexts, that the weak Harnack
inequality implies the H\"older estimate of Theorem
\ref{t:holder}. The weak Harnack inequality is applied iteratively to
obtain a decay of the oscillation of the function at a sequence of
scales around a point. In this case, we have to adapt this classical
procedure to the context of kinetic cylinders, and accounting for the
nonlocality of the diffusion operator $\cL_K$. There is no major new
obstruction to prove that Theorem \ref{t:whi} implies Theorem
\ref{t:holder}. The difficult part is to prove Theorem \ref{t:whi}.

De Giorgi's method has two main parts. The first one is typically
presented as a control of the $L^\infty$ norm of a solution to an
elliptic or parabolic PDE in terms of its $L^2$ norm. The second part
is an iterative improvement of oscillation leading to the H\"older
continuity of that solution. Both parts will look very different in
our context. One can still argue that the first part of our proof is
inspired by the first part of De Giorgi's method. The second part of
our proof, however, uses ideas from the proof of Krylov and Safonov
for parabolic equations in non-divergence form \cite{ks}.

\subsection{The first lemma of De Giorgi}

The very first ingredient that one needs for De Giorgi's proof is a
relation known as \emph{Cacciopoli's inequality}. It is essentially a
localized version of the energy dissipation inequality that applies to
nonnegative sub-solutions of the equation.

Let us start by analyzing the norms of $f$ that we can control by the
energy dissipation. If we multiply the equation \eqref{e:lin-kin-2} by
$f$ and integrate using the coercivity property of Proposition
\ref{p:coercivity}, we get
\begin{equation} \label{e:energy_dissipation}
 \sup_{t \in [0,T]} \|f(t,\cdot, \cdot)\|_{L^2}^2 + \int_0^T \|f(t,\cdot,\cdot)\|_{L^2_x \dot H^s_v}^2 \leq \|f_0\|_{L^2}^2 + \langle \text{lower order terms} \rangle.
\end{equation}
We were intentionally vague about the domains of the norms. They
depend on the domain of the equation. Naturally, there are also some
boundary terms involved when the equation \eqref{e:lin-kin-2} holds in
a bounded domain. These boundary terms may be a complicated due to the
nonlocality of the diffusion operator $\cL_K$. Let us ignore them at
this point.

The energy dissipation estimate \eqref{e:energy_dissipation} provides
us with an improvement of differentiability with respect to the
$v$-variable. Indeed, it involves the $\dot H^s$ norm in this
variable. In order to carry out De Giorgi's method, we will need to
turn the estimate \eqref{e:energy_dissipation} into an estimate for
$\|f\|_{L^p_{t,x,v}}$ for some $p>2$. The estimate
\eqref{e:energy_dissipation} does not suffice for that because it does
not encode any regularization with respect to the variable $x$.

Now, we need to invoke the hypoelliptic nature of the equation
\eqref{e:lin-kin-2}. A standard idea in kinetic equations would be to
apply averaging lemmas to achieve this (as in \cite{stokols2019} or \cite{gimv}). We
take a more classical approach, which is inspired by
\cite{pascucci2004}. We plug the function $f$ into the fractional
Kolmogorov equation \eqref{e:fractional_kolmogorov}.
\begin{equation} \label{e:forced_kolmogorov}
 \partial_t f + v \cdot \nabla_x f + (-\Delta)^s_v f = (-\Delta)^s_v f + \cL_K f + h.
\end{equation}
Admittedly, at first sight this idea looks artificial. There is no
cancellation in the right hand side. The gain comes from the
observation that the right hand side involves fractional
differentiation with respect to the $v$ variable only, and these are
the directions where we got some regularity estimates from the energy
dissipation \eqref{e:energy_dissipation}. We can derive various
estimates thanks to our precise knowledge of the fractional Kolmogorov
equation. It has a semi-explicit kernel from which we can compute
exactly its regularization effects in all directions.

From the inequality \eqref{e:energy_dissipation} we get an estimate
for $\|f\|_{L^2_{t,x}(H^s_v)}$. Combining this estimate with
Proposition~\ref{p:Hs2H-s}, we deduce that the right hand side in
\eqref{e:forced_kolmogorov} is in $L^2_{t,x}(H^{-s}_v)$. Then we can
estimate the $L^p_{t,x,v}$ norm of the convolution of this right hand
side with the kernel of the fractional Kolmogorov equation. The
fractional Kolmogorov equation is invariant by the same group of
scaling and Galilean transformations as our class of kinetic nonlocal
equations. We know an explicit formula for its solution. This formula
encodes the hypoelliptic interaction between the kinetic transport
terms and the fractional diffusion in the velocity variable.

Ultimately, the computation described above leads to an estimate of the $L^p_{t,x,v}$ norm of $f$, for some $p>2$, in terms of its $L^2_{t,x,v}$ norm. Applying this estimate to proper truncations of the solution $f$ and following De Giorgi's iteration leads to the following version of De Giorgi's first lemma.

\begin{lemma} \label{l:weaker_harnack}
Let $f : [-1,0] \times B_1 \times \R^d \to [0,\infty)$ be a super-solution of the equation in $Q_1$,
\[ \partial_t f + v \cdot \nabla_x f - \cL_K f \geq 0 \qquad \text{in $Q_1$}.\]
There exists an $\eps_0>0$ (depending only on dimension and the ellipticity parameters) so that if
\[ | \{ f < 2\} \cap Q_1| < \eps_0,\]
then $f \geq 1$ in $Q_{1/2}$.
\end{lemma}

Lemma \ref{l:weaker_harnack} is a simplified version of \cite[Lemma
6.6]{imbert2020weak}.

Lemma \ref{l:weaker_harnack} is a lower bound for nonnegative
super-solutions of the equation. It differs from the classical
presentation of De Giorgi's first lemma as an upper bound for
sub-solutions. It would be possible to write an upper bound for
sub-solutions under the stronger assumptions $K \approx |v-v'|^{-d-2s}$
(see \cite{stokols2019}), but it is impossible under our less
restrictive hypothesis on the kernel. This has already been observed
in the context of parabolic integro-differential equations with
degenerate kernels (see \cite{Dyda-Kassmann-2015}).

Let us compare Lemma \ref{l:weaker_harnack} with Theorem
\ref{t:whi}. Their geometric settings are of course different, but we
can see some similarity when we state them in the following way.
\begin{itemize}
\item Lemma \ref{l:weaker_harnack} says that if $f \leq 1$ at any
  point in $Q_{1/2}$, then
  $|\{f \geq 2\} \cap Q_1| \leq |Q_1| - \eps_0$.
\item Theorem \ref{t:whi} (with $h=0$) says that if $f \leq 1$ at any
  point in $Q^+$, then $|\{f \geq A\} \cap Q^-| \lesssim A^{-1/\eps}$
  for all $A>0$.
\end{itemize}

Theorem \ref{t:whi} is effectively an upper bound for the measure of
the level sets $|\{ f \geq A\}| \cap Q^-$, for $A$ large, for any
super-solution of the equation so that $f \leq 1$ at some point in
$Q^-$.

The second part of De Giorgi's proof consists of an estimate of the
measure between two level sets of the solution $f$. Ultimately, it
leads to a decay estimate for the measure $|\{f \geq A\}|$ for large
$A$'s. The classical method by De Giorgi involves an explicit
computation relating the measure of level sets of a function with its
$H^1$ norm. In \cite{chan2011}, the authors follow an alternative
idea for integro-differential equations using an estimate
depending essentially on a lower bound on the kernel
$K \gtrsim |v'-v|^{-d-2s}$. In \cite{gimv}, De Giorgi's original
computation is replaced with an elegant compactness argument. None of
those ideas apply in our context. The compactness argument in
\cite{gimv} can be applied to integro-differential equations, after
working out several technical difficulties, but only in the case
$s \geq 1/2$. For $s<1/2$, we employ a completely different approach
inspired by the ideas of Krylov and Safonov in \cite{ks}. Our method
for $s<1/2$ applies in the full range $s \in (0,1)$ if we assume in
addition the non-divergence symmetry condition
\eqref{e:symmetry_non-divergence}.

The Boltzmann kernel $K_f$ defined in \eqref{e:Kf_exact} always
satisfies \eqref{e:symmetry_non-divergence}. Thus, the method to prove
Theorem~\ref{t:whi} inspired by the ideas of Krylov and Safonov
suffices for the whole range of parameters. In \cite{imbert2020weak},
we also describe the method inspired by the ideas in \cite{gimv}, that
works only for $s \geq 1/2$, because it allows us to remove the
assumption \eqref{e:symmetry_non-divergence} in the statement of
Theorem \ref{t:whi} for general kinetic integro-differential equations.

\subsection{The propagation lemma}

Lemma \ref{l:weaker_harnack} involves a lower bound in the small
kinetic cylinder $Q_{1/2}$. One can effortlessly scale Lemma
\ref{l:weaker_harnack} to relate the level set
$\{f \geq 2\} \cap Q_r(z)$ with the minimum of $f$ in $Q_{r/2}(z)$,
provided that the equation holds in $Q_r(z)$ and $f$ is nonnegative
everywhere. Our next objective is to extend the set where we take the
minimum of $f$ to a larger kinetic cylinder than $Q_{r/2}(z)$. This is
achieved through the use of explicit barrier functions described in
the following lemma.
%--------------------------------------------------------------------
\begin{lemma}[barrier functions] \label{l:barriers} Let $\tau>0$, $R>1$, and $T>0$ be arbitrary parameters. There exist $\theta>0$
  and $\Rone>R>0$ depending on these parameters, dimension, $s$ and the
  ellipticity constants in \eqref{e:ellipticity_upper},
  \eqref{e:ellipticity_cone} so that the following statement is true.

There exists a function $\varphi: [0,\infty) \times \R^d \times \R^d \to [0,1]$ satisfying the following properties.
\begin{itemize}
\item $\varphi \in C^{1,1}$. Moreover, $\varphi$ is $C^\infty$ in the set $\{\varphi > 0\}$.
\item For any kernel $K$ that satisfies \eqref{e:ellipticity_upper} and \eqref{e:ellipticity_cone} (and also \eqref{e:symmetry_non-divergence} in the case $s \geq 1/2$), we have
\[ \partial_t \varphi + v \cdot \nabla_x \varphi - \cL_K \varphi \leq 0 \text{ in } [0,\infty) \times \R^d \times \R^d.\]
\item $\varphi(0,x,v) > 0$ only if $(x,v) \in B_1 \times B_1$.
\item $\varphi(t,x,v) \geq \theta$ if $(t,x,v) \in [\tau,T] \times B_{R^{1+2s}} \times B_R$.
\item $\varphi(t,x,v) = 0$ if $t \in [0,T]$ and $(x,v) \notin B_{\Rone^{1+2s}} \times B_{\Rone}$.
\end{itemize}
\end{lemma}
%-------------------------------------------------------------------

The function $\varphi$ is used as a lower barrier. It allows us to
propagate a lower bound for a super-solution $f$ on
$\{0\} \times B_{r^{1+2s}} \times B_r$ to an arbitrarily large kinetic
cylinder $[\tau,T] \times B_{(Rr)^{1+2s}} \times B_{Rr}$, provided that the
equation holds in a suitable larger domain containing
$[0,T] \times B_{(\Rone r)^{1+2s}} \times B_{\Rone r}$.

The proof of Lemma \ref{l:barriers} consists of a more or less
explicit computation. It is another proof where the hypoellipticity of
the equation plays a role. This time, in a more crude and explicit
manner. The computation leading to Lemma \ref{l:barriers} is explained
in \cite[Section 7]{imbert2020weak}.

By scaling and translating this construction, we derive the following corollary.

\begin{cor}
  \label{c:propagation_of_lower_bounds}
  Given $R_0>0$ and $r_0>0$ such that $R_0 - r_0 \geq R_1 r_0$, let $\Qext$ and $Q^-$ be as in Theorem~\ref{t:whi}, see Figure~\ref{fig:whi}.  Let $f$ be a super-solution
  \[
    \partial_t f + v \cdot \nabla_x f - \cL_K f \geq 0 \text{ in }  \Qext.
  \]
  Assume that $K$ satisfies \eqref{e:ellipticity_upper} and   \eqref{e:ellipticity_cone}. We also assume  \eqref{e:symmetry_non-divergence} when $s \geq 1/2$. If $f \geq A$
  in $Q_r(z_0)$ for some cylinder $Q_r(z_0)$ contained in $Q^-$, then for any arbitrarily large constants $T$ and $R$ as in Lemma~\ref{l:barriers}, there exists $\theta >0$ such that
  \[
    f \geq \theta A \quad \text{in } \quad  [t_0,\min(t_0 + Tr^{2s},0)] \times B_{(Rr)^{1+2s}} \times B_{Rr} .
\]
Moreover, the factor $\theta>0$ depends on $T$, $R$, $d$, $s$ and the ellipticity parameters of the kernel $K$.
\end{cor}

We can further combine Corollary \ref{c:propagation_of_lower_bounds}
with Lemma \ref{l:weaker_harnack} and get
\begin{cor}
  \label{c:propagated_weaker_harnack}
  Let $f$ be a super-solution
  \[
    \partial_t f + v \cdot \nabla_x f - \cL_K f \geq 0 \text{ in }  \Qext.
  \]
  Assume that $K$ satisfies \eqref{e:ellipticity_upper} and
  \eqref{e:ellipticity_cone}. We also assume
  \eqref{e:symmetry_non-divergence} when $s \geq 1/2$. Assume also
  that for some cylinder $Q_r(z_0)$ contained in $Q^-$,
  \[
    |\{ f \geq A\} \cap Q_r(z_0)| \geq (1-\eps_0) |Q_r|.
  \]
  Then
  \[
    f(t,x,v) \geq \theta A \quad \text{in} \quad [t_0,t_0 + Tr^{2s}] \times B_{(Rr)^{1+2s}} \times B_{Rr}.
  \]
  Here, $\theta$, $T$, $R$, $\Qext$ and $Q^-$ are as in Corollary \ref{c:propagation_of_lower_bounds}.
\end{cor}

Corollary \ref{c:propagated_weaker_harnack} tells us a nontrivial
relation between the level sets $\{f \geq A\}$ and
$\{f \geq \theta A\}$. It implies Theorem \ref{t:whi}. However, this
implication is nontrivial. It depends on a special covering argument
inspired in the crawling ink-spots lemma from \cite{ks}, that we describe below.

Here is a further corollary that will be used in the proof of Theorem
\ref{t:whi}.

\begin{cor} \label{c:inkspots_min_radius} Let a kernel $K$ satisfy
  \eqref{e:ellipticity_upper} and \eqref{e:ellipticity_cone} and, when
  $s \geq 1/2$, let it also satisfy \eqref{e:symmetry_non-divergence}.
  There exists a constant $\rzero>0$, that only depends on dimension
  $d$, $s$ and the ellipticity parameters of $K$, such that for all
  super-solution $f$ of
\[
  \partial_t f + v \cdot \nabla_x f - \cL_K f \geq 0 \text{ in }  \Qext
\]
such that $\min_{Q^+} f \leq 1$, if there exists some cylinder
$Q_r(z_0)$ contained in $Q^-$ such that,
\[
  |\{ f \geq A\} \cap Q_r(z_0)| \geq (1-\eps_0) |Q_r|,
\]
then $r < \rzero$.
\end{cor}

Corollary \ref{c:inkspots_min_radius} is an immediate consequence of
Corollary \ref{c:propagated_weaker_harnack}. Indeed, if $r$ was large,
Corollary \ref{c:propagated_weaker_harnack} would imply that $f > 1$
in $Q^+$.

\subsection{Ink-spots}

The following lemma is a type of covering result known as the
ink-spots lemma.

\begin{lemma}[Ink-spots]
  \label{l:inkspots}
  Let $E \subset F \subset B_1$
  be two measurable sets. Assume that for some constant $\delta>0$,
  $|E| \leq (1-\delta)|B_1|$ and whenever there is any ball
  $B \subset B_1$ such that $|E \cap B| \geq (1-\delta) |B|$, we must
  have $B \subset F$. Then, the following inequality holds
  $|E| \leq (1-c\delta)|F|$ for some constant $c>0$ depending on
  dimension only.
\end{lemma}

\begin{figure}[htb]
  \includegraphics[height=1.5in]{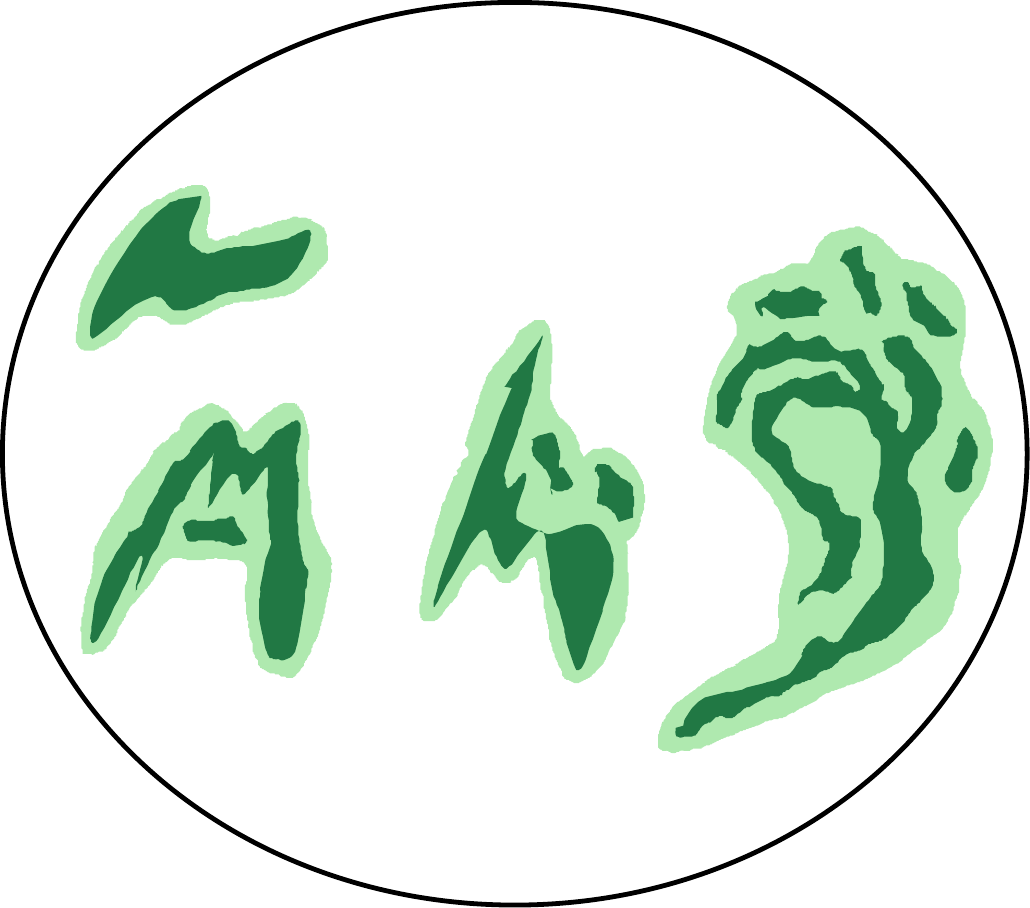}
  \caption{The set $F$ (greenish color) looks like the result of a growing ink stain which started as $E$ (dark color).}
  \label{f:inkspots}
\end{figure}

Lemma \ref{l:inkspots} is used to obtain the decay of level sets
in the elliptic version of the weak Harnack inequality by Krylov and
Safonov. It is relatively easy to prove Lemma \ref{l:inkspots}
as a consequence of Vitali's covering lemma.

In \cite{ks}, Krylov and Safonov describe a modification of
Lemma \ref{l:inkspots} that is suitable for parabolic
equations. They call it \emph{the crawling ink-spots lemma}. In that
case $E \subset F$ are sets in space-time. They assume that whenever
there is a parabolic cylinder where $E$ is very concentrated, then an
enlarged version of that cylinder, that takes place later in time, is
contained in $F$.

We need to further modify the covering lemma in \cite{ks} to fit the
setting of kinetic equations. In order to state our kinetic version of
the ink-spots lemma, we start with defining the \emph{stacked}
cylinder. Given any kinetic cylinder $Q=Q_r(z_0)$ (as in Section
\eqref{s:cylinders}), we define $Q^m$ as
\[ Q^m := \{ (t,x,v) : 0 < t-t_0 < m r^{2s} , |v-v_0| < r,
  |x-x_0-(t-t_0)v_0| < (m+2) r^{1+2s} \}.\] Note that $t_0$ is the
final time for $Q$ and the initial time for $Q^m$. The lapse of $Q^m$
is $m$ times the lapse of $Q$. Moreover, the space width of $Q^m$ is
enlarged by a factor $(m+2)$ with respect to the space width of $Q$.

The following result is our kinetic version of the ink-spots lemma. It
is proved in \cite[Section 10]{imbert2020weak}.

\begin{lemma}[Kinetic crawling ink-spots]
  \label{l:inkspots-kinetic}
  Let $E \subset F$ be two measurable sets. We make the following
  assumptions for some $\delta \in (0,1)$ and some $\rzero>0$.
  \begin{itemize}
  \item $E \subset Q_1$.
  \item Whenever a kinetic cylinder $Q \subset Q_1$ satisfies
    $|Q \cap E| \geq (1-\delta) |Q|$, then $Q^m \subset F$ and also
    $Q = Q_r(z)$ for some $r<\rzero$.
  \end{itemize}
  Then
  \[ |E| \leq \frac{m+1}m (1-c\delta) \left( |F \cap Q_1| + C m \rzero^{2s} \right).\]
\end{lemma}

With these ingredients, we can outline the proof of Theorem~\ref{t:whi}.

\begin{proof}[Sketch of proof of Theorem \ref{t:whi} for $h=0$]
Assume $f(z) \leq 1$ for some point $z \in Q^+$. We need to prove that
\[ |\{ f \geq A\} \cap Q^-| \lesssim A^{-1/\eps}.\]
This decay follows by a simple iteration once we established the inequality
\[
  |\{ f \geq A\} \cap Q^-| \leq (1-c \eps_0) |\{ f \geq \theta A\}  \cap Q^-|,
\]
for $\theta A > 2$, $\theta$ and $\eps_0$ as in
Corollary \ref{c:propagated_weaker_harnack}, and $c$ as in Lemma
\ref{l:inkspots-kinetic}.

This inequality between level sets follows from Lemma
\ref{l:inkspots-kinetic} applied to
$E = \{ f \geq \theta A\} \cap Q^-$ and $F = \{ f \geq A\} \cap
Q^-$. The assumptions of Lemma \ref{l:inkspots-kinetic} are fulfilled
thanks to Corollaries \ref{c:propagated_weaker_harnack} and
\ref{c:inkspots_min_radius}.
\end{proof}

\section{The Schauder estimate}
\label{s:schauder}

This section is devoted to the Schauder theory for kinetic equations
in non-divergence form with H\"older continuous coefficients. We
consider a function $f$ that solves \eqref{e:lin-kin-2} with a
H\"older continuous source term $h$ and a kernel $K$ that satisfies
the ellipticity conditions \eqref{e:ellipticity_upper} and
\eqref{e:ellipticity_cone} together with the \emph{non-divergence}
symmetry condition \eqref{e:symmetry_non-divergence}.

We must first make sense of the notion of \emph{H\"older continuous coefficients} for a kernel $K(t,x,v,v')$. We add the following assumption, which depends on a parameter $\alpha' \in (0,\min(1,2s))$.

\begin{assumption}[H\"older continuity of the kernel in $(t,x,v)$]
  \label{a:coef}
  There exists a positive constant $A_0$ such that whenever
  $z_1 =(t_1,x_1,v_1)$ and $z_2 = (t_2,x_2,v_2)$ belong to $Q_1 \cap Q_r(z_0)$, for any kinetic cylinder $Q_r(z_0)$, then
  \[
    \forall \rho>0, \qquad \int_{B_\rho} |K (t_1,x_1,v_1,v_1+w) - K (t_2,x_2,v_2,v_2+w)| |w|^2 \d w \le A_0 \rho^{2-2s} r^{\alpha'}.
  \]
\end{assumption}
%-----------------------------------------------------------------------

With this notion of $C^{\alpha'}$ \emph{H\"older coefficients}, we are ready to state the Schauder estimates for general kinetic integro-differential equations.

\begin{thm}[The Schauder estimate]
  \label{t:schauder}
  Let $s \in (0,1)$, $\alpha \in (0, \min (1,2s))$ and $\alpha' = 2s \alpha / (1+2s)$. Let
  $K\colon Q_1 \times \R^d \to \R$ be a nonnegative kernel such that
  \eqref{e:ellipticity_upper}, \eqref{e:ellipticity_cone},
  \eqref{e:symmetry_non-divergence} and Assumption~\ref{a:coef} hold
  true.  Let $h: Q_1 \to \R$ be $\alpha'$-H\"older continuous. Assume further that $2s+\alpha' \notin \{1,2\}$.

  If $f$ satisfies
  \eqref{e:lin-kin-2} in $Q_1$, then
  \[
    [f]_{\cCl^{2s+\alpha'} (Q_{1/2})} \le C (\| f \|_{\cCl^\alpha ((-1,0] \times B_1 \times \R^d)} + \|h\|_{\cCl^{\alpha'} (Q_1)}) .
  \]
  The constant $C$ only depends on dimension, the order $2s$ of the
  integral diffusion, ellipticity constants $\mu,\lambda,\Lambda$ and
  $A_0$ from Assumption~\ref{a:coef}.
\end{thm}
%------------------------------------

A slightly more general version of Theorem \ref{t:schauder} is obtained in \cite{imbert2016schauder}.

Since we do not assume the cancellation condition \eqref{e:symmetry_cancellation}, the notion of weak solutions in the sense of distributions does not make sense in the generality of
Theorem \ref{t:schauder}. This is the same situation as in the classical Schauder estimates for elliptic PDEs in non-divergence form. The natural framework under which the Schauder estimates apply
is that of viscosity solutions. Such a generalization would involve only some minor technical adjustments to our current proof. Likewise, if we assume in addition that the cancellation condition \eqref{e:symmetry_cancellation} holds (which is true in the case of the Boltzmann equation), then Theorem~\ref{t:schauder} would extend to weak solutions in the sense of distributions without any major additional difficulty.

To prove Theorem \ref{t:schauder}, we first analyze the simpler case
where $K$ depends only on $(v'-v)$. Then, we apply an interpolation
inequality to account for the variations of the kernel, similarly as
in the classical proof of the Schauder estimates. Note that if $K$
depends only on $(v'-v)$, Assumption \ref{a:coef} automatically holds
with $A_0 = 0$.

The case $K= K(v'-v)$ is proved using a blow-up technique, following
closely the ideas by  Serra \cite{serra2015}. We set up the
requirements for an iterative proof of the regularity result. We
proceed by contradiction by negating the main estimate. A
\emph{blow-up limit} leads to certain ancient solution of the
equation. The contradiction is reached through a Liouville type result
that rules out such a solution.

In the rest of this section, we outline the ideas involved in the
proof of Theorem~\ref{t:schauder}. The uninterested reader may move directly to
Section \ref{s:upper_bounds}.

\subsection{A Liouville type result}

A simple form of Liouville theorem that one can state for kinetic
integro-differential equations is the following: every bounded
solution of \eqref{e:lin-kin-2}, with $h=0$, in
$(-\infty,0] \times \R^d \times \R^d$, under the assumptions of
Theorem~\ref{t:holder}, must be constant. The proof of this statement
is simply to apply Theorem~\ref{t:holder} in $Q_R$ for large $R$'s. We
obtain that the H\"older seminorm of $f$ in $Q_{R/2}$ is
$\lesssim R^{-\alpha}$. We deduce that $f$ is constant taking
$R \to \infty$. A slightly more detailed analysis of the inequalities
of this proof reveals that the \emph{boundedness} hypothesis can be
relaxed to a slow enough algebraic rate. Let us state it in the
following proposition.
%--------------------------------------------------------
\begin{proposition}[Liouville -- I]  \label{p:liouville_easy} Let $f$ be a solution of
  \eqref{e:lin-kin-2}, with $h=0$, in
  $Q_\infty = (-\infty,0] \times \R^d \times \R^d$, where $K$
  satisfies \eqref{e:ellipticity_upper}, \eqref{e:ellipticity_cone}
  and \eqref{e:symmetry_cancellation}. Assume further that
  \[ \| f\|_{C^0(Q_R)} \leq C(1+R)^\delta,\] for some constant $C$ and
  some $\delta \geq 0$ smaller than the H\"older exponent $\alpha$ in
  Theorem \ref{t:holder}. Then $f$ is constant.
\end{proposition}

For the proof of Theorem \ref{t:schauder}, we need to allow more
growth at infinity for the solution $f$. The precise form of the
Liouville theorem that we use is the following. It applies to a
kernel $K=K(v'-v)$ that depends on $(v'-v)$ only.

%----------------------------------------------------------------------
\begin{proposition}[Liouville -- II]
  \label{p:liouville}
  Let $\alpha \in (0,\min (1,2s))$ and $\alpha' = \frac{2s}{1+2s}\alpha$
  and $\beta$ such that
  $\lfloor 2s + \alpha' \rfloor < 2s + \beta < 2s + \alpha'$ and   $\alpha' - \beta < \delta$, where $\delta>0$ is the H\"older exponent in Theorem \ref{t:holder}. Let $f$ be a function
  in $\cClloc^{2s+\beta} ((-\infty,0] \times \R^d \times \R^d)$  satisfying the following conditions.
  \begin{enumerate}
  \item \label{a:i} There exists a constant $C_1>0$ such that for all $R>0$,
    \[
      \forall \beta' \in [0,2s+\beta], \quad [f]_{\cCl^{\beta'}(Q_R)} \le C_1 R^{2s+ \alpha'- \beta'};
    \]
  \item \label{a:ii} For any $\xi = (h,y,w) \in \R^{1+2d}$ with $h \le 0$, the function $g (z):= f(\xi \circ z)-f(z)$ solves
    \[
      \partial_t g + v \cdot \nabla_x g - \cL_K g = 0
      \qquad \text{ in } (-\infty,0] \times \R^d \times \R^d
    \]
    where $K=K(v-v')$ satisfies \classK.
  \end{enumerate}
  Then $f$ is a polynomial of kinetic degree smaller than $2s+\alpha'$.
\end{proposition}
%----------------------------------------------------------------------

One advantage of kernels that depend on $(v'-v)$ only is that they are in divergence and non-divergence form at the same time. Indeed, if
$K=K(v'-v)$ satisfies \eqref{e:symmetry_non-divergence}, then it also
satisfies \eqref{e:symmetry_cancellation} with $\Lambda=0$. In this
case, $\cL_K$ is comparable to a second order elliptic operator with
constant coefficients.

Note that in the statement of Proposition \ref{p:liouville}, we do not
require $f$ to be the solution of an equation. The equation is imposed
on its increments $g$. The growth at infinity for the function $f$
would make the nonlocal operator $\cL_K f$ undefined. Indeed, with the
assumptions of Proposition \ref{p:liouville}, for any fixed value of
$t_0,x_0, v_0$, we only have
$|f(t_0,x_0,v') - f(t_0,x_0,v_0)| \leq C_{t_0,x_0,v_0}
(1+|v'|)^{2s+\alpha'}$. With this growth at infinity, the tail of the
integral in \eqref{e:integro-differential_operator} diverges. We state
the equation for $g$, and not for $f$, because the equation for $f$
may not make sense.

The proof of Proposition \ref{p:liouville} consists in applying
Proposition \ref{p:liouville_easy} to various increments of $f$. The
first step is to take increments in space $\xi=(0,y,0)$. These values
of $\xi$ are in the center of the Galilean Lie group, so all the
computations work like one would expect. Applying
Proposition~\ref{p:liouville_easy} to $g$, with $\xi = (0,y,0)$, we
conclude that the function $f$ must be constant in $x$. Once we
established that $f$ must be constant in $x$, the statement of
Proposition \ref{p:liouville} reduces to a Liouville theorem for
parabolic integro-differential equations. The rest of the proof
continues by applying Proposition \ref{p:liouville_easy} to increments
of $f$ in $t$ and $v$. The details of the proof are given in
\cite[Section 4]{imbert2020schauder}.

\subsection{The blow-up argument}

In this paragraph, we briefly explain how to derive the Schauder
estimate from the Liouville type result in the ``constant coefficient
case'', that is to say when the kernel $K$ depends only on $(v'-v)$.

\begin{proof}[Sketch of the proof of Theorem~\ref{t:schauder} for  $K=K(v'-v)$]
Following Serra \cite{serra2015}, the core of the proof of the Schauder estimate consists in proving that for some $\beta$ smaller than $\alpha'$, the following estimate holds true,
  \begin{equation} \label{e:schauder_proof} [f]_{\cCl^{2s+\alpha'}
      (Q_{1/2})} \lesssim \left(\|h\|_{\cCl^{\alpha'} (Q_1)} +
      \|f\|_{\cCl^{2s+\beta} ((-1,0]\times B_1 \times \R^d)} \right).
  \end{equation}
Using the interpolation inequalities from Proposition~\ref{p:interpol} together with some control on the tails of the integral operator, yields the Schauder estimate as stated in Theorem~\ref{t:schauder}. Let us focus on the derivation of the previous estimate.

Since we assumed $2s+\alpha'$ is not an integer, let us pick $\beta$ so that $\lfloor 2s+\alpha'\rfloor \leq 2s+\beta \leq 2s+\alpha'$. Note that since $\alpha' = 2s \alpha / (1+2s)$ with $\alpha < \min(1,2s)$, we know further that $(\N + 2s \N) \cap [2s+\beta,2s+\alpha'] = \emptyset$.

By an appropriate normalization, we assume without loss of generality that
  \[\| h \|_{\cCl^{\alpha'} (Q_1)} + \|f\|_{\cCl^{2s+\beta} ((-1,0]\times B_1 \times \R^d)} \le 1,\] and we
aim at bounding $[f]_{\cCl^{2s+\alpha'} (Q_{1/2})} $ from above.

Since $(\N + 2s \N) \cap [2s+\beta,2s+\alpha') = \emptyset$, we are able to estimate the seminorm in $C^{2s+\alpha'}$ if we establish the following inequality
  \[
    [f]_{\cCl^{2s+\beta} (Q_r (z))} \le C_0 r^{\alpha'-\beta}
    \ \text{ for all } z \in Q_{1/2}, r>0 \text{ such that } Q_r(z) \subset Q_1.
  \]
  Indeed, this inequality implies
  $[f]_{\cCl^{2s+\alpha'}(Q_{1/2})} \lesssim C_0$.

  Since we work with a smooth \footnote{The qualitative
    smoothness assumption on $f$ simplifies the proof but is not
    strictly necessary.}  solution $f$, we know that the following maximum is
  achieved at some point $z \in \overline Q_{1/2}$ and
  $r \in (0,1/2]$.
  \[
    \max_{r>0, z \in \overline{Q_{1/2}}, Q_r(z) \subset Q_1 } r^{\beta-\alpha'} [f]_{\cCl^{2s+\beta} (Q_r (z))} =: C_0.
  \]
  We must prove that $C_0$ is bounded from above in terms of $d$, $s$ and the ellipticity parameters only.

  The proof proceeds by contradiction. We assume that there exist
  sequences $f_j\in \cCl^{2s+\beta} ((-1,0] \times B_1 \times \R^d)$,
  $h_j \in \cCl^{\alpha'} (Q_1)$, $K_j \in \cK$ such that
  \begin{eqnarray}
    \label{e:estim-cj}
    \|h_j\|_{\cCl^{\alpha'} (Q_1)} + \|f_j \|_{\cCl^{2s+\beta} ((-1,0] \times B_1 \times \R^d)} \le 1, \\
    \nonumber
    (\partial_t + v \cdot \nabla_x) f_j - \cL_{K_j} f_j= h_j, \\
    \nonumber
    \sup_{r>0, z \in Q_{1/2}, Q_r(z) \subset Q_1} r^{\beta-\alpha'} [f_j]_{\cCl^{2s+\beta} (Q_r (z))} \nearrow +\infty \quad \text{ as } j \to +\infty.
   \end{eqnarray}
   Above, we observed that, for each $j$, the supremum is reached at
   some $r_j > 0$ and $z_j \in \overline{Q_{1/2}}$. We can see that
   necessarily $r_j \to 0$ as $j \to +\infty$ in the above sequence.

   In order to derive a contradiction, we rescale the sequence of
   functions $\{ f_j \}$ to map its values in $Q_{r_j} (z_j)$ to the
   values of $\tilde f_j$ in $Q_1$, to normalize the
   $\cCl^{2s+\beta}$-norm of $\tilde f_j$ in $Q_1$, and by removing
   the polynomial part $q_j$ (of kinetic degree strictly smaller than
   $2s+\beta$).  Recalling that $S_r$ denotes the natural scaling
   associated with the class of equations we work with, let
   $\tilde{f}_j$ be defined for all $z \in Q_1$ by,
   \[ \tilde{f}_j (z) = \frac{ (f_j -q_j) (z_j \circ S_{r_j} (z))}{r_j^{2s+\alpha'} F_j},\]
where
\[     F_j := r_j^{\beta-\alpha'} [f_j]_{\cCl^{2s+\beta} (Q_{r_j} (z_j))}. \]
   It then satisfies $\|\tilde{f}_j\|_{\cCl^{2s+\beta} (Q_1)} =1$ and
   \begin{eqnarray*}
     \forall R \in [1,c_s r_j^{-1}], && [\tilde{f}_j]_{\cCl^0 (Q_R)} \le R^{2s+\alpha'}, \\
     \forall R\in [1,c_s r_j^{-1}], && [\tilde{f}_j]_{\cCl^{2s+\beta}(Q_R)} \le R^{2s+\alpha'-\beta}.
   \end{eqnarray*}
   The constant $c_s >0$ in the previous estimates only depends on
   $s$. It is chosen so that if $z_j \in Q_{1/2}$ then
   $Q_{c_s} (z_j) \subset Q_1$.

   The remainder of the proof consists in getting a limit $f_\infty$
   of the (sub)sequence of $\{\tilde f_j\}$ and to apply the Liouville
   theorem to the function $f_\infty$. We prove that the increments of
   $f_\infty$ satisfy an equation, and Proposition \ref{p:liouville}
   implies that the \emph{blow-up} limit $f_\infty$ has to be a
   polynomial of degree less than $2s+\alpha'$. However, the
   polynomial expansion at the origin is zero for each $\tilde f_j$ by
   construction. We conclude that $f_\infty$ has to vanish,
   contradicting $\|\tilde{f}_\infty\|_{\cCl^{2s+\beta} (Q_1)} =1$.
\end{proof}

\section{Pointwise upper bounds}
\label{s:upper_bounds}

In this section we describe the pointwise upper bounds for the
Boltzmann equation that we obtained in \cite{imbert2020decay}, and we
briefly described in Section \ref{s:warm-up_bounds}.
%--------------------------------------------------------------------
\begin{thm}[Pointwise upper bounds]
  \label{t:upper_bounds}
  Let $\gamma \in (-d,1]$, $s \in (0,1)$ such that
  $\gamma +2s \in [0,2]$ and let $B$ be a collision kernel of the
  non-cutoff form \eqref{e:B}. Let $f$ be a solution of the Boltzmann
  equation \eqref{e:boltzmann} in $(0,T) \times \T^d \times \R^d$ such
  that $f(0,x,v) = \fin (x,v)$ in $\T^d \times \R^d$ and
  \eqref{e:hydro} holds. We obtain the following estimates.
\begin{enumerate}
\item \textbf{Propagation of upper bounds.} There exists a decay rate $q_0$
  depending on $s$, $\gamma$, $d$ and the parameters in
  \eqref{e:hydro} so that for any $q > q_0$, if
  $f_0 \le C (1 + |v|)^{-q}$ for some $C>0$, then there exists a
  constant $N$ depending on $C$, $d$, $s$, $d$ and the parameters in
  \eqref{e:hydro} such that
  \[
    \forall \, t \in [0,T], \ x \in \mathbb{T}^d, \ v \in \R^d, \quad
    f (t,x,v) \le N \left(1+|v|\right)^{-q}.
  \]
\item \textbf{Generation of upper bounds for hard potentials.} Assume
  $\gamma>0$. For any $q \geq 0$, there exist constants $N$ and
  $\beta>0$, depending on $d$, $s$, $d$, \textcolor{red}{$q$} and the parameters in
  \eqref{e:hydro} only, such that
  \[
    \forall \, t \in [0,T], \ x \in \mathbb{T}^d, \ v \in \R^d, \quad
    f (t,x,v) \le N (1+t^{-\beta}) \left(1+|v|\right)^{-q}.
  \]
\item \textbf{Generation of upper bounds for soft potentials.} Assume
  $\gamma \in [-2s,0]$. There exists a constant $N$, depending on $d$,
  $s$, $d$ and the parameters in \eqref{e:hydro} only, such that
  \[
    \forall \, t \in [0,T], \ x \in \mathbb{T}^d, \ v \in \R^d, \quad
    f (t,x,v) \le N \left(1+t^{-\frac{d}{2s}}\right) \left( 1 +
      |v|\right)^{-d-1-\frac{d \gamma}{2s}}.
  \]
\end{enumerate}
\end{thm}
% --------------------------------------------------------------------

Unlike the H\"older estimates in Theorem \ref{t:holder}, or the
Schauder estimates in Theorem \ref{t:schauder}, the upper bounds in
Theorem \ref{t:upper_bounds} do not apply to generic
integro-differential equations. They are specific to the Boltzmann
equation \eqref{e:boltzmann} with the collision operator $Q(f,f)$ and the non-cutoff kernel of the form
\eqref{e:B}, for solutions satisfying the assumption \eqref{e:hydro}.

The proof of Theorem~\ref{t:upper_bounds} has the basic structure of a
classical barrier argument for parabolic PDEs. We postulate that
$f(t,x,v) \leq U(t,v) := N A(t) (1+|v|)^{-q}$. We pick $A(t)$ to be
constant for the proof of the propagation of the upper bounds, and to
be of the form $A(t) = (1+t^{-\beta})$ for the generation of upper
bounds. In either case, the inequality holds trivially at the initial
time. If the conclusion of the theorem was false, there would be a
first crossing point $(t_0,x_0,v_0)$ so that
\begin{align*}
f(t_0,x_0,v_0) &= U(t_0,v_0), \\
f(t,x,v) &\leq U(t,v) \text{ whenever } t \leq t_0.
\end{align*}
At this first crossing point, the conditions above imply the following relations for the derivatives of $f$
\begin{align*}
\partial_t f(t_0,x_0,v_0) &\geq \partial_t U(t_0,v_0), \\
\nabla_x f(t_0,x_0,v_0) &= \nabla_x U(t_0,v_0) = 0.
\end{align*}
Unlike the classical method of barrier functions, there is no straight forward relationship between $Q(f,f)$ and $Q(U,U)$ at the first
crossing point $(t_0,x_0,v_0)$. The inequality $\cL_{K_f} U \geq \cL_{K_f} f$ holds at $(t_0,x_0,v_0)$, but is not
enough for the proof of our result. Instead, we shall use the inequality: $f(t_0,x_0,v) \leq U(t_0,v)$, for all $v \in \R^d$,
together with the upper bounds on mass and energy in \eqref{e:hydro}, and perform a delicate analysis of the quadratic integral operator
$Q(f,f)$ to deduce a (negative) upper bound for $Q(f,f)(t_0,x_0,v_0)$. This is the key computation that leads to the proof of Theorem \ref{t:upper_bounds}, and it is purely nonlocal. In
particular, the same reasoning as in the proof of Theorem \ref{t:upper_bounds} does not work with the Landau equation (see \cite{cameron2018} for a result on upper bounds for the Landau
equation using different methods).

The precise computation for estimating the value of
$Q(f,f)(t_0,x_0,v_0)$ relies on the integro-differential structure of
the collision operator,
\[
Q(f,f)(t_0,x_0,v_0) = \mathcal{L}_{K_f} f (t_0,x_0,v_0) +  f (t_0,x_0,v_0)(f(t_0,x_0,\cdot) \ast c_b |\cdot|^\gamma)(v_0) .
\]
We take advantage of the cone of nondegeneracy
described in Subection~\ref{s:cone} and the upper bounds on mass and
energy in \eqref{e:hydro}. The integral structure
of $\mathcal L_{K_f}$ allows us to immediately transfer the
information given by the bounds $M_0$ and $E_0$ in \eqref{e:hydro} to
an estimate of the value of $\mathcal{L}_{K_f}$ at the point
$(t_0,x_0,v_0)$. There is no analog of this method for the usual,
\emph{local}, partial differential equations.

This is the only step in our program where we use the periodicity assumption of $f$ with respect to the $x$ variable. It is used for convenience and only to ensure the existence of the first crossing point $(t_0,x_0,v_0)$. There are many alternative structures for global solutions $f : [0,T] \times \R^d \times \R^d \to [0,\infty)$ that may also be considered. Note that our estimates do not depend on the length of the period.

From all the steps in our regularity program, Theorem
\ref{t:upper_bounds} is the only one which would be hard to generalize
to weak solutions. The idea of evaluating the equation at the first
crossing point is not compatible with the notion of solution in the
sense of distributions. The natural generalized framework for this
method is that of \emph{viscosity sub-solutions}, as developed by
Crandall and Lions for first order, elliptic and parabolic
equations. The notion of viscosity solutions have not attracted
practically any attention in the context of the Boltzmann equations to
date. The definition of viscosity super-solutions is given in
\cite{imbert2019lower}. The definition of viscosity sub-solution, is
identical to the one of super-solution but with reversed
inequalities. The definition would only make sense for a solution $f$
that is at least locally bounded, so it would still be quite
restrictive in terms of the qualitative properties of the solutions we
start with.

\section{The change of variables}
\label{s:change_of_vars}

After getting the upper bounds of Theorem \ref{t:upper_bounds}, the
H\"older estimate from Theorem \ref{t:holder} and the Schauder
estimates from Theorem \ref{t:schauder} imply local H\"older bounds
for the Boltzmann equation as soon as we verify that their hypothesis
hold for the Boltzmann kernel $K_f$ of \eqref{e:Kf_exact} in terms of
the parameters of \eqref{e:hydro} only.

Like we explained in Section \ref{s:cone}, the cone of nondegeneracy
in condition \eqref{e:ellipticity_cone} is automatically satisfied by
the kernel $K_f$ of the Boltzmann equation provided that we restrict
our attention to some bounded set of velocities $v \in B_R$. In that
case, we have that the set $A=A(v) \subset S^{d-1}$ is symmetric with
respect to the origin and concentrated in a band of width
$\lesssim (1+|v|)^{-1}$ around the equator perpendicular to $v$. Its
measure is $\approx(1+|v|)^{-1}$, and
\begin{equation} \label{e:cone_of_nondeg_cov}
 K_f(v,v') \geq \lambda (1+|v|)^{1+\gamma+2s} |v'-v|^{-d-2s}, \qquad \text{ whenever } \frac{(v'-v)}{|v'-v|} \in A.
\end{equation}
Naturally, the condition \eqref{e:ellipticity_cone} is satisfied for
$K_f$, but only if we restrict our analysis to some bounded set of
velocities $v \in B_R$. We make the same observation for the other
conditions \eqref{e:ellipticity_upper} and
\eqref{e:symmetry_cancellation}. For example, the following inequality
is proved in \cite{silvestre2016new}, which justifies that
\eqref{e:ellipticity_upper} holds for $v \in B_R$:
\[
  \int_{\R^d \setminus B_R} K(v,v') \d v' \lesssim R^{-2s} \left(
    \int_{\R^d} f(v+w) |w|^{\gamma+2s} \d w \right) \leq C R^{-2s}
  (1+|v|)^{\gamma+2s},
\]
for a constant $C$ depending only on $M_0$ and $E_0$ in
\eqref{e:hydro}, provided that $\gamma+2s \in [0,2]$.

A similar inequality holds for
\eqref{e:symmetry_cancellation}. Therefore, Theorem \ref{t:holder}
applies locally to any solution of the Boltzmann equation that
satisfies \eqref{e:hydro}. It tells us that any such solution is
H\"older continuous in $(\tau,\infty) \times \R^d \times B_R$, for any
positive time $\tau>0$ and any bounded value for $R$. Now we would
like to apply our H\"older estimate for $f$ to deduce Assumption
\ref{a:coef} for $K_f$, and in that way be able to apply the Schauder
estimates of Theorem~\ref{t:schauder}. Later on, we would want to keep
applying the Schauder estimates to increments and derivatives of $f$
to deduce higher regularity. However, in order to do so, it is
essential to obtain global H\"older and Schauder estimates, that are
not restricted to bounded velocities only. That difficulty is overcome
by the change of variables described in this section.

For any value of $v_0 \notin B_1$, we define the linear transformation
$T_0$ by the formula
\[
  T_{v_0}(a v_0 + w) := \frac{a}{|v_0|} v_0 + w \qquad \text{  whenever } w \perp v_0.
\]
Note that $T_0$ maps the unit ball $B_1$ into an ellipsoid that is flattened by the factor $1/|v_0|$ in the direction of $v_0$. If $v_0 \in B_1$, we simply take $T_{v_0}$ to
be the identity. Given any $z_0 = (t_0,x_0,v_0)$, we further define
\[
  \mathcal T_{z_0} := z_0 \circ \left(|v_0|^{-\gamma-2s} t,
    |v_0|^{-\gamma-2s} T_{v_0} x, T_{v_0} v \right).
\]
Here, $\circ$ is the Galilean group operator in $\R^{1+2d}$. This transformation $\mathcal T_{z_0}$ maps $Q_1$ into a neighborhood of $z_0$ that is scaled anisotropically,
by flattening the direction of $v_0$.

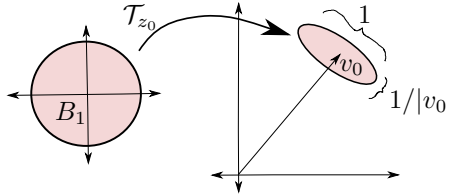
\begin{figure}[ht] \label{f:Tv0}
\setlength{\unitlength}{1in}
 \begin{center}
\begin{picture}(2.052 ,1.000)
\put(0,0){\includegraphics[height=1.000in]{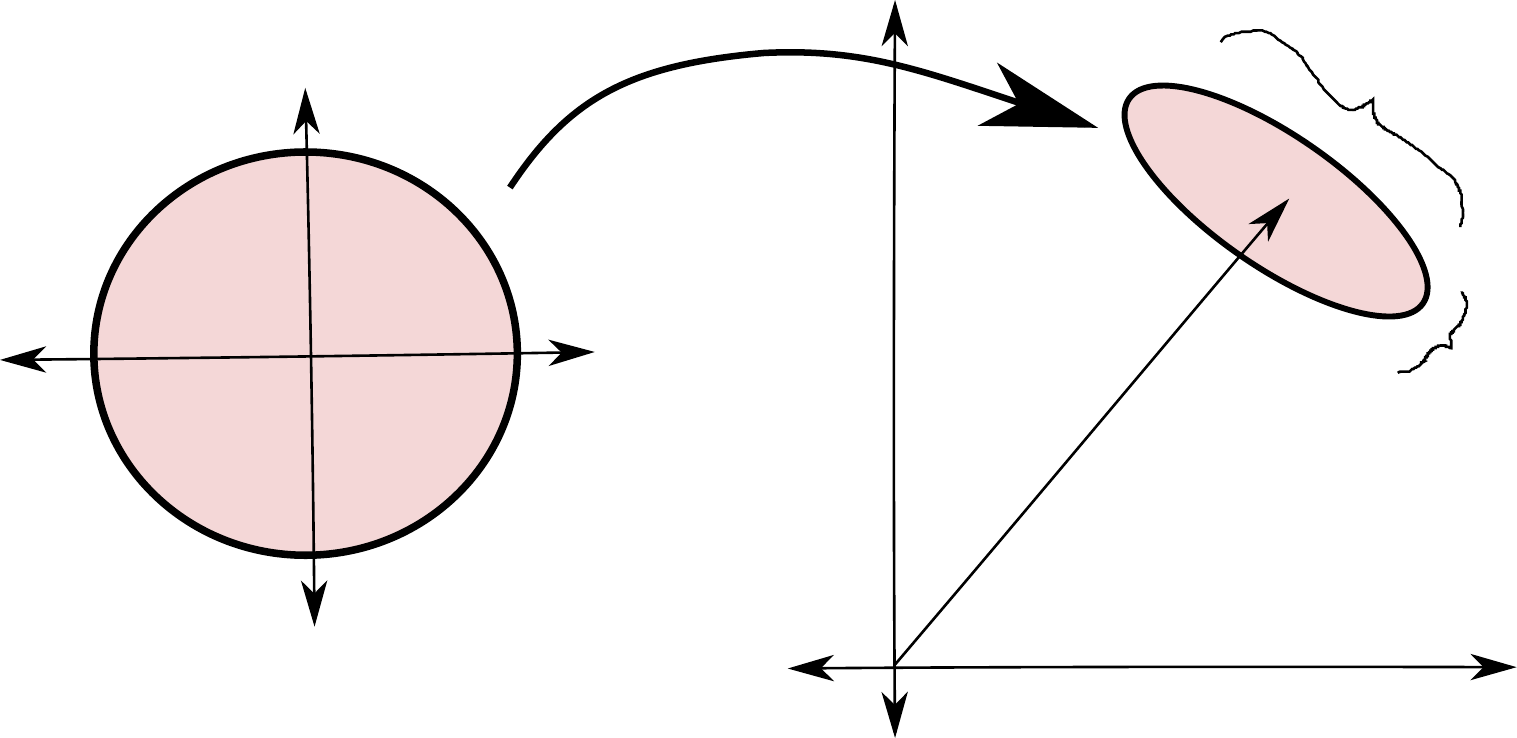}}
\put(1.74,0.64){$v_0$}
\put(1.839,0.893){$1$}
\put(2,0.447){$1 / |v_0|$}
\put(0.606,0.889){$\mathcal T_{z_0}$}
\put(0.252,0.385){$B_1$}
\end{picture}
\caption{The linear transformation $\mathcal{T}_{z_0}$ maps the velocities inside the unit ball $B_1$ to the velocities inside an ellipsoid centered at $v_0$.}
\end{center}
\end{figure}

When $f$ is a solution to the Boltzmann equation, then
$\bar f (z)=f (\mathcal T_{z_0} (z))$ solves a modified
equation in $Q_1$
\[ \partial_t \bar f + v \cdot \nabla_x \bar f - \cL_{\bar K_f} f = \bar h,\]
where
\[ \bar h(z) = c |v_0|^{-\gamma-2s} \bar f(z) (f \ast_v |\cdot|^\gamma)(\mathcal T_{z_0} z),\]
and
\[ \bar K_f(t,x,v,v') = |v_0|^{-\gamma-2s} K(\mathcal T_{z_0} z, v_0 + T_{v_0} v').\]

The benefit of this transformed equation is that the kernel $\bar K$
satisfies the ellipticity conditions \eqref{e:ellipticity_upper},
\eqref{e:ellipticity_cone} and \eqref{e:symmetry_cancellation} in
$Q_1$ with constants depending on the parameters of \eqref{e:hydro}
but independent of the value of $|v_0|$. Thus, Theorems \ref{t:holder}
and \ref{t:schauder} are applied to $\bar f$ in $Q_1$ with fixed
parameters. Then, from the explicit change of variables, we deduce
explicit H\"older estimates that hold globally for $v \in \R^d$.

The formula for the change of variables is motivated by the
nondegeneracy cone described above. Since the set $A$ is concentrated
in a band of width $\approx (1+|v|)^{-1}$ around the equator
perpendicular to $v$, the purpose of the linear transformation
$T_{v_0}$ is of course to stretch this band to make it of width
$\approx 1$. Then, the factor $|v_0|^{-\gamma-2s}$ applied to the time
and space variables normalizes the lower bound
\eqref{e:cone_of_nondeg_cov} by removing the factor that depends on
$(1+|v|)$. The change of variables $\mathcal T_{z_0}$, by design,
transforms the estimate \eqref{e:cone_of_nondeg_cov} into the
condition \eqref{e:ellipticity_cone} for $\bar K$ in $Q_1$, for values
of $\mu$ and $\lambda$ that are uniform with respect to $|v_0|$. It is
remarkable that the same transformation $\mathcal T_{z_0}$ also give
us \eqref{e:ellipticity_upper} and \eqref{e:symmetry_cancellation} for
$\bar K$ in $Q_1$, in terms of $M_0$ and $E_0$ of \eqref{e:hydro},
uniformly in $|v_0|$ as well.

\begin{lemma} \label{l:change_of_variables} Let
  $f:[0,T] \times \R^d \times \R^d$ be a nonnegative function
  satisfying \eqref{e:hydro}. Let $\bar K$ be the kernel described
  above. Then $\bar K$ satisfies \eqref{e:ellipticity_upper},
  \eqref{e:ellipticity_cone} and \eqref{e:symmetry_cancellation} with
  constants depending on the parameters in \eqref{e:hydro}, $s$,
  $\gamma$ and dimension, but \textbf{not} on $z_0$.
\end{lemma}

The proof of Lemma \ref{l:change_of_variables} is given in
\cite{imbert2019regularity}. It involves relatively lengthy technical
computations.

With Lemma \ref{l:change_of_variables} at hand, the H\"older estimate
from Theorem \ref{t:holder}, becomes an explicit global H\"older
estimate for solutions of the Boltzmann equation. The Schauder
estimates from Theorem \ref{t:schauder} give us explicit higher order
estimates for the solutions to the Boltzmann equation with an explicit
(although somewhat complicated) asymptotic behavior for large
velocities. The same logic can be applied to practically any local
regularity estimate. For example combining Proposition
\ref{p:coercivity} with the change of variables of Lemma
\ref{l:change_of_variables}, we recover the sharp coercivity estimate
with respect to the anisotropic distance of Gressman and Strain
\cite{gressman2011}.

\section{The bootstrap argument}

\label{s:bootstrap}

The last step in the proof of Theorem \ref{t:main} is to iteratively
obtain estimates in higher and higher order H\"older norms by applying
the Schauder estimates to increments and derivatives of the solution
$f$ of equation \eqref{e:boltzmann}.

In order to keep track of global H\"older norms that decay as
$|v|\to \infty$, we introduce the following definition.
\[
  [f]_{C^\alpha_{\ell,q}( [\tau,T] \times \R^d \times \R^d)} := \sup
  \left\{ (1+|v|)^q [f]_{\cCl^\alpha(Q_r(z))} : r \in (0,1] \text{ and
    } Q_r(z) \subset [\tau,T] \times \R^d \times \R^d \right\}.\]
The seminorm $[\cdot]_{C^\alpha_{\ell,q}}$ encodes a decay of the form
$(1+|v|)^{-q}$ for the $[\cdot]_{\cCl^\alpha}$ norm as
$|v| \to \infty$.

The upper bound of Theorem \ref{t:upper_bounds} gives us literally an
a priori estimate for
$[f]_{C^0_{\ell,q}( [\tau,T] \times \R^d \times \R^d)}$ for every
value of $q>0$. Next, we apply the H\"older estimates of Theorem
\ref{t:holder} combined with the change of variables of Section
\ref{s:change_of_vars}. We obtain the following conclusion (see
\cite[Proposition 7.1]{imbert2019regularity}).
\begin{proposition} \label{p:global_holder} Let
  $f : [0,T] \times \R^d \times \R^d \to [0,\infty)$ be a solution to
  \eqref{e:boltzmann} so that \eqref{e:hydro} holds. Then, for
  $q > d+\gamma+2s$, we have,
\[ [f]_{C^\alpha_{\ell,q-\gamma_+}} \leq C [f]_{C^0_{\ell,q}( [\tau,T] \times \R^d \times \R^d)}\]
for some $\alpha>0$ and $C$ that depend on the parameters of \eqref{e:hydro}, $d$, $s$ and $\gamma$ only.
\end{proposition}

Proposition \ref{p:global_holder} improves the upper bounds from
Theorem \ref{t:upper_bounds} into a global H\"older estimate.

Since the kernel $K_f$ depends on the function $f$ through the formula
\eqref{e:Kf_exact}, a regularity estimate for $f$ translates into
certain kind of regularity for $K_f$. The H\"older estimate of
Proposition \ref{p:global_holder} implies that Assumption \ref{a:coef}
holds (at least locally) for the kernel $K_f$.

The next step is to improve the smoothness of $f$ by applying the
Schauder estimate, first to $f$ itself, and then to its increments and
derivatives. The following proposition is the result of combining
Theorem~\ref{t:schauder} with the change of variables of Section
\ref{s:change_of_vars}. We state it in terms of a generic function $g$
because we will apply it to $g=f$ and also to $g$ equal to several
combinations of derivatives and increments of $f$. The next
proposition is taken from \cite[Proposition
7.5]{imbert2019regularity}.

\begin{proposition} \label{p:global_schauder}
Let $f : [0,T] \times \R^d \times \R^d \to [0,\infty)$ be a solution to \eqref{e:boltzmann} so that \eqref{e:hydro} holds. Let $g,h :  [0,T] \times \R^d \times \R^d \to [0,\infty)$ so that the following equation holds
\[ \partial_t g + v\cdot \nabla_x g - \cL_{K_f} g = h.\]
Let $\tau>0$, $\alpha \in (0,\min(1,2s))$, and $\alpha'= 2s/(1+2s) \, \alpha$. The following inequality holds
\[ [g]_{C^{2s+\alpha'}_{\ell,q}( [\tau,T]\times \R^d \times \R^d)} \leq C \left( \|g\|_{C^\alpha_{\ell,q+2s+\alpha}( (0,T]\times \R^d \times \R^d)} + \|h\|_{C^{\alpha'}_{\ell,q+2s+\alpha}( (0,T]\times \R^d \times \R^d)} \right),\]
where the constant $C$ depends on the parameters in \eqref{e:hydro}, $d$, $s$, $\gamma$ and $q$.
\end{proposition}

Thus, we can summarize the first three steps of the proof of Theorem \ref{t:main} as the following.
\begin{enumerate}
\item Apply Theorem \ref{t:upper_bounds} to obtain an a priori estimate for $\|f\|_{C^0_{\ell,q}( [\tau,T]\times \R^d \times \R^d)}$, for any values of $\tau>0$ (arbitrarily small) and $q>0$ (arbitrarily large).
\item Apply Proposition \ref{p:global_holder} to obtain an a priori estimate for $\|f\|_{C^\alpha_{\ell,q}( [\tau,T]\times \R^d \times \R^d)}$, for some small $\alpha>0$ and for any values of $\tau>0$ (arbitrarily small) and $q>0$ (arbitrarily large).
\item Apply Proposition \ref{p:global_schauder}, with $g=f$ and $h = (f \ast_v |\cdot|^\gamma) \, f$,  to obtain an a priori estimate for $\|f\|_{C^{2s+\alpha'}_{\ell,q}( [\tau,T]\times \R^d \times \R^d)}$, for an even smaller $\alpha'>0$, and for any values of $\tau>0$ (arbitrarily small) and $q>0$ (arbitrarily large).
\end{enumerate}

The next steps in order to complete the proof of Theorem \ref{t:main} is to apply Proposition \ref{p:global_schauder} with $g$ equal to incremental quotients and derivatives of the solution $f$. To that end, we have to compute the equation satisfied by each of these functions. In each case there are a number of error terms that are absorbed into the right hand side $h$. The precise computations are carried out in \cite{imbert2019regularity}, together with several technical estimates relating the H\"older norms $C^\alpha_{\ell,q}$ with the Boltzmann collision operator $Q$, derivatives and incremental quotients.

\section{Possible generalizations and implications}
\label{s:extensions}

In this section, we review some of the natural implications of Theorem
\ref{t:main} and discuss their possible generalizations.

\subsection{A continuation criteria}

A continuation criteria can easily be derived by combining Theorem
\ref{t:main} with a compatible short time existence result. Indeed,
the short-time existence result says that a smooth solution exists for
some period of time depending on some regularity norm of the initial
data. Theorem \ref{t:main} says that for as long as \eqref{e:hydro}
holds, then no regularity norm of the solution will blow
up. Consequently, the smooth solution can be extended indefinitely.

The existence of smooth solutions for a short period of time for the
non-cutoff Boltzmann equation was first established in
\cite{amuxy2010} for initial data with five derivatives in $L^2_{loc}$
and Gaussian decay (in the $L^2$ sense). This result is not compatible
with our Theorem \ref{t:main} because we cannot ensure the persistence
of the Gaussian decay. Our upper bounds in Theorem
\ref{t:upper_bounds} ensure the persistence of arbitrarily large
algebraic decay rate, but not precisely the Gaussian decay.

The first short time existence result that requires an algebraic decay
rate for the initial data was given in \cite{morimoto2015local} for
$s \in (0,1/2)$ and $\gamma \in (-3/2,0]$. More recently, the result
was extended to the full range of parameters in
\cite{henderson2019local}. These results are compatible with Theorem
\ref{t:main}. Thus, we effectively get the following continuation
condition.

\begin{cor} \label{c:continuation} Assume that
  $f : (0,T) \times \R^d \times \R^d \to [0,\infty)$ is a smooth
  function, periodic in $x$, so that $(1+|v|)^q f(t,x,v)$ is bounded for all $q >0$, and
  $f$ satisfies the Boltzmann equation \eqref{e:boltzmann}. If such a
  solution $f$ cannot be continued to a larger time interval, then one
  of the following events must occur:
\begin{align*}
\lim_{t \to T} \inf_{x \in \R^d} \int_{\R^d} f(t,x,v) \d v &= 0, \\
\lim_{t \to T} \sup_{x \in \R^d} \int_{\R^d} f(t,x,v) \d v &= +\infty, \\
\lim_{t \to T} \sup_{x \in \R^d} \int_{\R^d} |v|^2 f(t,x,v) \d v &= +\infty, \\
\lim_{t \to T} \sup_{x \in \R^d} \int_{\R^d} f \log f(t,x,v) \d v &= +\infty.
\end{align*}
\end{cor}

Our Theorem \ref{t:main} gives us the extra piece of information that the estimate \eqref{e:main_estimate} holds
uniformly for all times (away from zero). Thus, any solution $f$ for
which \eqref{e:hydro} holds is uniformly smooth as $t \to \infty$.

In a very recent preprint \cite{henderson2020self}, Henderson, Snelson and Tarfulea relax the continuation criteria of Corollary~\ref{c:continuation} to only an upper bound for the mass or energy densities. They assume that the initial value $f_0$ is bounded below in a solid ball somewhere (which is always true if $f_0$ is continuous and not identically zero). By a barrier argument, they propagate this lower bound to all positive times and all other points in space. In this way, they prove that a lower bound for the mass density holds for all time, even though it might degenerate as $t \to \infty$. They also use this lower barrier to establish the cone of nondegeneracy for the kernel $K_f$ as in \eqref{e:cone_nondegeneracy}. They observe that the only place where we crucially use the upper bound for the entropy in our program was in establishing this cone of nondegeneracy. In this way, they are able to improve Corollary \ref{c:continuation} and deduce that if a singularity ever occurs at time $T$ for the Boltzmann equation without cutoff, then one of the following two events must occur:
\begin{align*}
\lim_{t \to T} \sup_{x \in \R^d} \int_{\R^d} f(t,x,v) \d v &= +\infty, \\
\lim_{t \to T} \sup_{x \in \R^d} \int_{\R^d} |v|^2 f(t,x,v) \d v &= +\infty.
\end{align*}
Note that using the methods in \cite{henderson2020self}, we cannot obtain global regularity estimates that are uniform for all time, depending only on an upper bound for the mass and energy estimates. Their estimates hold in any interval of time $(t,x,v) \in [0,T] \times \R^d \times \R^d$, but they degenerate as $T  \to \infty$.

\subsection{Convergence to equilibrium}

A very well known result by L. Desvillettes and C. Villani \cite{dv2005} says that
the solution $f$ to the Boltzmann equation \eqref{e:boltzmann}
converges to a Maxwellian provided that the following two conditions
hold.
\begin{enumerate}
\item The function $f$ is $C^\infty$ with uniform bounds as $t \to \infty$.
\item The function $f$ is bounded below by a fixed Maxwellian.
\end{enumerate}

Our Theorem \ref{t:main} says that the first condition of the Theorem
of Desvillettes and Villani holds as soon as \eqref{e:hydro}
holds. The Maxwellian lower bound is obtained also as a consequence
of \eqref{e:hydro} in our joint work with Cl\'ement Mouhot
\cite{imbert2019lower}. Thus, combining all these results,
\eqref{e:hydro} becomes the only condition necessary to deduce the
convergence to equilibrium of the solution $f$ as $t \to \infty$.

\subsection{Weak solutions}

Our Theorem \ref{t:main} is stated as an a priori estimate for smooth
classical solutions. It is natural to wonder if the same regularity
estimates would hold for weak solutions as well. But what is a
\emph{weak solution} exactly?

The term \emph{weak solution} typically means that the function $f$ is
not necessarily smooth, and the equation is to be understood in the
sense of distributions. This notion of solution, in the sense of
distributions, works very well for many linear equations (like the
Laplace, heat or wave equations). However, it fails miserably for most
nonlinear equations. For example, it is well known that solutions in
the sense of distributions have very undesirable properties for
nonlinear conservation law equations, including
non-uniqueness. \emph{Entropy solutions} is the right notion of
solution for scalar conservation laws. Distributional solutions are
also completely unsuitable for the study of the Hamilton-Jacobi
equations, or fully nonlinear parabolic equations. In these latter
cases, the right notion of solution is in the \emph{viscosity sense}.

It is currently very unclear what kind of generalized solution we
should study for the Boltzmann equation \eqref{e:boltzmann}. There is
no notion of weak solution for which we can prove both existence and
uniqueness in fair generality.

The only global existence result, far from equilibrium, known so far
for the non-cutoff Boltzmann equation~\eqref{e:boltzmann} is that of
Alexandre and Villani \cite{av2002}. They prove the existence of a
renormalized solution with defect measure. The uniqueness of solutions
of this kind is not known, and arguably not expected.

It is still an open problem to determine whether the result of Theorem \ref{t:main} holds for the class of solutions defined by Alexandre and Villani. What we mean is that if $f$ is a (not necessarily smooth) renormalized solution with defect measure so that \eqref{e:hydro} holds almost everywhere, then one might expect $f$ to be necessarily smooth and satisfy the estimates in Theorem \ref{t:main}.

Using entropy dissipation estimates from \cite{advw}, we see that for any renormalized solution $f$ of \eqref{e:boltzmann} for which the assumption \eqref{e:hydro} holds almost everywhere, we have $\sqrt{f} \in L^2_{t,x} H_{v,loc}^s$. The main difficulty for moving forward with our regularity program is in establishing the upper bounds of Theorem \ref{t:upper_bounds}. In its proof, we evaluate the equation at a \emph{single point}, where the function $f$ would first cross certain upper barriers.

Except from the upper bounds of Theorem \ref{t:upper_bounds}, the
other steps of the proof are relatively easy to adapt to weak
solutions. We commented in Sections \ref{s:DG} and \ref{s:schauder}
how the H\"older and Schauder estimates would be extended to weak
solutions.

It is debatable whether it is interesting to explore the applicability
of Theorem \ref{t:main} to an intermediate kind of solution. That is,
to a notion of solution stronger than that of Alexandre and Villani,
but weaker than classical smooth solutions. In terms of our proofs,
there seems to be no major obstruction to extend Theorem \ref{t:main}
to hold for solutions in the sense of distributions that are bounded
with a decay of the form $f \leq C (1+|v|)^{q_0}$ for a $q_0$ large
enough as described in \cite[Section 5]{imbert2020decay}. Naturally,
the regularity estimates would not depend on this constant $C$.

\subsection{Mild solutions}

Even though we work with a very strong notion of solution, Theorem \ref{t:main} can be applied for any class of solution where we can establish \emph{stability}. We mean any class of solutions where we can approximate the initial data with a smooth function and pass to the limit. A great example is the global mild solutions to the Boltzmann equation \eqref{e:boltzmann} without cutoff constructed by Duan, Liu, Sakamoto and Strain in \cite{duan2019global}.

These solutions  stay close to a Maxwellian in a suitable norm that does not impose any regularity with respect to the $v$ variable. In this perturbative regime, they establish the global well posedness of the problem. A priori, these solutions could be very rough, but it can be verified by a direct computation that they satisfy our assumption \eqref{e:hydro}. We explain here how to apply Theorem \ref{t:main} to conclude that they are actually $C^\infty$.

If the initial data $f_0$ is smooth and rapidly decaying as $|v| \to \infty$, in addition to satisfying the hypothesis in \cite{duan2019global}, the continuation criteria tells us that it can only blow up when one of the conditions in Corollary~\ref{c:continuation} holds. The smooth solution with initial data $f_0$ must coincide with the solution constructed in \cite{duan2019global} for as long as this smooth solution exists, since uniqueness holds in both regimes. But the solution constructed in \cite{duan2019global} satisfies \eqref{e:hydro} globally, so it can never blow up. Thus, it will be a global smooth solution and the estimates in Theorem \ref{t:main} will hold.

If the initial data $f_0$ satisfies the hypothesis of \cite{duan2019global} but is not smooth, we can approximate it with a sequence $f_0^\varepsilon$ of smooth ones, and rapidly decaying. For each $f_0^\varepsilon$, the solution $f^\varepsilon$ is $C^\infty$ and the result of our Theorem \ref{t:main} holds. Our estimates for $f$ do not depend quantitatively on the smoothness of $f_0^\varepsilon$. Thus, we get uniform regularity estimates that pass to the limit and apply for any (non-smooth) initial data $f_0$. We can legitimately take the limit as $\eps \to 0$ thanks to the well posedness result in \cite{duan2019global}.

We should point out that when $\gamma \leq 0$, the assumptions in \cite{duan2019global} do not imply our decay condition \eqref{e:initial_decay_assumption}. Thus, for $\gamma \leq 0$, we would only apply the result of Theorem \ref{t:main} to the solutions in \cite{duan2019global} if the initial data satisfies \eqref{e:initial_decay_assumption} in addition to their assumptions.

\begin{remark}
It is conceivable that under the strong moment estimates from \cite{duan2019global}, pointwise upper bounds should follow as in Theorem \ref{t:upper_bounds} without any additional hypothesis on $f_0$ even for $\gamma\leq 0$. However, it does not follow directly from any result currently in the literature.
\end{remark}

\section{Open problems}
\label{s:open_problems}

Some open questions were already discussed in the previous section,
including the possible extension of the estimates in Theorem
\ref{t:main} to renormalized solutions with defect measure. In this
last section, we propose several other open problems related to the
material in this survey.

\subsection{Control on the hydrodynamic quantities}

The most obvious open problem after Theorem~\ref{t:main} is whether
the hypothesis in \eqref{e:hydro} can be ensured in any way. This
would imply the unconditional global solvability of the Boltzmann
equation. It is a remarkable open problem. In the first chapter of
\cite{cedric2012theoreme}, C\'edric Villani recounts a lively discussion
with Cl\'ement Mouhot from several years ago about this issue. Based on
the discussion included in the introduction of this survey, we believe
that the global existence of smooth solutions to the non-cutoff
Boltzmann equation is an open problem that, if true, would be harder
to prove than the global solvability of the Navier-Stokes equations
(also if true). The latter is one of the famous Millennium problems.

A more plausible project in the near term would be to remove or weaken
some of the conditions in \eqref{e:hydro}. The well known Prodi-Serrin
condition for the Navier-Stokes equation suggests that maybe only an
upper bound in some suitable $L^p$ space would suffices for the mass,
energy and/or entropy densities. Perhaps only a subset of the
inequalities from \eqref{e:hydro} suffice to obtain Theorem
\ref{t:main}. There are several possibilities. We do not have any
precise conjecture in this direction.

\subsection{Very soft potentials}
\label{s:very_soft}

The non-cutoff collision kernel \eqref{e:B} makes sense for
$s \in (0,1)$ and $\gamma > -d$. Yet, we only present our results in
the range $\gamma+2s \in [0,2]$. The case $\gamma+2s > 2$ is covered
in \cite{cameron2019velocity}. The case $\gamma+2s < 0$ remains open
and would require new ideas.

The case $\gamma+2s < 0$ is commonly referred to as \emph{very soft
  potentials}. In particular, $\gamma=-3$ and $s \to 1$ corresponds to
the Landau-Coulomb equation in three dimensions which is most relevant
for the study of plasma dynamics. Our proof of the upper bounds in
Theorem \ref{t:upper_bounds} fails in this range.

When $\gamma+2s<0$, regularity estimates fail even for the
space-homogeneous Boltzmann equation. We do not currently know any
global $L^\infty$ estimate for the solution $f$ of the non-cutoff
Boltzmann equation when $\gamma+2s<0$, even in the space-homogeneous
case and for initial data $f_0$ in the Schwartz class.

The main difficulty for the very soft potential range is to control
the lower order term in \eqref{e:boltzmann_rewritten}, that is
\[ f(v) \, \left( \int_{\R^d} f(v-w) |w|^\gamma \d w \right),\] with
the diffusion term $\cL_{K_f} f$. Naturally, the more negative
$\gamma$ is, the more singular the lower order term in
\eqref{e:boltzmann} becomes. The cone of nondegeneracy described in
Section \ref{s:cone}, together with the upper bounds on mass, energy
and entropy, allows us to control this lower order term with
$\cL_{K_f} f$ only when $\gamma+2s \geq 0$. It seems that in order to
succeed in proving an $L^\infty$ estimate for $\gamma+2s < 0$, we
would need some further understanding on $\cL_{K_f} f$ sharper than
the information we get from its cone of nondegeneracy.

The problem of $L^\infty$ bounds for very soft potentials is also open
in the context of the space homogenous Landau equation. See for
example
\cite{gualdani2016estimates,silvestre2017landau,golse2019partial}.

\subsection{Bounded domains.}

We stated Theorem \ref{t:main} for periodic boundary conditions in the
space variable $x$. It would be straight forward to extend the result
to other variants of space domains without boundary. For example, we
may assume that as $|x| \to \infty$ the solution $f(t,x,v)$ converges
uniformly to a fixed Maxwellian $M(v)$. Our proof would work
\textit{mutatis mutandis} under this alternative formulation.

The case of domains with boundary is of course of physical
relevance. Several types of boundary conditions are considered in the
literature of kinetic equations: diffuse reflection, specular
reflection and bounce back reflection. An extension of Theorem
\ref{t:main} for solutions $f$ in a bounded domain in space, with any
of these boundary conditions, requires further work. There are several
subtleties involved on the effects of the boundary on the regularity
of $f$. See for example \cite{guo2010decay,kim2011thesis,guo2016bv,guo2017regularity,kim2018specular},
for analysis of boundary effects on solutions of the Boltzmann
equation. No analysis has been made yet concerning the possibility to
extend some form of Theorem \ref{t:main} to any domain with boundary,
for any of the physical boundary conditions.

Another possibility would be to extend the result of Theorem
\ref{t:main} as an interior regularity condition. That is, if a
function $f: [0,T] \times \Omega \times \R^d \to [0,\infty)$ solves
the non-cutoff Boltzmann equation \eqref{e:boltzmann} and satisfies
\eqref{e:hydro}, we would expect the same estimates of Theorem
\ref{t:main} to hold in any subdomain of the form
$[\tau,T] \times K \times \R^d$ for any $\tau > 0$ and $K$ compactly
contained in $\Omega$. This is an open problem in a bounded domain
that is independent of the subtleties involved in the analysis of
physical boundary conditions.

\subsection{H\"older estimates for kinetic equations with diffusion in non-divergence form}

Our theorem~\ref{t:holder} is a kinetic non-local version of the
classical result of De Giorgi, Nash and Moser. There is no regularity
assumption on the kernel $K$. The cancellation condition
\eqref{e:symmetry_cancellation} is a nonlocal form of the divergence
structure of elliptic operators. A natural question is: can we replace
the cancellation condition \eqref{e:symmetry_cancellation} by the
\emph{non-divergence} symmetry condition
\eqref{e:symmetry_non-divergence}?

The Boltzmann kernel $K_f$ naturally satisfies the symmetry condition
\eqref{e:symmetry_non-divergence}. It also satisfies the cancellation
condition \eqref{e:symmetry_cancellation}, but it takes more trouble
to verify it. The reason why we presented Theorem~\ref{t:holder} under
the cancellation condition \eqref{e:symmetry_cancellation} instead of
the symmetry condition \eqref{e:symmetry_non-divergence} is simply
because that is what we are able to prove. It does not mean that the
alternative is false. In fact, we believe it is probably true as well.

The reason why we succeed to prove Theorem \ref{t:holder} under the
cancellation condition \eqref{e:symmetry_cancellation} but not under
the symmetry condition \eqref{e:symmetry_non-divergence} has its roots
in the study of kinetic equations with second order diffusion.

Let us recall the setting of the classical H\"older estimates for
parabolic equations with rough coefficients. The theorem of De Giorgi,
Nash and Moser gives us estimates in H\"older spaces for solutions of
an equation of the form
\[ u_t - \partial_{x_i}\left( a_{ij}(t,x) \partial_{x_j} u \right) =
  0,\] only under the uniform ellipticity assumption
$\lambda I \leq \{a_{ij}\} \leq \Lambda I$. There is also an analogous
result by Krylov and Safonov for parabolic equations in non-divergence
form
\[ u_t - a_{ij}(t,x) \partial_{x_i x_j} u = 0.\] The techniques
involved in the proofs of the theorem of De Giorgi, Nash and Moser are
very different from those involved in the proof of the theorem by Krylov
and Safonov. A kinetic version of these two types of equations would
be
\[
\begin{cases}
\text{(divergence form)} & \partial_t f + v \cdot \nabla_x f - \partial_{v_i} \left( a_{ij}(t,x,v) \partial_{v_j} \right) f = 0, \\
\text{(non-divergence form)} & \partial_t f + v \cdot \nabla_x f - a_{ij}(t,x,v) \partial_{v_i v_j} f = 0.
\end{cases}
\]

In either case, we should assume the uniform ellipticity condition:
$\lambda I \leq \{a_{ij}(t,x,v)\} \leq \Lambda I$. No further
smoothness assumptions should be made on the coefficients. One would
expect that if the equation holds in a kinetic cylinder $Q_1$, then
the following estimate holds
\[ \|f\|_{C_\ell^\alpha(Q_{1/2})} \leq C \|f\|_{C_\ell^0(Q_1)}.\]

Such an estimate is known to be true in the divergence case (see
\cite{wang2009,wang2011,gimv}). However, it is still an open problem
for the non-divergence case.

\subsection{Coercivity estimates for integro-differential operators}

In Proposition \ref{p:coercivity}, we describe a coercivity estimate
for the Boltzmann collision operator. This coercivity condition is
well known and has a long history in the Boltzmann literature. The
upper bound on the entropy in \eqref{e:hydro} can be replaced by a
lower bound on the temperature tensor, or a more general condition
described in \cite{gressman2011}. However, these results do not follow
from the general condition in \cite{chaker2019coercivity} and rely on
the specific structure of the Boltzmann kernel $K_f$ described in
\eqref{e:Kf_exact}.

The result in \cite{chaker2019coercivity} is a coercivity estimate for
general integro-differential operators of the form $\cL_K$ as in
\eqref{e:integro-differential_operator}, not necessarily related to
the Boltzmann equation. The question is to determine simple sufficient
conditions on a kernel $K(v,v')$ so that the following inequality
holds for some constant $c>0$.
\begin{equation} \label{e:question_coercivity}
 \iint_{\R^d \times \R^d} (f(v) - f(v'))^2 K(v,v') \d v' \d v \geq c \|f\|_{\dot H^s}^2.
\end{equation}

It is obvious that \eqref{e:question_coercivity} holds when
$K(v,v') \gtrsim |v'-v|^{-d-2s}$.  Proposition \ref{p:coercivity} says
that \eqref{e:question_coercivity} follows as a consequence of
\eqref{e:ellipticity_cone}. The assumptions in
\cite{chaker2019coercivity} are more general, but there are still some
simple cases where \eqref{e:question_coercivity} holds even though the
assumptions in \cite{chaker2019coercivity} do not apply. The simplest
example, in two dimensions, is to consider
$\cL f = -((-\partial_{11})^s + (-\partial_{22})^s ) f$. This operator
corresponds to a kernel $K$ consisting of a singular measure
concentrated on $v_1-v'_1 = 0$ and $v_2-v'_2=0$. It is easy to verify
that in this case
\[ -\int_{\R^d} (\cL f) f \d v \gtrsim \|f\|_{\dot H^s}^2,\]
however, the kernel in $\cL$ is a singular measure. It is equal to zero almost everywhere.

A general condition that has been suggested is the following: there
exists $\lambda>0$ so that for every $v \in \R^d$, $r>0$ and
$e \in S^{d-1}$ we have
\[ \int_{B_r(v)} [(v'-v) \cdot e ]_+^2 K(v,v') \d v' \geq \lambda
  r^{2-2s}.\] This condition is satisfied by every kernel that is
currently known to satisfy \eqref{e:question_coercivity}. Whether it
actually implies \eqref{e:question_coercivity} is still an open
problem. See \cite{Dyda-Kassmann-2015} and
\cite{Bux-Kassmann-Schulze-2017} for more references concerning this question.

\subsection{Regularity estimates for moderately soft potentials whose initial data does not decay}

In the case of moderately soft potentials ($\gamma \in [-2s,0]$), our estimates in Theorem \ref{t:main} depend on the decay of the initial data $f_0$ though the values of $N_r$. We require \eqref{e:initial_decay_assumption} to hold for every $r \geq 0$.

We do not expect to have any gain of moments in the case of soft potentials. Without the assumption~\eqref{e:initial_decay_assumption}, we only expect a moderate decay for large velocities at positive time, as described in the third case of Theorem \ref{t:upper_bounds}.

It is natural to wonder whether $C^\infty$ regularity estimates may hold in the case $\gamma \leq 0$, for solutions that do not enjoy a fast decay for large velocities. The iteration leading to our proof of Theorem \ref{t:main} does not suffice to prove it. In each step we gain some differentiation at the expense of some decay power. If we do not start by a fast decay as in \eqref{e:initial_decay_assumption}, our iteration would stop after finitely many steps.

\subsection{Conditional propagation of regularity for the cutoff Boltzmann equation}

With the cutoff assumption on the collision kernel $B$, there is no
regularization effect in the Boltzmann equation. However, it is still
plausible to expect propagation of regularity. If the initial data
$f_0$ is smooth and rapidly decaying as $|v| \to \infty$, it is
conceivable that the solution of the Boltzmann equation
\eqref{e:boltzmann} \textbf{with cutoff} stays smooth for as long as
\eqref{e:hydro} holds, at least in the case of hard potentials.

This problem is of a very different nature compared with the
regularization estimates studied in these notes. We have not done any
work on the cutoff case so far.

\bibliographystyle{smfplain}
\bibliography{survey}

\def\cprime{$'$}
\providecommand{\bysame}{\leavevmode ---\ }
\providecommand{\og}{``}
\providecommand{\fg}{''}
\providecommand{\smfandname}{\&}
\providecommand{\smfedsname}{\'eds.}
\providecommand{\smfedname}{\'ed.}
\providecommand{\smfmastersthesisname}{M\'emoire}
\providecommand{\smfphdthesisname}{Th\`ese}
\begin{thebibliography}{100}

\bibitem{advw}
{\scshape R.~Alexandre, L.~Desvillettes, C.~Villani {\normalfont \smfandname}
  B.~Wennberg} -- {\og Entropy dissipation and long-range interactions\fg},
  \emph{Arch. Ration. Mech. Anal.} \textbf{152} (2000), no.~4, p.~327--355.

\bibitem{AlexandreII}
{\scshape R.~Alexandre, Y.~Morimoto, S.~Ukai, C.-J. Xu {\normalfont
  \smfandname} T.~Yang} -- {\og The {B}oltzmann equation without angular cutoff
  in the whole space: {II}, {G}lobal existence for hard potential\fg},
  \emph{Anal. Appl. (Singap.)} \textbf{9} (2011), no.~2, p.~113--134.

\bibitem{amuxy2011ultimate}
\bysame , {\og The {B}oltzmann equation without angular cutoff in the whole
  space: qualitative properties of solutions\fg}, \emph{Arch. Ration. Mech.
  Anal.} \textbf{202} (2011), no.~2, p.~599--661.

\bibitem{AlexandreI}
\bysame , {\og The {B}oltzmann equation without angular cutoff in the whole
  space: {I}, {G}lobal existence for soft potential\fg}, \emph{J. Funct. Anal.}
  \textbf{262} (2012), no.~3, p.~915--1010.

\bibitem{AlexandreSome}
{\scshape R.~Alexandre} -- {\og Some solutions of the {B}oltzmann equation
  without angular cutoff\fg}, \emph{J. Statist. Phys.} \textbf{104} (2001),
  no.~1-2, p.~327--358.

\bibitem{AlexandreReview}
\bysame , {\og A review of {B}oltzmann equation with singular kernels\fg},
  \emph{Kinet. Relat. Models} \textbf{2} (2009), no.~4, p.~551--646.

\bibitem{AlexandreLittlewood1}
{\scshape R.~Alexandre {\normalfont \smfandname} M.~El~Safadi} -- {\og
  Littlewood-{P}aley theory and regularity issues in {B}oltzmann homogeneous
  equations. {I}. {N}on-cutoff case and {M}axwellian molecules\fg}, \emph{Math.
  Models Methods Appl. Sci.} \textbf{15} (2005), no.~6, p.~907--920.

\bibitem{AlexandreLittlewood2}
{\scshape R.~Alexandre {\normalfont \smfandname} M.~Elsafadi} -- {\og
  Littlewood-{P}aley theory and regularity issues in {B}oltzmann homogeneous
  equations. {II}. {N}on cutoff case and non {M}axwellian molecules\fg},
  \emph{Discrete Contin. Dyn. Syst.} \textbf{24} (2009), no.~1, p.~1--11.

\bibitem{amuxy2010}
{\scshape R.~Alexandre, Y.~Morimoto, S.~Ukai, C.-J. Xu {\normalfont
  \smfandname} T.~Yang} -- {\og Regularizing effect and local existence for the
  non-cutoff {B}oltzmann equation\fg}, \emph{Archive for Rational Mechanics and
  Analysis} \textbf{198} (2010), no.~1, p.~39--123.

\bibitem{alexandre2012smoothing}
\bysame , {\og Smoothing effect of weak solutions for the spatially homogeneous
  {B}oltzmann equation without angular cutoff\fg}, \emph{Kyoto J. Math.}
  \textbf{52} (2012), no.~3, p.~433--463.

\bibitem{av2002}
{\scshape R.~Alexandre {\normalfont \smfandname} C.~Villani} -- {\og On the
  {B}oltzmann equation for long-range interactions\fg}, \emph{Communications on
  pure and applied mathematics} \textbf{55} (2002), no.~1, p.~30--70.

\bibitem{arsenio2012}
{\scshape D.~Ars\'{e}nio {\normalfont \smfandname} N.~Masmoudi} -- {\og
  Regularity of renormalized solutions in the {B}oltzmann equation with
  long-range interactions\fg}, \emph{Comm. Pure Appl. Math.} \textbf{65}
  (2012), no.~4, p.~508--548.

\bibitem{bardos1991}
{\scshape C.~Bardos, F.~Golse {\normalfont \smfandname} D.~Levermore} -- {\og
  Fluid dynamic limits of kinetic equations. {I}. {F}ormal derivations\fg},
  \emph{J. Statist. Phys.} \textbf{63} (1991), no.~1-2, p.~323--344.

\bibitem{barlow2009non}
{\scshape M.~T. Barlow, R.~F. Bass, Z.-Q. Chen {\normalfont \smfandname}
  M.~Kassmann} -- {\og Non-local {D}irichlet forms and symmetric jump
  processes\fg}, \emph{Trans. Amer. Math. Soc.} \textbf{361} (2009), no.~4,
  p.~1963--1999.

\bibitem{bass2005harnack}
{\scshape R.~F. Bass {\normalfont \smfandname} M.~Kassmann} -- {\og Harnack
  inequalities for non-local operators of variable order\fg},
  \emph{Transactions of the American Mathematical Society} (2005), p.~837--850.

\bibitem{bass2005holder}
\bysame , {\og H{\"o}lder continuity of harmonic functions with respect to
  operators of variable order\fg}, \emph{Communications in Partial Differential
  Equations} \textbf{30} (2005), no.~8, p.~1249--1259.

\bibitem{bass2002harnack}
{\scshape R.~F. Bass {\normalfont \smfandname} D.~A. Levin} -- {\og Harnack
  inequalities for jump processes\fg}, \emph{Potential Analysis} \textbf{17}
  (2002), no.~4, p.~375--388.

\bibitem{basslevin2002}
{\scshape R.~F. Bass {\normalfont \smfandname} D.~A. Levin} -- {\og Transition
  probabilities for symmetric jump processes\fg}, \emph{Trans. Amer. Math.
  Soc.} \textbf{354} (2002), no.~7, p.~2933--2953 (electronic).

\bibitem{bjorland2012}
{\scshape C.~Bjorland, L.~Caffarelli {\normalfont \smfandname} A.~Figalli} --
  {\og Non-local gradient dependent operators\fg}, \emph{Adv. Math.}
  \textbf{230} (2012), no.~4-6, p.~1859--1894.

\bibitem{buckmaster2019formation2}
{\scshape T.~Buckmaster, S.~Shkoller {\normalfont \smfandname} V.~Vicol} --
  {\og {Formation of point shocks for 3D compressible Euler}\fg},
  arXiv:1912.04429, 2019.

\bibitem{buckmaster2019formation}
\bysame , {\og {Formation of shocks for 2D isentropic compressible Euler}\fg},
  arXiv:1907.03784, 2019.

\bibitem{Bux-Kassmann-Schulze-2017}
{\scshape K.~uwe {Bux}, M.~{Kassmann} {\normalfont \smfandname} T.~{Schulze}}
  -- {\og {Quadratic forms and Sobolev spaces of fractional order}\fg},
  \emph{{Proc. Lond. Math. Soc. (3)}} \textbf{119} (2019), no.~3, p.~841--866
  (English).

\bibitem{chan2011}
{\scshape L.~Caffarelli, C.~H. Chan {\normalfont \smfandname} A.~Vasseur} --
  {\og Regularity theory for parabolic nonlinear integral operators\fg},
  \emph{J. Amer. Math. Soc.} \textbf{24} (2011), no.~3, p.~849--869.

\bibitem{caffarelli2009regularity}
{\scshape L.~Caffarelli {\normalfont \smfandname} L.~Silvestre} -- {\og
  Regularity theory for fully nonlinear integro-differential equations\fg},
  \emph{Communications on Pure and Applied Mathematics} \textbf{62} (2009),
  no.~5, p.~597--638.

\bibitem{caffarelli2010drift}
{\scshape L.~Caffarelli {\normalfont \smfandname} A.~Vasseur} -- {\og Drift
  diffusion equations with fractional diffusion and the quasi-geostrophic
  equation\fg}, \emph{Ann. of Math} \textbf{171} (2010), no.~3, p.~1903--1930.

\bibitem{cameron2018}
{\scshape S.~Cameron, L.~Silvestre {\normalfont \smfandname} S.~Snelson} --
  {\og Global a priori estimates for the inhomogeneous {L}andau equation with
  moderately soft potentials\fg}, \emph{Ann. Inst. H. Poincar\'{e} Anal. Non
  Lin\'{e}aire} \textbf{35} (2018), no.~3, p.~625--642.

\bibitem{cameron2019velocity}
{\scshape S.~Cameron {\normalfont \smfandname} S.~Snelson} -- {\og {Velocity
  Decay Estimates for Boltzmann equation with hard potentials}\fg},
  arXiv:1910.01613, 2019.

\bibitem{MR1555365}
{\scshape T.~Carleman} -- {\og Sur la th\'{e}orie de l'\'{e}quation
  int\'{e}grodiff\'{e}rentielle de {B}oltzmann\fg}, \emph{Acta Math.}
  \textbf{60} (1933), no.~1, p.~91--146.

\bibitem{MR0098477}
\bysame , \emph{Probl\`emes math\'{e}matiques dans la th\'{e}orie cin\'{e}tique
  des gaz}, Publ. Sci. Inst. Mittag-Leffler. 2, Almqvist \& Wiksells
  Boktryckeri Ab, Uppsala, 1957.

\bibitem{chaker2019coercivity}
{\scshape J.~Chaker {\normalfont \smfandname} L.~Silvestre} -- {\og Coercivity
  estimates for integro-differential operators\fg}, arXiv:1904.13014, 2019.

\bibitem{davila2012nonsymmetric}
{\scshape H.~Chang~Lara {\normalfont \smfandname} G.~D{\'a}vila} -- {\og
  Regularity for solutions of nonlocal, nonsymmetric equations\fg}, \emph{Ann.
  Inst. H. Poincar\'e Anal. Non Lin\'eaire} \textbf{29} (2012), no.~6,
  p.~833--859.

\bibitem{ChenHe2011}
{\scshape Y.~Chen {\normalfont \smfandname} L.~He} -- {\og Smoothing estimates
  for {B}oltzmann equation with full-range interactions: spatially homogeneous
  case\fg}, \emph{Arch. Ration. Mech. Anal.} \textbf{201} (2011), no.~2,
  p.~501--548.

\bibitem{ChenHe2012}
\bysame , {\og Smoothing estimates for {B}oltzmann equation with full-range
  interactions: {S}patially inhomogeneous case\fg}, \emph{Arch. Ration. Mech.
  Anal.} \textbf{203} (2012), no.~2, p.~343--377.

\bibitem{christodoulou2007}
{\scshape D.~Christodoulou} -- \emph{The formation of shocks in 3-dimensional
  fluids}, EMS Monographs in Mathematics, European Mathematical Society (EMS),
  Z\"{u}rich, 2007.

\bibitem{constantin2012}
{\scshape P.~Constantin {\normalfont \smfandname} V.~Vicol} -- {\og Nonlinear
  maximum principles for dissipative linear nonlocal operators and
  applications\fg}, \emph{Geom. Funct. Anal.} \textbf{22} (2012), no.~5,
  p.~1289--1321.

\bibitem{dv2005}
{\scshape L.~Desvillettes {\normalfont \smfandname} C.~Villani} -- {\og On the
  trend to global equilibrium for spatially inhomogeneous kinetic systems: the
  {B}oltzmann equation\fg}, \emph{Invent. Math.} \textbf{159} (2005), no.~2,
  p.~245--316.

\bibitem{desvillettes2009stability}
{\scshape L.~Desvillettes {\normalfont \smfandname} C.~Mouhot} -- {\og
  Stability and uniqueness for the spatially homogeneous {B}oltzmann equation
  with long-range interactions\fg}, \emph{Arch. Ration. Mech. Anal.}
  \textbf{193} (2009), no.~2, p.~227--253.

\bibitem{desvillettes2005smoothness}
{\scshape L.~Desvillettes {\normalfont \smfandname} B.~Wennberg} -- {\og
  Smoothness of the solution of the spatially homogeneous {B}oltzmann equation
  without cutoff\fg}, \emph{Comm. Partial Differential Equations} \textbf{29}
  (2004), no.~1-2, p.~133--155.

\bibitem{dfp}
{\scshape M.~Di~Francesco {\normalfont \smfandname} S.~Polidoro} -- {\og
  Schauder estimates, {H}arnack inequality and {G}aussian lower bound for
  {K}olmogorov-type operators in non-divergence form\fg}, \emph{Adv.
  Differential Equations} \textbf{11} (2006), no.~11, p.~1261--1320.

\bibitem{tj2018}
{\scshape H.~Dong, T.~Jin {\normalfont \smfandname} H.~Zhang} -- {\og Dini and
  {S}chauder estimates for nonlocal fully nonlinear parabolic equations with
  drifts\fg}, \emph{Anal. PDE} \textbf{11} (2018), no.~6, p.~1487--1534.

\bibitem{duan2019global}
{\scshape R.~Duan, S.~Liu, S.~Sakamoto {\normalfont \smfandname} R.~M. Strain}
  -- {\og {Global mild solutions of the Landau and non-cutoff Boltzmann
  equations}\fg}, arXiv:1904.12086, 2019.

\bibitem{Dyda-Kassmann-2015}
{\scshape B.~{Dyda} {\normalfont \smfandname} M.~{Kassmann}} -- {\og
  {Regularity estimates for elliptic nonlocal operators}\fg}, \emph{{Anal.
  PDE}} \textbf{13} (2020), no.~2, p.~317--370 (English).

\bibitem{eidelman1998}
{\scshape S.~D. Eidelman, S.~D. Ivasyshen {\normalfont \smfandname} H.~P.
  Malytska} -- {\og A modified {L}evi method: development and application\fg},
  \emph{Dopov. Nats. Akad. Nauk Ukr. Mat. Prirodozn. Tekh. Nauki} (1998),
  no.~5, p.~14--19.

\bibitem{Felsinger2013}
{\scshape M.~Felsinger {\normalfont \smfandname} M.~Kassmann} -- {\og Local
  {R}egularity for {P}arabolic {N}onlocal {O}perators\fg}, \emph{Comm. Partial
  Differential Equations} \textbf{38} (2013), no.~9, p.~1539--1573.

\bibitem{gimv}
{\scshape F.~Golse, C.~Imbert, C.~Mouhot {\normalfont \smfandname} A.~F.
  Vasseur} -- {\og Harnack inequality for kinetic {F}okker-{P}lanck equations
  with rough coefficients and application to the {L}andau equation\fg},
  \emph{Ann. Sc. Norm. Super. Pisa Cl. Sci. (5)} \textbf{19} (2019), no.~1,
  p.~253--295.

\bibitem{golse2019partial}
{\scshape F.~Golse, M.~P. Gualdani, C.~Imbert {\normalfont \smfandname}
  A.~Vasseur} -- {\og {Partial regularity in time for the space homogeneous
  Landau equation with Coulomb potential}\fg}, arXiv:1906.02841, 2019.

\bibitem{gs2011}
{\scshape P.~T. Gressman {\normalfont \smfandname} R.~M. Strain} -- {\og Global
  classical solutions of the {B}oltzmann equation without angular cut-off\fg},
  \emph{J. Amer. Math. Soc.} \textbf{24} (2011), no.~3, p.~771--847.

\bibitem{gressman2011}
\bysame , {\og Sharp anisotropic estimates for the {B}oltzmann collision
  operator and its entropy production\fg}, \emph{Adv. Math.} \textbf{227}
  (2011), no.~6, p.~2349--2384.

\bibitem{gualdani2016estimates}
{\scshape M.~Gualdani {\normalfont \smfandname} N.~Guillen} -- {\og {Estimates
  for radial solutions of the homogeneous Landau equation with Coulomb
  potential}\fg}, \emph{Analysis \& PDE} \textbf{9} (2016), no.~8,
  p.~1773--1810.

\bibitem{guo2016bv}
{\scshape Y.~Guo, C.~Kim, D.~Tonon {\normalfont \smfandname} A.~Trescases} --
  {\og B{V}-regularity of the {B}oltzmann equation in non-convex domains\fg},
  \emph{Arch. Ration. Mech. Anal.} \textbf{220} (2016), no.~3, p.~1045--1093.

\bibitem{guo2010decay}
{\scshape Y.~Guo} -- {\og Decay and continuity of the {B}oltzmann equation in
  bounded domains\fg}, \emph{Arch. Ration. Mech. Anal.} \textbf{197} (2010),
  no.~3, p.~713--809.

\bibitem{guo2017regularity}
{\scshape Y.~Guo, C.~Kim, D.~Tonon {\normalfont \smfandname} A.~Trescases} --
  {\og Regularity of the {B}oltzmann equation in convex domains\fg},
  \emph{Invent. Math.} \textbf{207} (2017), no.~1, p.~115--290.

\bibitem{he2016sharp}
{\scshape L.-B. He} -- {\og Sharp bounds for {B}oltzmann and {L}andau collision
  operators\fg}, \emph{Ann. Sci. \'{E}c. Norm. Sup\'{e}r. (4)} \textbf{51}
  (2018), no.~5, p.~1253--1341.

\bibitem{he2012wellposedness}
{\scshape L.~He} -- {\og Well-posedness of spatially homogeneous {B}oltzmann
  equation with full-range interaction\fg}, \emph{Comm. Math. Phys.}
  \textbf{312} (2012), no.~2, p.~447--476.

\bibitem{henderson2019local}
{\scshape C.~Henderson, S.~Snelson {\normalfont \smfandname} A.~Tarfulea} --
  {\og {Local well-posedness of the Boltzmann equation with polynomially
  decaying initial data}\fg}, arXiv:1910.07138, 2019.

\bibitem{henderson2020self}
\bysame , {\og {Self-generating lower bounds and continuation for the Boltzmann
  equation}\fg}, arXiv:2005.13668, 2020.

\bibitem{huo2008regularity}
{\scshape Z.~Huo, Y.~Morimoto, S.~Ukai {\normalfont \smfandname} T.~Yang} --
  {\og Regularity of solutions for spatially homogeneous {B}oltzmann equation
  without angular cutoff\fg}, \emph{Kinet. Relat. Models} \textbf{1} (2008),
  no.~3, p.~453--489.

\bibitem{imbert2016schauder}
{\scshape C.~Imbert, T.~Jin {\normalfont \smfandname} R.~Shvydkoy} -- {\og
  Schauder estimates for an integro-differential equation with applications to
  a nonlocal {B}urgers equation\fg}, \emph{Ann. Fac. Sci. Toulouse Math. (6)}
  \textbf{27} (2018), no.~4, p.~667--677.

\bibitem{imbert2019lower}
{\scshape C.~Imbert, C.~Mouhot {\normalfont \smfandname} L.~Silvestre} -- {\og
  {Gaussian lower bounds for the Boltzmann equation without cut-off}\fg},
  HAL-02078069, Mar 2019.

\bibitem{imbert2020decay}
\bysame , {\og {Decay estimates for large velocities in the Boltzmann equation
  without cutoff}\fg}, \emph{{Journal de l'{\'E}cole polytechnique -
  Math{\'e}matiques}} \textbf{7} (2020), p.~143--184.

\bibitem{imbert2019regularity}
{\scshape C.~Imbert {\normalfont \smfandname} L.~Silvestre} -- {\og {Global
  regularity estimates for the Boltzmann equation without cut-off}\fg},
  arXiv:1909.12729, 2019.

\bibitem{imbert2020schauder}
\bysame , {\og {The Schauder estimate for kinetic integral equations}\fg},
  \emph{{Analysis and PDEs}} (2020), {To appear}.

\bibitem{imbert2020weak}
\bysame , {\og {The weak Harnack inequality for the Boltzmann equation without
  cut-off}\fg}, \emph{{Journal of the European Mathematical Society}}
  \textbf{22} (2020), no.~2, p.~pp. 507--592.

\bibitem{tj2015}
{\scshape T.~Jin {\normalfont \smfandname} J.~Xiong} -- {\og Schauder estimates
  for solutions of linear parabolic integro-differential equations\fg},
  \emph{Discrete Contin. Dyn. Syst.} \textbf{35} (2015), no.~12, p.~5977--5998.

\bibitem{kassmann2009priori}
{\scshape M.~Kassmann} -- {\og A priori estimates for integro-differential
  operators with measurable kernels\fg}, \emph{Calculus of Variations and
  Partial Differential Equations} \textbf{34} (2009), no.~1, p.~1--21.

\bibitem{kassmann2013intrinsic}
{\scshape M.~{Kassmann} {\normalfont \smfandname} A.~{Mimica}} -- {\og
  {Intrinsic scaling properties for nonlocal operators}\fg}, \emph{{J. Eur.
  Math. Soc. (JEMS)}} \textbf{19} (2017), no.~4, p.~983--1011 (English).

\bibitem{rang2013h}
{\scshape M.~Kassmann, M.~Rang {\normalfont \smfandname} R.~W. Schwab} -- {\og
  Integro-differential equations with nonlinear directional dependence\fg},
  \emph{Indiana Univ. Math. J.} \textbf{63} (2014), no.~5, p.~1467--1498.

\bibitem{kassmann2013regularity}
{\scshape M.~Kassmann {\normalfont \smfandname} R.~W. Schwab} -- {\og
  Regularity results for nonlocal parabolic equations\fg}, \emph{Riv. Math.
  Univ. Parma (N.S.)} \textbf{5} (2014), no.~1, p.~183--212.

\bibitem{kim2011thesis}
{\scshape C.~Kim} -- {\og Formation and propagation of discontinuity for
  {B}oltzmann equation in non-convex domains\fg}, \emph{Comm. Math. Phys.}
  \textbf{308} (2011), no.~3, p.~641--701.

\bibitem{kim2018specular}
{\scshape C.~Kim {\normalfont \smfandname} D.~Lee} -- {\og The {B}oltzmann
  equation with specular boundary condition in convex domains\fg}, \emph{Comm.
  Pure Appl. Math.} \textbf{71} (2018), no.~3, p.~411--504.

\bibitem{kolm}
{\scshape A.~Kolmogoroff} -- {\og Zuf\"allige {B}ewegungen (zur {T}heorie der
  {B}rownschen {B}ewegung)\fg}, \emph{Ann. of Math. (2)} \textbf{35} (1934),
  no.~1, p.~116--117.

\bibitem{komatsu1995}
{\scshape T.~Komatsu} -- {\og Uniform estimates for fundamental solutions
  associated with non-local {D}irichlet forms\fg}, \emph{Osaka J. Math.}
  \textbf{32} (1995), no.~4, p.~833--860.

\bibitem{ks}
{\scshape N.~V. Krylov {\normalfont \smfandname} M.~V. Safonov} -- {\og A
  property of the solutions of parabolic equations with measurable
  coefficients\fg}, \emph{Izv. Akad. Nauk SSSR Ser. Mat.} \textbf{44} (1980),
  no.~1, p.~161--175, 239.

\bibitem{zbMATH03277871}
{\scshape O.~A. {Ladyzhenskaya}, V.~A. {Solonnikov} {\normalfont \smfandname}
  N.~N. {Ural'tseva}} -- \emph{{Linear and quasi-linear equations of parabolic
  type. Translated from the Russian by S. Smith.}}, vol.~23, American
  Mathematical Society (AMS), Providence, RI, 1968 (English).

\bibitem{davila2014parabolic}
{\scshape H.~C. Lara {\normalfont \smfandname} G.~D{\'a}vila} -- {\og
  Regularity for solutions of non local parabolic equations\fg}, \emph{Calc.
  Var. Partial Differential Equations} \textbf{49} (2014), no.~1-2,
  p.~139--172.

\bibitem{Luk2018}
{\scshape J.~Luk {\normalfont \smfandname} J.~Speck} -- {\og Shock formation in
  solutions to the 2{D} compressible {E}uler equations in the presence of
  non-zero vorticity\fg}, \emph{Invent. Math.} \textbf{214} (2018), no.~1,
  p.~1--169.

\bibitem{lunardi1997}
{\scshape A.~Lunardi} -- {\og Schauder estimates for a class of degenerate
  elliptic and parabolic operators with unbounded coefficients in {$\R^n$}\fg},
  \emph{Ann. Scuola Norm. Sup. Pisa Cl. Sci. (4)} \textbf{24} (1997), no.~1,
  p.~133--164.

\bibitem{manfredini}
{\scshape M.~Manfredini} -- {\og The {D}irichlet problem for a class of
  ultraparabolic equations\fg}, \emph{Adv. Differential Equations} \textbf{2}
  (1997), no.~5, p.~831--866.

\bibitem{merle2019smooth}
{\scshape F.~Merle, P.~Raphael, I.~Rodnianski {\normalfont \smfandname}
  J.~Szeftel} -- {\og {On smooth self similar solutions to the compressible
  Euler equations}\fg}, arXiv:1912.10998, 2019.

\bibitem{merle2019implosion}
\bysame , {\og On the implosion of a three dimensional compressible fluid\fg},
  arXiv:1912.11009, 2019.

\bibitem{MP}
{\scshape R.~Mikulevicius {\normalfont \smfandname} H.~Pragarauskas} -- {\og On
  the {C}auchy problem for integro-differential operators in {H}{\"o}lder
  classes and the uniqueness of the martingale problem\fg}, \emph{Potential
  Anal.} \textbf{40} (2014), no.~4, p.~539--563.

\bibitem{morimoto2010gevrey}
{\scshape Y.~Morimoto {\normalfont \smfandname} S.~Ukai} -- {\og Gevrey
  smoothing effect of solutions for spatially homogeneous nonlinear {B}oltzmann
  equation without angular cutoff\fg}, \emph{J. Pseudo-Differ. Oper. Appl.}
  \textbf{1} (2010), no.~1, p.~139--159.

\bibitem{morimoto2009regularity}
{\scshape Y.~Morimoto, S.~Ukai, C.-J. Xu {\normalfont \smfandname} T.~Yang} --
  {\og {Regularity of solutions to the spatially homogeneous Boltzmann equation
  without angular cutoff}\fg}, \emph{Discrete and Continuous Dynamical
  Systems-Series A} \textbf{24} (2009), no.~1, p.~187--212.

\bibitem{morimoto2015local}
{\scshape Y.~Morimoto {\normalfont \smfandname} T.~Yang} -- {\og {Local
  existence of polynomial decay solutions to the Boltzmann equation for soft
  potentials}\fg}, \emph{Analysis and Applications} \textbf{13} (2015), no.~06,
  p.~663--683.

\bibitem{MouhotExpl}
{\scshape C.~Mouhot} -- {\og Explicit coercivity estimates for the linearized
  {B}oltzmann and {L}andau operators\fg}, \emph{Comm. Partial Differential
  Equations} \textbf{31} (2006), no.~7-9, p.~1321--1348.

\bibitem{pascucci2004}
{\scshape A.~Pascucci {\normalfont \smfandname} S.~Polidoro} -- {\og The
  {M}oser's iterative method for a class of ultraparabolic equations\fg},
  \emph{Commun. Contemp. Math.} \textbf{6} (2004), no.~3, p.~395--417.

\bibitem{radkevich2008}
{\scshape E.~V. Radkevich} -- {\og Equations with nonnegative characteristic
  form. {II}\fg}, \emph{Sovrem. Mat. Prilozh.} (2008), no.~56,
  Differentsial\cprime nye Uravneniya s Chastnymi Proizvodnymi, p.~3--147.

\bibitem{schwab2016}
{\scshape R.~W. Schwab {\normalfont \smfandname} L.~Silvestre} -- {\og
  Regularity for parabolic integro-differential equations with very irregular
  kernels\fg}, \emph{Anal. PDE} \textbf{9} (2016), no.~3, p.~727--772.

\bibitem{sergio2004recent}
{\scshape P.~Sergio} -- {\og Recent results on {K}olmogorov equations and
  applications\fg}, in \emph{Workshop on Second Order Subelliptic Equations and
  Applications}, vol.~3, GRAFICOM Edizioni, 2004, p.~129--143.

\bibitem{serra2015}
{\scshape J.~Serra} -- {\og Regularity for fully nonlinear nonlocal parabolic
  equations with rough kernels\fg}, \emph{Calc. Var. Partial Differential
  Equations} \textbf{54} (2015), no.~1, p.~615--629.

\bibitem{sideris1985}
{\scshape T.~C. Sideris} -- {\og Formation of singularities in
  three-dimensional compressible fluids\fg}, \emph{Comm. Math. Phys.}
  \textbf{101} (1985), no.~4, p.~475--485.

\bibitem{silvestre2006holder}
{\scshape L.~Silvestre} -- {\og {Holder estimates for solutions of
  integro-differential equations like the fractional Laplace}\fg},
  \emph{Indiana University mathematics journal} \textbf{55} (2006), no.~3,
  p.~1155--1174.

\bibitem{silvestre2011differentiability}
\bysame , {\og {On the differentiability of the solution to the
  Hamilton--Jacobi equation with critical fractional diffusion}\fg},
  \emph{Advances in mathematics} \textbf{226} (2011), no.~2, p.~2020--2039.

\bibitem{silvestre2016new}
\bysame , {\og {A new regularization mechanism for the Boltzmann equation
  without cut-off}\fg}, \emph{Communications in Mathematical Physics}
  \textbf{348} (2016), no.~1, p.~69--100.

\bibitem{silvestre2017landau}
\bysame , {\og Upper bounds for parabolic equations and the {L}andau
  equation\fg}, \emph{J. Differential Equations} \textbf{262} (2017), no.~3,
  p.~3034--3055.

\bibitem{song2004}
{\scshape R.~Song {\normalfont \smfandname} Z.~Vondra{\v{c}}ek} -- {\og Harnack
  inequality for some classes of {M}arkov processes\fg}, \emph{Math. Z.}
  \textbf{246} (2004), no.~1-2, p.~177--202.

\bibitem{stokols2019}
{\scshape L.~F. Stokols} -- {\og H\"{o}lder continuity for a family of nonlocal
  hypoelliptic kinetic equations\fg}, \emph{SIAM J. Math. Anal.} \textbf{51}
  (2019), no.~6, p.~4815--4847.

\bibitem{cedric2012theoreme}
{\scshape C.~Villani} -- {\og Th{\'e}or{\`e}me vivant\fg}, 2012.

\bibitem{sayro1971}
{\scshape J.~I. \v{S}atyro} -- {\og The smoothness of the solutions of certain
  degenerate second order equations\fg}, \emph{Mat. Zametki} \textbf{10}
  (1971), p.~101--111.

\bibitem{wang2009}
{\scshape W.~Wang {\normalfont \smfandname} L.~Zhang} -- {\og The {$C^\alpha$}
  regularity of a class of non-homogeneous ultraparabolic equations\fg},
  \emph{Sci. China Ser. A} \textbf{52} (2009), no.~8, p.~1589--1606.

\bibitem{wang2011}
{\scshape W.~Wang {\normalfont \smfandname} L.~Zhang} -- {\og The {$C^\alpha$}
  regularity of weak solutions of ultraparabolic equations\fg}, \emph{Discrete
  Contin. Dyn. Syst.} \textbf{29} (2011), no.~3, p.~1261--1275.

\bibitem{zhang2012gevrey}
{\scshape T.-F. Zhang {\normalfont \smfandname} Z.~Yin} -- {\og Gevrey
  regularity of spatially homogeneous {B}oltzmann equation without cutoff\fg},
  \emph{J. Differential Equations} \textbf{253} (2012), no.~4, p.~1172--1190.

\end{thebibliography}
\end{document}